\documentclass[11pt]{article}
\usepackage[a4paper, margin=1in, footskip=0.25in]{geometry}
\usepackage{graphicx}
\usepackage{amsmath,mathrsfs,mathtools}
\usepackage{float}
\usepackage{amssymb}
\usepackage{amsthm}
\usepackage[toc,page]{appendix}
\usepackage[mathscr]{euscript}
\usepackage{upgreek}
\usepackage{dsfont}
\usepackage[authoryear]{natbib}
\usepackage[T1]{fontenc}
\usepackage{newtxtext,newtxmath} % Modern Times Roman alternative; tighter and cleaner than Computer Modern
\usepackage[patch=none]{microtype} % Improves character protrusion and font expansion (reduces hyphenations and saves space)
\usepackage{url}
\usepackage[dvipsnames]{xcolor}
\usepackage[colorlinks=true, linkcolor=black, citecolor=OliveGreen, urlcolor=black]{hyperref}
\usepackage{fontawesome}
\usepackage{epsdice}
\usepackage{booktabs}
\usepackage{subcaption}

\usepackage{setspace}
\allowdisplaybreaks % Allows page breaks inside long align environments
\theoremstyle{plain}
\newtheorem{theorem}{Theorem}[section]
\newtheorem{lemma}[theorem]{Lemma}
\newtheorem{proposition}[theorem]{Proposition}
\newtheorem{corollary}[theorem]{Corollary}

\theoremstyle{definition}
\newtheorem{definition}[theorem]{Definition}

\newtheorem{assumption}[theorem]{Assumption} % Critical for financial math

\theoremstyle{remark}
\newtheorem{remark}[theorem]{Remark}

\numberwithin{equation}{section}

\newcommand{\Mmark}{\mathrm M}

\newcommand{\dif}{\mathrm{d}} % Upright 'd' for derivatives/stochastic integrals e.g., \D W_t

\newcommand{\us}{s_1}
\newcommand{\bs}{s_2}

\usepackage{authblk} % Professional author/affiliation block

\title{\textbf{Optimal Dividends and Voluntary Liquidation under Omega-Type Regulatory Termination}}

\author[1]{Ran Xu\thanks{Email: \texttt{ran.xu@xjtlu.edu.cn}}}
\author[2]{Bo Li\thanks{Email: \texttt{libo@nankai.edu.cn}}}
\author[3]{Xin Xu\thanks{Email: \texttt{Xin.Xu03@xjtlu.edu.cn}}}
\author[4]{Zbigniew Palmowski \thanks{Email: \texttt{zbigniew.palmowski@pwr.edu.pl}}}
\affil[1]{Department of Financial and Actuarial Mathematics, Xi'an Jiaotong-Liverpool University, Suzhou,  China.}
\affil[2]{Department of Mathematics and LPMC, Nankai University, Tianjin, China.}
\affil[3]{Department of Foundational Mathematics, Xi'an Jiaotong-Liverpool University, Suzhou,  China.}
\affil[4]{Department of Applied Mathematics, Wrocław University of Science and Technology,
Wybrzeże Wyspiańskiego 27 50-370 Wrocław, Poland}
\date{October 29, 2026}

\begin{document}
\maketitle

%\doublespacing

\begin{abstract}
We study a joint optimal dividend and voluntary liquidation problem for a company whose surplus follows a spectrally negative L\'evy process and is subject to Omega-type regulatory termination. Dividend payments are modeled as a singular control, while voluntary liquidation is an optimal stopping decision with a state-dependent liquidation payoff. When the surplus falls below a regulatory solvency threshold, an Omega clock accumulates the time spent in distress and may terminate the firm before voluntary liquidation occurs.
Under dynamic programming and continuity assumptions for an auxiliary Markov subproblem, we characterize its value as the smallest controlled stochastic supersolution of the associated quasi-variational inequality. We prove that any regular stochastic solution with a dividend
barrier coincides with the value function over all admissible dividend and liquidation policies, and hence is unique within this class. Using $\omega$-scale functions, we construct single-band, singleton, and no-liquidation candidates, establish sufficient conditions for global optimality, and give conditions under which these candidates supply the unique regular barrier stochastic solution. Brownian and Cram\'er--Lundberg examples illustrate how regulatory parameters affect dividend and liquidation boundaries and generate transitions between policy regimes.

\end{abstract}

\section{Introduction}\label{sec:introduction}

In corporate financial management and actuarial risk theory, optimal dividend policy requires balancing immediate distributions to
shareholders against the value of preserving the firm's future operating capacity. Classical models often terminate the firm when its
surplus first crosses a prescribed ruin level. In practice, however, financial distress need not lead to immediate termination. A firm may continue operating while undercapitalized, although prolonged distress can expose it to an increasing risk of regulatory intervention. This distinction between temporary distress and immediate failure motivates occupation-time-based termination mechanisms. The Omega-clock framework provides a natural representation: regulatory exposure accumulates while the surplus remains below a designated solvency threshold, and termination occurs once the accumulated exposure exceeds an independent exponential trigger; see, for example,
\cite{rodosthenous2018beating}.

This paper studies the interaction among dividend distribution, voluntary liquidation, and such regulatory termination. Whenever the surplus falls below the regulatory solvency threshold $y$, the firm enters the regulatory distress region and is exposed to termination with intensity $q$. We assume that regulatory
termination yields zero payoff to shareholders. Before that event occurs, however, shareholders may voluntarily liquidate the firm and receive a state-dependent liquidation payoff. They therefore face a trade-off between securing the current liquidation value and continuing operations in the hope of recovery and future dividend payments. In particular, regulatory distress does not necessarily trigger immediate voluntary liquidation, as shareholders may optimally continue operating at sufficiently low surplus levels while choosing liquidation over an intermediate range of surplus levels.

The classical dividend literature focuses primarily on the payout decision. In the seminal formulation of \cite{DeFinetti-1957}, shareholders maximize the expected discounted dividends until ruin. In a discrete risk model, De Finetti showed that the optimal policy can take a barrier form, i.e., the surplus above a constant level is immediately distributed as dividends. When dividends are absolutely continuous and paid at a bounded rate, \cite{jeanblanc1995optimization} and \cite{AT97} established threshold-type optimal strategies for diffusion models, while \cite{GS06} obtained related results for the Cram\'er--Lundberg model. Other foundational contributions include \cite{hojgaard99}, \cite{schmidli2007stochastic}, and the references therein.

For spectrally negative L\'evy risk processes, \cite{avram2007optimal} used fluctuation identities and scale functions to derive the value of barrier dividend strategies and sufficient conditions for their optimality. \cite{loeffen2008optimality}
subsequently established barrier optimality under complete monotonicity of the L\'evy density. These results provide the fluctuation-theoretic foundation for dividend control under spectrally negative L\'evy risk processes.

A limitation of the classical formulation is that solvency considerations enter only through the  ruin event. The model does not explicitly assign a cost to operating in a distressed state before ruin occurs. This has motivated several extensions of dividend optimization problems incorporating solvency and risk management considerations. For example, \cite{paulsen2003optimal} studies dividend optimization under solvency constraints for diffusion processes, where a solvency related minimum dividend barrier is imposed to perturb the optimization problem. Another widely used approach is to introduce costs or penalties associated with ruin. Under the Cram\'er--Lundberg model and its diffusion approximation, \cite{ThonhauserandAlbrecher-2007} consider dividend optimization with the time value of ruin and study both bounded and unbounded dividend rates. For spectrally negative L\'evy processes, \cite{loeffen2009optimal} considers a fixed penalty at ruin, while \cite{loeffen2010finetti} study affine penalties and establish optimality results for barrier and ``take-money-and-run''
strategies under suitable regularity conditions. Related extensions incorporating capital injections, and other solvency
mechanisms include, among others, \cite{yao2011optimal}, \cite{avram2015gerber},  \cite{xu2020optimal}, and \cite{wang2024optimality}.

These extensions enrich the treatment of insolvency, but corporate exit remains tied only to an exogenously specified ruin event. On the other hand, endogenous liquidation has a separate but related literature.
An early contribution is \cite{porteus1977optimal}, who studies dividend, reinvestment, financing, and liquidation decisions jointly in a dynamic model of a firm. At a more general level, \cite{davis1994problem} formulate a singular stochastic control problem with discretionary stopping, providing a classical precedent for models in which a finite-variation control and an endogenous stopping decision are chosen jointly.
A related probabilistic approach is developed by \cite{de2020optimal}, who connects dividend control under partial information to an auxiliary optimal stopping problem.
In a related setting, \cite{gajek2017complete} consider firm valuation under a spectrally negative L\'evy model with dividend and financing decisions together with the possibility of liquidation. \cite{lindensjo2020optimal} also incorporate an endogenous liquidation option and obtain policies combining dividend control with a liquidation
threshold/constraint. In a model with stochastic profitability, \cite{reppen2020optimal} show that dividend and voluntary
liquidation decisions can generate nontrivial control and stopping regions. These studies demonstrate that dividends payout and exit decisions may interact substantially.

In this paper, we introduce an additional source of termination associated with prolonged financial distress. A firm that falls below a regulatory solvency threshold need not be closed immediately, but continued operation in this region may expose it to regulatory termination. We model this exposure through an Omega clock: while the surplus remains below a prescribed threshold $y>0$, regulatory hazard accumulates at a rate of $q\ge0$, so that termination risk increases
with the accumulated time spent in the regulatory distress region. Shareholders must therefore choose dividends and voluntary liquidation while facing an additional exogenous termination mechanism.

Occupation-time and state-dependent killing mechanisms provide a natural framework for modeling such delayed regulatory termination in financial
distress. Related early developments in finance include the Parisian barrier framework of \cite{chesney1997brownian}, in which activation depends on the duration of excursions beyond a
prescribed barrier, and the state-dependent intensity approach of \cite{carr2000valuation} for valuing executive stock options.
In actuarial risk models, the Omega mechanism developed in \cite{albrecher2011optimal}  allows a company to continue operating in a distressed surplus region while accumulating exposure to bankruptcy or termination. For spectrally negative L\'evy processes, \cite{rodosthenous2018beating} study optimal stopping under a random horizon generated by an Omega clock, and \cite{rodosthenous2020sell} extend occupation-time-based random discounting to a drawdown-dependent setting.

Meanwhile, level-dependent killing for spectrally negative L\'evy processes was developed by
\cite{li2018fluctuations} and \cite{li2018weighted}. The corresponding exit identities and resolvents are expressed in terms of
$\omega$-scale functions defined through Volterra integral equations. This theory, recently, has been extended to one-sided Markov additive processes through $\omega$-scale matrices \cite{czarna2020fluctuation} and to level-dependent spectrally negative
L\'evy processes in \cite{palmowski2025fluctuations}. These fluctuation identities for Omega-killed spectrally negative L\'evy processes provide the analytical tools used later in this paper to construct the candidate value functions under Omega-type regulatory termination.
Recent work by \cite{bodnariu2026optimal} studies dividend control with classical ruin and Omega-type distress default in a drifted Brownian model. Our formulation additionally includes voluntary liquidation with a state-dependent payoff, a solvency floor on dividend payments, and spectrally negative L\'evy surplus dynamics.

Building on these developments, we formulate a joint optimal dividend and voluntary liquidation problem for a firm whose surplus follows a spectrally negative L\'evy process. Dividend payments are modeled as a singular control, whereas voluntary liquidation is represented by an optimal stopping
decision with a state-dependent liquidation payoff. Whenever the surplus falls below an exogenously specified solvency threshold $y$, the firm is exposed to an Omega-type regulatory termination mechanism: the
occupation time spent below $y$ is compared with an independent exponential trigger, and shareholders receive zero terminal payoff if
regulatory termination occurs first. The regulatory threshold and the endogenous voluntary liquidation of shareholders are therefore distinct, allowing the analyzing of regulatory distress, voluntary liquidation, and dividend distribution jointly.

The contributions of this paper are summarized below.

First, we characterize the interaction between dividend distribution,
voluntary liquidation and Omega-type regulatory termination
for a firm with spectrally negative L\'evy surplus and a
state-dependent liquidation payoff. We give
sufficient conditions under which shareholders optimally continue
operating at low surplus, liquidate over an intermediate surplus
band and pay dividends at an upper barrier. The lower liquidation
boundary governs exit following recovery from severe distress,
whereas the upper boundary governs exit following deterioration
from the operating region. This distinction explains the coexistence
of recovery incentives and voluntary liquidation, subject to a
solvency floor on dividend payments.

Second, we develop a probabilistic characterization and verification
framework for the mixed singular-control and stopping problem.
Under the auxiliary dynamic programming and continuity assumptions,
we characterize the value function over an auxiliary Markov control family
as the smallest controlled stochastic supersolution. We formulate
stochastic solutions through controlled supermartingale inequalities
and martingale identities in the strict no-action region. Within
the regular barrier class of these solutions, we construct an
attaining policy by reflecting surplus at the dividend barrier
and liquidating on the obstacle contact set. A verification
argument over all admissible dividend and liquidation policies
identifies each such solution
with both the auxiliary Markov value and the value of the original
problem. This establishes global optimality of the constructed
policy and uniqueness of the solution function within the regular
barrier class. We further give sufficient conditions for the
verified scale-function candidates to belong to this class and
establish existence of a regular barrier stochastic solution under
explicit Brownian parameter conditions with an affine liquidation
payoff and no voluntary liquidation.

Third, we construct single-band, singleton and no-liquidation
candidates using $\omega$-scale functions and derive conditions
selecting their optimal liquidation and dividend boundaries. For each candidate policy, we give sufficient conditions for global optimality among all admissible dividend and liquidation policies. Numerical illustrations in Brownian and
Cram\'er--Lundberg benchmark models show how a higher regulatory
threshold can activate voluntary liquidation and how claim
 jump frequency  can produce non-monotone dividend-barrier
in the reported insurance benchmark model.

The remainder of the paper is organized as follows. Section \ref{sec:formulation} formulates the joint dividend-control and
optimal-liquidation problem under Omega-type regulatory termination. Section \ref{sec:preliminaries} introduces the required scale
functions and fluctuation identities and establishes preliminary properties of the control problem. Section \ref{sec:hjb} introduces
the auxiliary Markov formulation and its associated quasi-variational inequality. Section \ref{sec:SS} develops the controlled stochastic-supersolution characterization, while Section \ref{sec:barrier} establishes attainability and uniqueness in the regular barrier class.
Section \ref{sec:analytical} constructs the semi-analytical barrier candidates, including single-band, singleton, and no-liquidation
regimes, and gives sufficient conditions for global optimality. Section \ref{sec:numerical} provides numerical illustrations and
comparative statics for the Brownian and Cram\'er--Lundberg benchmark
models. Section \ref{sec:conclude} concludes the paper. Some technical proofs are collected in Appendix.

\section{The Model and Problem Formulation}\label{sec:formulation}
Let $X=(X_t)_{t\ge 0}$ be a spectrally negative L\'evy process with L\'evy triplet $(\gamma,\sigma,\nu)$, where $\gamma >0$, $\sigma\ge 0$, and $\nu$ is a L\'evy measure concentrated on $(0,\infty)$ satisfying
\begin{equation}\label{sigma-finite-1}
\int_{(0,\infty)}\left(1\wedge x^2\right)\nu(\dif x) <\infty.
\end{equation}
By defining the L\'evy measure on the positive real line, representing the magnitudes of the downward jumps, the Laplace exponent of $X$ is then given by
\begin{equation}\label{def: psi}
\psi (\theta) := \gamma \theta + \frac{1}{2}\sigma^2\theta^2 - \int_{(0,\infty )}(1-e^{-\theta x} - \theta x\mathds{1}_{\{0< x < 1\}})\nu(\dif x),\qquad  \theta\ge 0,
\end{equation}
where $\mathds{1}_{\{\cdot\}}$ denotes the indicator function.

We formulate the joint optimal dividend and voluntary liquidation problem under Omega-type regulatory termination with respect to the information generated by the L\'evy surplus process.
Let $(\Omega^X,\mathcal F^X,\mathbb F^X,\mathbb P_x^X)$ be a filtered probability space satisfying the usual conditions and carrying a spectrally negative L\'evy process $X=(X_t)_{t\ge0}$ with $X_0=x$, where $\mathbb F^X$ is the augmented natural filtration of $X$. To consider the exogenous regulatory termination and the associated exponential trigger (see e.g., \cite{rodosthenous2018beating}), let $(\Omega^e,\mathcal E,\mathbb P^e)$ support a unit-mean exponential random variable $\mathbf{e}$, and define the product probability space
\[
(\Omega,\mathcal F,\mathbb P_x) := (\Omega^X\times\Omega^e,\,\mathcal F^X\otimes\mathcal E,\, \mathbb P_x^X\otimes\mathbb P^e).
\]
The enlarged decision filtration $\mathbb F=(\mathcal F_t)_{t\ge0}$ in our problem is the usual
augmentation of $\mathcal F_t^X\otimes\{\varnothing,\Omega^e\}$.
In particular, $\mathbf{e}$ is independent of $\mathcal{F}_\infty$ and is not observable by the decision maker.

In this paper, we consider a singular dividend control strategy. To be specific, the dividend control strategy is a non-decreasing, c\`adl\`ag, $\mathbb F$-predictable process $D=(D_t)_{t\ge0}$, with $D_{0-}=0$.
The corresponding controlled surplus is
\begin{equation}\label{def controlprocess}
U^D_t = X_t - D_t,\quad t\ge 0.
\end{equation}

We consider, in addition to the singular control of dividend payments, an endogenous voluntary liquidation decision. Instead of
continuing to operate the firm (or retaining the insurance portfolio), the decision-maker may liquidate the business and receive an immediate state-dependent payoff. Let  $g:\mathbb R\to\mathbb R$ denote the liquidation payoff function. We impose the following
structural condition.

\begin{quote}
\medskip
\noindent \textbf{(A.1)}
The function $g:\mathbb R\to\mathbb R$ is increasing, concave, and continuously differentiable. When $\sigma>0$, it is twice continuously differentiable away from finitely many points. Moreover, there exist
constants $C_0,C_1>0$ such that
\begin{equation}\label{g-bound}
|g(x)|\le C_0+C_1|x|,
\qquad x\in\mathbb R,
\end{equation}
and there exists $\bar x<\infty$ such that $g'(x)<1$, for $x\ge\bar{x}$.
\medskip
\end{quote}

\begin{remark}
Assumption \textbf{(A.1)} imposes an economically natural shape on the liquidation payoff. Concavity captures diminishing marginal
liquidation value as the firm's surplus increases, while $g'(x)<1$ for sufficiently large $x$ ensures that, at high surplus
levels, the marginal liquidation value is smaller than the unit marginal value of an admissible dividend payment. Together with the
linear-growth condition, these assumptions rule out liquidation payoffs that dominate dividend distribution for sufficiently large surplus.
\end{remark}

We next introduce Omega-type regulatory termination. Fix a regulatory solvency threshold $y>0$ and a regulatory termination intensity
$q\ge0$. Whenever the surplus lies below $y$, the firm is in the regulatory distress region and is exposed to exogenous termination at
intensity $q$. Equivalently, if $\mathbf e$ is an independent unit-exponential random variable, regulatory termination occurs when
the accumulated occupation exposure $ q\int_0^t\mathds{1}_{\{U_s^D<y\}}\,\dif s $
first exceeds $\mathbf e$. Thus the regulatory clock runs only when the firm is below the solvency threshold, and a larger value of $q$ corresponds to more severe regulatory termination risk during periods of distress. Our objective is to determine the optimal dividend and voluntary liquidation policies prior to this exogenous termination event.

We follow a solvency constraint for dividend payments similar to \cite{paulsen2003optimal}; see also \cite{lindensjo2020optimal}.
For a given regulatory solvency threshold $y>0$, a dividend strategy $D$ is admissible if
\begin{equation}\label{adimissble-D}
    \int_{[0,t]}
    \mathds{1}_{\{U_s^D<y\}}\,\dif D_s=0, \qquad t\ge0, \quad \text{a.s.}
\end{equation}
 Condition \eqref{adimissble-D} requires dividend payments to leave the controlled surplus at or above the regulatory
solvency threshold. In particular, within the class of barrier strategies considered in later analysis, the dividend barrier must satisfy $b\ge y$. We denote the set of admissible dividend strategies by $\mathcal{D}(y)$.

Let $\mathcal{T}$ denote the set of $\mathbb F$-stopping times taking values in $[0,\infty]$. An element $\tau\in\mathcal{T}$ represents the time of voluntary liquidation. For a fixed $D\in\mathcal{D}(y)$,
define the cumulative regulatory hazard by
\begin{equation}\label{eq:Gamma}
\Gamma_t^{D,y} := q\int_0^t \mathds{1}_{\{U_s^D<y\}}\,\dif s,
\qquad t\ge0,
\end{equation}
where $q\ge0$ is the regulatory termination intensity.  The associated Omega-type regulatory termination time is \begin{equation}\label{eq:zeta}
\zeta^{D,y} := \inf\left\{ t\ge0: \Gamma_t^{D,y}\ge\mathbf{e} \right\},
\end{equation}
with the convention $\inf\varnothing=\infty$. Thus the regulatory clock runs only while the controlled surplus lies below $y$, and a
larger value of $q$ corresponds to greater regulatory termination risk during periods of distress. When $q=0$, regulatory termination is absent.

The random time $\zeta^{D,y}$ may not be an $\mathbb F$-stopping time; it is introduced to give a probabilistic interpretation of the
regulatory termination mechanism. By the independence of $\mathbf{e}$ and $\mathbb{F}$, we have
\begin{equation}\label{eq:survival-omega}
\mathbb P_x \left( \zeta^{D,y}>t\,\middle|\,\mathcal F_t \right) = e^{-\Gamma_t^{D,y}}.
\end{equation}

For $(D,\tau)\in\mathcal{D}(y)\times\mathcal{T}$, define the performance functional by
\begin{align}
J(x;y;D,\tau):=\mathbb E_x\bigg[\int_{[0,\tau)}e^{-rt}\mathds{1}_{\{t<\zeta^{D,y}\}}\,\dif D_t +e^{-r\tau}g(U_\tau^D)\mathds{1}_{\{\tau<\zeta^{D,y}\}}\mathds{1}_{\{\tau<\infty\}}\bigg].
\label{eq:performance}
\end{align}
The interval $[0,\tau)$ in the dividend term specifies that no dividend payment and voluntary liquidation payoff are collected simultaneously.

For $D\in\mathcal D(y)$, define the cumulative discount along the controlled surplus:
\begin{equation}\label{eq:A-D}
A_t^{D,y}:=\int_0^t\omega_y(U_s^D)\,\dif s
=rt+q\int_0^t\mathds{1}_{\{U_s^D<y\}}\,\dif s
=rt+\Gamma_t^{D,y},\qquad
\omega_y(x):=r+q\mathds{1}_{\{x<y\}}.
\end{equation}
We write $A^{0,y}$ for $D\equiv0$ and $A^{b,y}:=A^{D^b,y}$ for the barrier control introduced in Section~\ref{sec:barrier}. The shorthand $A^y$ is used only when the controlled path is unambiguous.

Conditioning on $\mathcal F_t$ and using \eqref{eq:survival-omega}, the performance functional admits the
equivalent representation
\begin{equation}\label{eq:perform}
J(x;y;D,\tau) = \mathbb E_x\left[\int_{[0,\tau)}e^{-A_t^{D,y}}\,\dif D_t+e^{-A_\tau^{D,y}}g(U_\tau^D)\mathds{1}_{\{\tau<\infty\}}\right].
\end{equation}
Therefore, the value function of the joint dividend and voluntary liquidation problem is
\begin{equation}\label{eq:valuefunc}
V(x;y) := \sup_{D\in\mathcal D(y),\,\tau\in\mathcal T} J(x;y;D,\tau).
\end{equation}
Equivalently, using \eqref{eq:perform},
\[
V(x;y) = \sup_{D\in\mathcal D(y),\,\tau\in\mathcal T} \mathbb E_x\left[ \int_{[0,\tau)}e^{-A_t^{D,y}}\,\dif D_t
+ e^{-A_\tau^{D,y}} g(U_\tau^D) \mathds{1}_{\{\tau<\infty\}}\right].
\]
Throughout the paper, $r>0$, $y>0$, and $q\ge0$. The threshold $y$ is fixed within each optimization problem and varied in the comparative statics of Proposition~\ref{prop:valuefunction}(ii) and Section~\ref{sec:numerical}.

To conclude this section, we impose the following condition on the L\'evy measure.

\begin{quote}
  \medskip
\noindent\textbf{(A.2)} The L\'evy measure satisfies:
\begin{enumerate}
\item[(i).] $\displaystyle\int_{(0,\infty)}z\,\nu(\dif z)<\infty$;
\item[(ii).] $\nu$ has a completely monotone density on $(0,\infty)$.
\end{enumerate}

\end{quote}

\begin{remark}
The finite first moment in \textbf{(A.2)} gives the integrability bounds used in verification and finite variation of the jump component. When $\sigma=0$, the continuous drift is therefore written as
$\delta_X:=\gamma+\int_{(0,1)}z\,\nu(\dif z)>0$.
Complete monotonicity supplies the scale-function regularity used in Lemma~\ref{lem:G-continuity}, the differentiated Volterra formulas, and the benchmark tail calculations.

\end{remark}

\section{Preliminaries on Scale Functions and Properties of the Value Function}\label{sec:preliminaries}
In this section, we first introduce the scale functions for spectrally one-sided L\'evy processes, which serve as crucial tools for deriving (semi-)analytical results for the value function.

For each $r\ge0$, let $W^{(r)}$ be the $r$-scale function of $X$, which is a mapping from $\mathbb{R}$ to $[0,\infty)$, where $W^{(r)}(x)=0$ for $x<0$, and for $x\ge0$, it is uniquely defined through its Laplace transform
\begin{equation}\label{def: Wq}
\int_{0}^{\infty}e^{-\theta x}W^{(r)}(x) \dif x = \frac{1}{\psi(\theta) - r},\quad \text{for }\theta >\Phi(r),
\end{equation}
where
\begin{equation}\label{def: phi}
\Phi(r):= \sup\{\theta\ge0: \psi(\theta)=r \},
\end{equation}
and $\psi(\theta)$ is the Laplace exponent given in \eqref{def: psi}.
Another frequently used scale function $Z^{(r)}(x)$ is defined as
\begin{equation}\label{def: Zq}
Z^{(r)}(x) := 1 +r \int_{0}^{x} W^{(r)}(z)\dif z.
\end{equation}
For discussions on the differentiability of the scale function $W^{(r)}$, refer to e.g., \cite{kuznetsov2012theory} and the references therein.

In addition, since our optimal control problem involves Omega killing, we shall also introduce the  $\omega$-scale functions for spectrally negative L\'evy processes (see e.g., \cite{li2018fluctuations}). Let $\omega: \mathbb{R} \to \mathbb{R}_{+}$ be a locally bounded nonnegative measurable function. Then the abovementioned $r$-scale functions can be extended to the version of $\omega$-scale functions, namely when the process $X$ is killed exponentially with state-dependent killing intensity $\omega(\cdot)$, then we let
\begin{equation}\label{eq:omegaW}
\mathcal{W}^{(\omega)}(x) = W(x) + \int_{0}^{x}W(x-y)\omega(y)\mathcal{W}^{(\omega)}(y)\dif y,
\end{equation}
\begin{equation}\label{eq:omegaZ}
\mathcal{Z}^{(\omega)}(x) = 1 + \int_{0}^{x}W(x-y)\omega(y)\mathcal{Z}^{(\omega)}(y)\dif y,
\end{equation}
where $W(x)$ is a classical scale function for $r=0$. The two-sided exit problems with Omega-killing can be stated as follows: for $x\in[0,c]$ and the first exit times $\tau^-_0$ and $\tau^+_c$, we have
\begin{align*}
& \mathbb{E}_x\Bigg[\exp\left(-\int_{0}^{\tau^+_c}\omega(X_t)\dif t \right); \tau^+_c <\tau^-_0\Bigg] = \frac{\mathcal{W}^{(\omega)}(x)}{\mathcal{W}^{(\omega)}(c)},\\
 &  \mathbb{E}_x\Bigg[\exp\left(-\int_{0}^{\tau^-_0}\omega(X_t)\dif t \right); \tau^-_0 <\tau^+_c\Bigg] = \mathcal{Z}^{(\omega)}(x) - \frac{\mathcal{W}^{(\omega)}(x)}{\mathcal{W}^{(\omega)}(c)}\mathcal{Z}^{(\omega)}(c).
\end{align*}
Moreover, the above results can be extended to more general exit identities for any interval $[z,c]$. For $z\le x\le c$,
\begin{align}
&\mathbb{E}_x\Bigg[\exp\left(-\int_{0}^{\tau^+_c} \omega(X_s)\dif s\right); \tau^+_c< \tau^-_z \Bigg] = \frac{\mathcal{W}^{(\omega)}(x,z)}{\mathcal{W}^{(\omega)}(c,z)},\label{eq:exitup}\\
& \mathbb{E}_x\Bigg[\exp\left(-\int_{0}^{\tau^-_z} \omega(X_s)\dif s\right); \tau^-_z < \tau^+_c \Bigg] = \mathcal{Z}^{(\omega)}(x,z) -\frac{\mathcal{W}^{(\omega)}(x,z)}{\mathcal{W}^{(\omega)}(c,z)}\mathcal{Z}^{(\omega)}(c,z),\label{eq:exitdown}
\end{align}
where for any $\delta \ge 0$, \eqref{eq:omegaW} and \eqref{eq:omegaZ} can be expressed as
\begin{align*}
\mathcal{W}^{(\omega)}(x,z) &= W^{(\delta)}(x-z) + \int_{z}^{x}W^{(\delta)}(x-w)(\omega(w)-\delta)\mathcal{W}^{(\omega)}(w,z)\dif w,\\
\mathcal{Z}^{(\omega)}(x,z) &= Z^{(\delta)}(x-z) + \int_{z}^{x} W^{(\delta)}(x-w)(\omega(w)-\delta)\mathcal{Z}^{(\omega)}(w,z)\dif w.
\end{align*}
For $\omega(w)=r+q\mathds{1}_{\{w<y\}}$ and $\delta=r$, these equations become
\begin{align}
\mathcal{W}^{(r,q)}(x,z) &= W^{(r)}(x-z) + q\int_{z}^{x}W^{(r)}(x-w)\mathds{1}_{\{w<y\}}\mathcal{W}^{(r,q)}(w,z)\dif w,\label{omega-scalew}\\
\mathcal{Z}^{(r,q)}(x,z) &= Z^{(r)}(x-z) + q \int_{z}^{x} W^{(r)}(x-w)\mathds{1}_{\{w<y\}}\mathcal{Z}^{(r,q)}(w,z)\dif w.\label{omega-scalez}
\end{align}
Note that Lemma 2.1 of \cite{li2018fluctuations} guarantees that the $\omega$-scale functions are well-defined.
\begin{remark}\label{rem:omega-derivatives}
The functions $\mathcal W^{(r,q)}$ and $\mathcal Z^{(r,q)}$ depend on the fixed threshold $y$. We set $\mathcal W^{(r,q)}(x,z)=0$ for $x<z$ and use
$\mathcal W_1^{(r,q)}(x,z):=\partial_x\mathcal W^{(r,q)}(x,z)$ on open smooth intervals; $\mathcal Z_1^{(r,q)}$ is defined similarly.  For $z<y$, differentiating the Volterra equation gives
\[
\mathcal W_1^{(r,q)}(y+,z)-\mathcal W_1^{(r,q)}(y-,z)
=-qW^{(r)}(0)W^{(r+q)}(y-z).
\]
In bounded variation, $W^{(r)}(0)=1/\delta_X$, giving the derivative jump at threshold $y$, and similarly a value jump in the zero  when letting $x=z$. In the Gaussian case, $W^{(r)}(0)=0$ and the first derivatives agree at $y$.
\end{remark}

\begin{remark}
    when $z>y$, \eqref{omega-scalew} (similar for \eqref{omega-scalez}) reduces to the classical $r$-scale function $W^{r}(x-z)$; and when $z\le x<y$, let $u=x-z$ and $Q(u) = \mathcal{W}^{(r,q)}(u+z,z)$, we can rewrite \eqref{omega-scalew} as
    \[
Q(u) = W^{(r)}(u) + q\int_{0}^{u}W^{(r)}(u-w)Q(w)\dif w.
    \]
Hence, by taking the Laplace transform on both sides, we can easily verify that $\mathcal{W}^{(r,q)}(x,z) = W^{(r+q)}(x-z)$.
\end{remark}

Define $\mathcal I^{(r,q)}$ on $\mathbb R$ by
\begin{equation}\label{eq:Irqx-v2}
\mathcal I^{(r,q)}(x):=
\begin{cases}
q\displaystyle\int_0^\infty e^{-\Phi(r+q)u}W^{(r)}(u+x)\,\dif u,&q>0,\\[1mm]
e^{\Phi(r)x},&q=0.
\end{cases}
\end{equation}
For $q>0$, the integral is finite, and this normalization is $q$ times that in \cite[Eq.~(2.11)]{rodosthenous2018beating}. For $x\le0$, $\mathcal I^{(r,q)}(x)=e^{\Phi(r+q)x}$ in both cases. For $x\ge0$ it also satisfies
\[
\mathcal I^{(r,q)}(x)=e^{\Phi(r+q)x}
-q\int_0^xW^{(r+q)}(x-h)\mathcal I^{(r,q)}(h)\,\dif h.
\]
For $x\le a$, the one-sided exit identity is
\begin{equation}\label{eq:onesideexitH}
\mathbb E_x\!\left[e^{-A^{0,y}_{\tau_a^+}};\tau_a^+<\infty\right]
=\frac{\mathcal I^{(r,q)}(x-y)}{\mathcal I^{(r,q)}(a-y)}.
\end{equation}
For $a>y$ this follows from \cite[Theorem~2.5]{li2018fluctuations}; for $a\le y$ the path stays below $y$ until upward exit, and both sides reduce to $e^{-\Phi(r+q)(a-x)}$. When $q=0$, it is the usual $r$-discounted exit formula.

\begin{proposition}\label{prop:valuefunction}
    Let Assumption (\textbf{A.1}) and integrability condition $\int_{0}^{\infty}z\nu(\dif z) <\infty$
    in Assumption (\textbf{A.2}) hold. Let $r>0, q\ge 0$, and $y>0$. Then:
    \begin{enumerate}
        \item[(i)] For every fixed $y>0$,  the mapping $x\mapsto V(x;y)$ is finite, nonnegative, increasing, and continuous on $\mathbb{R}$. Furthermore,
        \begin{equation}\label{eq:lowerbound-V}
            V(x_1;y) - V(x_2;y) \ge x_1-x_2,\quad \text{for all } x_1>x_2\ge y.
        \end{equation}
        \item[(ii)] For every fixed $x\in\mathbb{R}$, the mapping $ y \mapsto V(x;y)$ is non-increasing on $(0,\infty)$.
        \item[(iii)] There exist constants $C_2, C_3>0$, independent of $x$ and $y$, such that
        \begin{equation}\label{V-bound}
            0\le V(x;y)\le C_2 + C_3 x^+, \qquad x\in\mathbb{R}, y>0,
        \end{equation}
        where $x^+ = \max\{x,0\}$. In particular, the value function has at most affine growth on $\mathbb{R}_+$.
        \item[(iv)] For every fixed $y>0$, the family
        \begin{equation}\label{VD-UI}
            \Bigg\{e^{-A^{D,y}_t}V(U_t;y) + \int_{[0,t]}e^{-A^{D,y}_s}\dif D_s: t\ge 0, D\in\mathcal{D}(y)  \Bigg\}
        \end{equation}
        is uniformly integrable.
    \end{enumerate}
   \end{proposition}
   \begin{proof}
       See Appendix \ref{App.1}.
   \end{proof}

The following estimate provides a common integrable bound for
discounted dividends and terminal values, and will be used in
the verification arguments below.

\begin{lemma}\label{lem:common-envelope}
Under Assumption \textbf{(A.2)}(i), set
$X_t^*:=\sup_{0\le s\le t}|X_s|$ and
$K_r:=r\int_0^\infty e^{-rs}(1+X_s^*)\,\dif s$.
Then $\mathbb E_x[K_r]<\infty$. For every admissible $D$, every finite stopping time $\theta$, and every $f$ with $|f(u)|\le C(1+|u|)$,
\[
\int_{[0,\theta]}e^{-A_s^{D,y}}\,\dif D_s\le K_r,
\qquad e^{-A_\theta^{D,y}}|f(U_\theta^D)|\le2CK_r.
\]
Moreover, $e^{-rt}(1+X_t^*)\to0$ almost surely and in $L^1$ as $t\to\infty$.
\end{lemma}
\begin{proof}
Finite first jump moment gives
$X_t=x+\delta_Xt+\sigma B_t-S_t$, where $S$ is the subordinator of downward jump magnitudes and $\mathbb E[S_t]=t\int z\nu(\dif z)$.
Then, $\mathbb E[X_t^*]\le |x|+C_1t+C_2\sqrt t$ and $\mathbb E[K_r]<\infty$.
Admissibility gives $D_t\le\overline X_t^+\le X_t^*$ and $|U_t^D|\le2X_t^*$.
Since $X_t^*$ is non-decreasing,
\[
e^{-r\theta}(1+X_\theta^*)\le
r\int_\theta^\infty e^{-rs}(1+X_s^*)\,\dif s\le K_r.
\]
Integration by parts bounds the dividend integral by
$e^{-r\theta}X_\theta^*+r\int_0^\theta e^{-rs}X_s^*\,\dif s\le K_r$.
The displayed tail integral tends to zero almost surely and in $L^1$, proving the final assertion.
\end{proof}

\section{Auxiliary Dynamic Programming and Quasi-Variational Inequality}\label{sec:hjb}

The value function \eqref{eq:valuefunc} is defined over $\mathcal D(y)\times\mathcal T$. For dynamic programming principle (DPP) and the stochastic supersolution characterization, we introduce an auxiliary family $\mathcal D^{\Mmark}(y)\subseteq\mathcal D(y)$, and define
\begin{equation}\label{eq:Vm}
V^{\Mmark}(x;y):=\sup_{D\in\mathcal D^{\Mmark}(y),\,\tau\in\mathcal T}J(x;y;D,\tau),
\end{equation}
under the standard Bellman framework of \cite{Fleming-2006}.

\begin{assumption}\label{ass:aux-dpp}
For each fixed $y>0$, $\mathcal D^{\Mmark}(y)$ contains the zero-dividend policy and all barrier policies with $b>y$, and satisfies the standard stability and measurable selection conditions, including admissible initial interventions.
%The conditional Bellman identities \eqref{eq:DPP} and \eqref{eq:continuation-DPP} below hold, and 
$V^{\Mmark}(\cdot;y)$ is 
continuous on $\mathbb R$.
\end{assumption}

Since the zero-dividend policy with $\tau=\infty$ has value zero,
\begin{equation}\label{eq:Vm_le_V}
0\le V^{\Mmark}(x;y)\le V(x;y),\qquad x\in\mathbb R.
\end{equation}
The growth and integrability bounds of Proposition~\ref{prop:valuefunction} therefore hold for $V^{\Mmark}$.
Using classical arguments one can show the following DPP.
\begin{proposition}[Auxiliary DPP]\label{prop:DPP}
Under Assumption~\ref{ass:aux-dpp}, for every $x\in\mathbb R$, $y>0$, and bounded $\mathbb F$-stopping time $\theta$,
\begin{align}\label{eq:DPP}
V^{\Mmark}(x;y)=\sup_{D\in\mathcal D^{\Mmark}(y),\,\tau\in\mathcal T}\mathbb E_x\Bigg[
&\int_{[0,\theta]}e^{-A_t^{D,y}}\mathds{1}_{\{t<\tau\}}\,\dif D_t
+e^{-A_\tau^{D,y}}g(U_\tau^D)\mathds{1}_{\{\tau\le\theta,\ \tau<\infty\}}\nonumber\\
&+e^{-A_\theta^{D,y}}V^{\Mmark}(U_\theta^D;y)\mathds{1}_{\{\theta<\tau\}}\Bigg].
\end{align}
\end{proposition}

For $\varepsilon>0$, define the strict non-liquidation set and its exit time by
\begin{align*}
\mathscr N_\varepsilon^{\Mmark}&:=\{x:V^{\Mmark}(x;y)>g(x)+\varepsilon\},\\
\tau_\varepsilon^D&:=\inf\{t\ge0:U_t^D\notin\mathscr N_\varepsilon^{\Mmark}\}\\
&=\inf\{t\ge0:V^{\Mmark}(U_t^D;y)\le g(U_t^D)+\varepsilon\}.
\end{align*}
Continuity makes $\mathscr N_\varepsilon^{\Mmark}$ open. Admissible dividends may be paid before its exit.

\begin{corollary}[Auxiliary local DPP]\label{cor:continuation-DPP}
Under Assumption~\ref{ass:aux-dpp}, let $\varepsilon>0$ and let
$h^D\le\tau_\varepsilon^D$ be bounded stopping times. Then
\begin{equation}\label{eq:continuation-DPP}
V^{\Mmark}(x;y)=\sup_{D\in\mathcal D^{\Mmark}(y)}\mathbb E_x\left[
\int_{[0,h^D]}e^{-A_s^{D,y}}\,\dif D_s
+e^{-A_{h^D}^{D,y}}V^{\Mmark}(U_{h^D}^D;y)\right].
\end{equation}
\end{corollary}

\begin{proof}
See Appendix~\ref{App:A+1}.
\end{proof}

\begin{proposition}\label{prop:Vm-supermartingale}
Under Assumption~\ref{ass:aux-dpp}, for every $D\in\mathcal D^{\Mmark}(y)$,
\begin{equation}\label{eq:Vm-process}
M_t^{V^{\Mmark},D}:=e^{-A_t^{D,y}}V^{\Mmark}(U_t^D;y)
+\int_{[0,t]}e^{-A_s^{D,y}}\,\dif D_s,\qquad t\ge0,
\end{equation}
is a uniformly integrable $\mathbb F$-supermartingale, and $M_0^{V^{\Mmark},D}\le V^{\Mmark}(x;y)$.
\end{proposition}
\begin{proof}
The common integrable bound in Lemma~\ref{lem:common-envelope} and \eqref{eq:Vm_le_V} give uniform integrability, including at stopped times. For $0\le s<t$, the auxiliary DPP for a fixed policy  yields
\begin{equation}\label{eq:p4.3-1}
\mathbb E_x\left[e^{-A_t^{D,y}}V^{\Mmark}(U_t^D;y)
+\int_{(s,t]}e^{-A_u^{D,y}}\,\dif D_u\,\middle|\,\mathcal F_s\right]
\le e^{-A_s^{D,y}}V^{\Mmark}(U_s^D;y).
\end{equation}
Adding the accrued dividends on $[0,s]$ proves the supermartingale property. The continuation control has increments $D_{s+u}-D_s$. Finally, Assumption~\ref{ass:aux-dpp} gives
$V^{\Mmark}(x;y)\ge D_0+V^{\Mmark}(x-D_0;y)=M_0^{V^{\Mmark},D}$.
\end{proof}

For sufficiently regular test functions with finite jump integral 
%put $\omega_y(x):=r+q\mathds{1}_{\{x<y\}}$ and 
define
\begin{equation}\label{operatorL}
\mathcal L\varphi(x)=\gamma\varphi'(x)+\frac{\sigma^2}{2}\varphi''(x)
+\int_{(0,\infty)}\left[\varphi(x-z)-\varphi(x)+z\varphi'(x)\mathds{1}_{\{z<1\}}\right]\nu(\dif z).
\end{equation}
Assuming that $V^{\Mmark}$ is smooth enough, the local Bellman equation gives $(\mathcal L-\omega_y)V^{\Mmark}=0$ where neither action is exercised. Initial admissible dividends give
$V^{\Mmark}(x;y)-V^{\Mmark}(z;y)\ge x-z$ for $x\ge z\ge y$, and hence $V_x^{\Mmark}\ge1$ wherever that derivative exists above $y$. Immediate liquidation gives $V^{\Mmark}\ge g$. Thus the quasi-variational inequality (QVI) away from the solvency threshold is
\begin{equation}\label{eq:QVI}
\begin{cases}
\displaystyle\max\{(\mathcal L-r)V^{\Mmark}(x),\ 1-V^{\Mmark\prime}(x),\ g(x)-V^{\Mmark}(x)\}=0,&x>y,\\[1mm]
\displaystyle\max\{(\mathcal L-r-q)V^{\Mmark}(x),\ g(x)-V^{\Mmark}(x)\}=0,&x<y.
\end{cases}
\end{equation}
At $y$, we have $\omega_y(y)=r$.

\section{Stochastic Solution and Supersolution Verification}\label{sec:SS}

The following stochastic solution formulation uses the finite-dividend inequality:
\begin{equation}\label{eq:finite-dividend-ineq}
\varphi(x)-\varphi(z)\ge x-z,\qquad x\ge z\ge y.
\end{equation}
This condition covers lump-sum and continuous dividends. To define the strict no-action set, put
\begin{equation}\label{eq:Delta-y}
\Delta_y\varphi(x):=\inf_{y\le z<x}\frac{\varphi(x)-\varphi(z)}{x-z},\quad x>y,
\qquad\Delta_y\varphi(y):=\liminf_{x\downarrow y}\Delta_y\varphi(x).
\end{equation}
For $x>y$, $\Delta_y\varphi(x)>1$ means that the finite-dividend quotients have a uniform positive gap above one.

\begin{definition}\label{def:stochastic_solution}
Fix $y>0$. A continuous function $\varphi:\mathbb R\to\mathbb R$ of at most affine growth is a \emph{controlled stochastic supersolution} if:
\begin{enumerate}
\item[(i)] $\varphi\ge g$ on $\mathbb R$;
\item[(ii)] it satisfies \eqref{eq:finite-dividend-ineq};
\item[(iii)] for every $D\in\mathcal D^{\Mmark}(y)$ with $D_0=0$, the process
\begin{equation}\label{eq:phi-controlled-super}
M_t^{\varphi,D}:=e^{-A_t^{D,y}}\varphi(U_t^D)
+\int_{(0,t]}e^{-A_s^{D,y}}\,\dif D_s
\end{equation}
is a supermartingale after stopping at the exit time from each compact state interval, and its resulting stopped families are uniformly integrable.
\end{enumerate}
Define its \emph{strict no-action set} by
\begin{equation}\label{eq:Cphi}
\mathscr C_\varphi^{\mathrm{str}}:=\{x<y:\varphi(x)>g(x)\}
\cup\{x\ge y:\varphi(x)>g(x),\ \Delta_y\varphi(x)>1\}.
\end{equation}
A controlled stochastic supersolution is a \emph{stochastic solution} if, for every compact interval $I\subset\mathscr C_\varphi^{\mathrm{str}}$, the stopped process
\begin{equation}\label{eq:phi-sub}
\left\{e^{-A_{t\wedge\tau_I}^{0,y}}\varphi(X_{t\wedge\tau_I}):t\ge0\right\},\qquad
\tau_I:=\inf\{t\ge0:X_t\notin I\},
\end{equation}
is a uniformly integrable $\mathbb F$-martingale.
\end{definition}

\begin{remark}\label{rem:controlled-vs-uncontrolled}
Following the stochastic formulation of \cite{avram2015gerber}, we use the controlled surplus $U^D$ and its discount factor $A^{D,y}$ in the supersolution property. We refer to a controlled stochastic supersolution simply as a stochastic supersolution below.
\end{remark}

\begin{proposition}\label{prop:min-stochasticsolution}
Under Assumption~\ref{ass:aux-dpp}, $V^{\Mmark}(\cdot;y)$ is the pointwise smallest stochastic supersolution. Writing $\mathscr G_y^{\Mmark}$ for this class,
\begin{equation}\label{eq:min-super}
V^{\Mmark}(x;y)=\min_{\varphi\in\mathscr G_y^{\Mmark}}\varphi(x),\qquad x\in\mathbb R.
\end{equation}
\end{proposition}
\begin{proof}
Continuity is given in Assumption~\ref{ass:aux-dpp}. Immediate liquidation, the initial-intervention property, and Proposition~\ref{prop:Vm-supermartingale} give the three supersolution conditions for $V^{\Mmark}$; its growth is inherited from \eqref{eq:Vm_le_V}.

Let $\varphi\in\mathscr G_y^{\Mmark}$. An initial dividend is separated using $\varphi(x)\ge D_0+\varphi(x-D_0)$, so first consider $D_0=0$. Put $\rho_n^D:=\inf\{t\ge0:U_t^D\notin[-n,n]\}$. Optional sampling of the stopped supersolution at $t\wedge\tau\wedge\rho_n^D$ gives
\[
\varphi(x)\ge\mathbb E_x\left[\int_{(0,t\wedge\tau\wedge\rho_n^D]}e^{-A_s^{D,y}}\,\dif D_s
+e^{-A_{t\wedge\tau\wedge\rho_n^D}^{D,y}}\varphi(U_{t\wedge\tau\wedge\rho_n^D}^D)\right].
\]
Lemma~\ref{lem:common-envelope} gives an integrable bound uniform in $n$ and in the stopped times. Dominated convergence allows $n\to\infty$ at fixed $t$. On $\{\tau\le t\}$,
\[
\int_{[0,\tau]}e^{-A_s^{D,y}}\,\dif D_s+e^{-A_\tau^{D,y}}\varphi(U_\tau^D)
\ge\int_{[0,\tau)}e^{-A_s^{D,y}}\,\dif D_s+e^{-A_\tau^{D,y}}g(U_\tau^D).
\]
The residual terminal term on $\{\tau>t\}$ tends to zero in $L^1$ by Lemma~\ref{lem:common-envelope} as well. Letting $t\to\infty$ and restoring an initial intervention gives $\varphi(x)\ge J(x;y;D,\tau)$. Taking the auxiliary supremum yields $\varphi\ge V^{\Mmark}$.
\end{proof}

For a closed interval $I$, define
$\tau_I^D:=\inf\{t\ge0:U_t^D\notin\operatorname{int}(I)\}$ and
$\rho^D:=\tau_I^D\wedge\tau_\varepsilon^D$.

\begin{definition}[$\varepsilon$-controlled local supersolution]\label{def:controlsuper}
A continuous function $\phi$ of at most affine growth is an $\varepsilon$-controlled local supersolution on $I$ if it satisfies \eqref{eq:finite-dividend-ineq},
\begin{equation}\label{eq:eps-boundary}
\phi(x)\ge V^{\Mmark}(x;y),\qquad x\notin\operatorname{int}(I)\cap\mathscr N_\varepsilon^{\Mmark},
\end{equation}
and, for every $D\in\mathcal D^{\Mmark}(y)$ with $D_0=0$,
\begin{equation}\label{eq:Mtcontrolledsuper}
\left\{e^{-A_{t\wedge\rho^D}^{D,y}}\phi(U_{t\wedge\rho^D}^D)
+\int_{(0,t\wedge\rho^D]}e^{-A_s^{D,y}}\,\dif D_s:t\ge0\right\}
\end{equation}
is a uniformly integrable supermartingale.

\end{definition}

\begin{proposition}\label{prop:controlsuper}
Under Assumption~\ref{ass:aux-dpp}, for every $\varepsilon>0$ and closed interval $I$,
\begin{equation}\label{eq:epsilon-min}
V^{\Mmark}(x;y)=\min_{\phi\in\mathscr L_I^{\varepsilon,y,\Mmark}}\phi(x),\qquad x\in I,
\end{equation}
where $\mathscr L_I^{\varepsilon,y,\Mmark}$ is the class in Definition~\ref{def:controlsuper}.
\end{proposition}
\begin{proof}
See Appendix~\ref{app:proof-controlsuper}.
\end{proof}

\begin{remark}\label{rem:G-subset-L}
For every $\varepsilon>0$ and closed interval $I$,
\begin{equation}\label{eq:G-subset-L}
\mathscr G_y^{\Mmark}\subseteq\mathscr L_I^{\varepsilon,y,\Mmark}.
\end{equation}
Indeed, the upper-bound argument in Proposition~\ref{prop:min-stochasticsolution} gives the exterior boundary inequality. Further stopping the compact-state supersolution at $\tau_\varepsilon^D$ preserves its supermartingale property, and Lemma~\ref{lem:common-envelope} supplies the required integrability.
\end{remark}

The following verification theorem gives sufficient conditions
for a candidate function $\varphi$ to satisfy
$\varphi(x)\ge V(x;y)$ for all $x\in\mathbb R$.

\begin{theorem}
\label{prop:direct-verification}
Let
$\varphi:\mathbb{R}\to\mathbb{R}$ be locally Lipschitz
and of at most affine growth.
Suppose that, for some finite set $\mathcal{B}$ containing $y$\footnote{$\mathcal B$ is a finite set containing $y$ and every breakpoint at which the smoothness required in the Theorem fails. For example, in the single-band case we discuss in Section \ref{sec:analytical}, $\mathcal B$ can be taken as $\{s_1,s_2,y,b\}$ together with the relevant nonsmooth points of $g$.},
$\varphi\in C^1(\mathbb{R}\setminus\mathcal{B})$
and has finite one-sided first derivatives at $\mathcal{B}$.
If $\sigma>0$, assume additionally that
$\varphi\in C^2(\mathbb{R}\setminus\mathcal{B})$
and is locally a difference of two convex functions.
Assume that:
\begin{enumerate}
\item[(a)]
$\varphi\ge g$ on $\mathbb{R}$;
\item[(b)]
$\varphi(x)-\varphi(z)\ge x-z$
for all $x\ge z\ge y$;
\item[(c)]
\begin{equation}
\label{eq:gen-super}
(\mathcal{L}-\omega_y)\varphi(x)\le 0,
\qquad x\in\mathbb{R}\setminus\mathcal{B};
\end{equation}
\item[(d)]
if $\sigma>0$,
\begin{equation}
\label{eq:verification-local-time}
\varphi'(a+)-\varphi'(a-)\le 0,
\qquad a\in\mathcal{B}.
\end{equation}
\end{enumerate}
Then
\begin{equation}
\label{eq:phi_ge_V}
\varphi(x)\ge V(x;y),
\qquad x\in\mathbb{R}.
\end{equation}
\end{theorem}

\begin{proof}
See Appendix~\ref{App:A+2}.
\end{proof}

\section{Barrier Martingales, Attainment, and Uniqueness}\label{sec:barrier}

We first establish martingale identities for surplus reflected
at a dividend barrier. We then prove that every regular barrier
stochastic solution equals the value function over all admissible
dividend and liquidation policies, with an optimal policy obtained
by reflection at the barrier and liquidation where the solution
equals the liquidation payoff. Consequently, there is at most one
such solution. Section~\ref{sec:regular-existence} gives sufficient
conditions under which the explicit candidates are regular
barrier stochastic solutions.

For $b>y$, define
\begin{equation}\label{eq:Db}
D_t^b:=\sup_{0\le s\le t}(X_s-b)^+,\qquad Y_t^b:=X_t-D_t^b\le b.
\end{equation}
Under $\mathbb P_x$, $D_0^b=(x-b)^+$ and $Y_0^b=x\wedge b$. Since $X$ has no positive jumps, $D^b$ is continuous on $(0,\infty)$ and increases only at $Y^b=b$. It is predictable and belongs to $\mathcal D(y)$, and Assumption~\ref{ass:aux-dpp} includes it in $\mathcal D^{\Mmark}(y)$.
\begin{equation}\label{eq:Aby}
A_t^{b,y}:=A_t^{D^b,y}=rt+q\int_0^t\mathds{1}_{\{Y_s^b<y\}}\,\dif s.
\end{equation}
For a candidate $h$ with barrier $b$, define
\begin{equation}\label{eq:S-C-barrier}
\mathscr S_h:=\{x:h(x)=g(x)\},\qquad
\mathscr K_{h,b}:=(-\infty,b)\setminus\mathscr S_h,\qquad
\mathscr E_{h,b}:=\mathscr K_{h,b}\cup\{b\}.
\end{equation}
The first set is the obstacle contact set; $\mathscr K_{h,b}$ is the interior continuation set, and $\mathscr E_{h,b}$ includes the reflecting boundary where dividends may be paid. The associated liquidation time is $\tau_h^b:=\inf\{t\ge0:Y_t^b\in\mathscr S_h\}$.

\begin{lemma}\label{lem:Ymartingale}
Fix $b>y$, $I=[a,c]$ with $a<c\le b$, and $x\in I$.
Let $h:\mathbb R\to\mathbb R$ satisfy the regularity and growth
assumptions stated in Theorem~\ref{prop:direct-verification}, with a finite breakpoint
set $\mathcal B_h$ containing $y$ whenever $y\in(a,c)$.
Assume
\begin{equation}\label{eq:h-harmonic}
(\mathcal L-\omega_y)h(u)=0,
\qquad u\in(a,c)\setminus\mathcal B_h.
\end{equation}
If $\sigma>0$, assume also that $h'(u-)=h'(u+)$ for every $u\in\mathcal B_h\cap(a,c)$..

\begin{enumerate}
\item[(i)]
If $c<b$, let
$T_I:=\inf\{t\ge0:X_t\notin(a,c)\}$.
Then
$\{e^{-A_{t\wedge T_I}^{0,y}}h(X_{t\wedge T_I}):t\ge0\}$
is a uniformly integrable martingale under $\mathbb P_x$.

\item[(ii)]
If $c=b$, assume
\begin{equation}\label{eq:smooth-b}
h'(b-)=1,
\end{equation}
and, when $\sigma=0$,
\begin{equation}\label{eq:barrier-BV-harmonic}
\delta_Xh'(b-)-rh(b)
+\int_{(0,\infty)}[h(b-z)-h(b)]\nu(\dif z)=0.
\end{equation}
Let $\tau_I^b:=\inf\{t\ge0:Y_t^b\le a\}$.
Then
\begin{equation}\label{eq:reflected-martingale}
\left\{
e^{-A_{t\wedge\tau_I^b}^{b,y}}
h(Y_{t\wedge\tau_I^b}^b)
+\int_{[0,t\wedge\tau_I^b]}e^{-A_s^{b,y}}\,\dif D_s^b
:t\ge0
\right\}
\end{equation}
is a uniformly integrable martingale under $\mathbb P_x$.
\end{enumerate}
\end{lemma}

\begin{proof}
See Appendix~\ref{App:B}.
\end{proof}

\begin{definition}\label{def:regular-barrier-solution}
A stochastic solution $h$ is a \emph{regular barrier stochastic solution at $b>y$} if $\mathscr S_h\subset(-\infty,b)$ and:
\begin{enumerate}
\item[(i)] $h(x)=h(b)+x-b$ for $x\ge b$, $h'(b-)=1$, and, in the bounded-variation case, \eqref{eq:barrier-BV-harmonic} holds;
\item[(ii)] $\Delta_yh(x)>1$ for $x\in[y,b)\setminus\mathscr S_h$;
\item[(iii)]
$h$ satisfies the regularity assumptions stated at the beginning
of Theorem~\ref{prop:direct-verification}, with a finite breakpoint
set $\mathcal B_h$ containing $y$, and
\[
(\mathcal L-\omega_y)h(x)=0,
\qquad x\in\mathscr K_{h,b}\setminus\mathcal B_h.
\]
When $\sigma>0$, assume additionally that $h\in C^1(\mathbb R)$.
\end{enumerate}
\end{definition}

\begin{proposition}\label{prop:barrier-attain}
If $h$ is a regular barrier stochastic solution at $b>y$, then
\begin{equation}\label{eq:h-performance}
h(x)=J(x;y;D^b,\tau_h^b),\qquad x\in\mathbb R.
\end{equation}
\end{proposition}
\begin{proof}
See Appendix~\ref{App:C+1}.
\end{proof}

\begin{theorem}\label{thm:barrier-uniqueness}
Under \textbf{(A.1)}--\textbf{(A.2)}, let the auxiliary family contain the zero-dividend policy and every interior barrier policy, as in Assumption~\ref{ass:aux-dpp}. Every regular barrier stochastic solution $h$ at some $b>y$ satisfies
\begin{equation}\label{eq:h=Vm}
h(x)=V^{\Mmark}(x;y)=V(x;y),\qquad x\in\mathbb R.
\end{equation}
Therefore, the regular barrier class contains at most one function, and $(D^b,\tau_h^b)$ is globally optimal for every member of this class.
\end{theorem}
\begin{proof}
As a stochastic supersolution, $h$ satisfies $h\ge g$ and
the finite-dividend inequality \eqref{eq:finite-dividend-ineq}.
Its growth and regularity satisfy the assumptions of
Theorem~\ref{prop:direct-verification} by
Definition~\ref{def:regular-barrier-solution}.

Let $x\notin\mathcal B_h$ and choose a compact interval
$I\subset\mathbb R\setminus\mathcal B_h$ containing $x$
in its interior. Set
$T_I:=\inf\{t\ge0:X_t\notin I\}$.
The supermartingale property under the zero-dividend policy
and the discounted Dynkin formula give
\[
\mathbb E_x\left[
\int_0^{t\wedge T_I}e^{-A_s^{0,y}}
(\mathcal L-\omega_y)h(X_{s-})\,\dif s
\right]\le0.
\]
Dividing by $t$ and letting $t\downarrow0$, continuity of
$(\mathcal L-\omega_y)h$ on $I$ yields
$(\mathcal L-\omega_y)h(x)\le0$.

When $\sigma>0$, the condition $h\in C^1(\mathbb R)$ gives
$h'(a+)-h'(a-)=0$ at every breakpoint in $\mathcal B_h$.
All conditions of Theorem~\ref{prop:direct-verification}
therefore hold, and $h\ge V$.

Since $D^b\in\mathcal D^{\Mmark}(y)$,
Proposition~\ref{prop:barrier-attain} gives
\[
h(x)=J(x;y;D^b,\tau_h^b)
\le V^{\Mmark}(x;y)\le V(x;y)\le h(x),
\qquad x\in\mathbb R.
\]
Thus $h=V^{\Mmark}=V$ and $(D^b,\tau_h^b)$ is globally
optimal. Every regular barrier stochastic solution equals
$V$, proving uniqueness within this class.
\end{proof}

\section{Semi-analytical Construction and Global Verification}\label{sec:analytical}
We construct the values of single-band, singleton, and no-liquidation policies using $\omega$-scale functions. Theorem~\ref{prop:direct-verification} establishes global optimality when the stated verification conditions hold. We then give conditions under which these candidates provide the unique regular barrier stochastic solution of Section~\ref{sec:barrier}.
We first treat interior dividend barriers $b>y$ and conclude with the scenario $b=y$ when $q=0$.

Fix initially
\begin{equation}\label{eq:fixed-triplet}
s_1\le s_2<b,\qquad b>y,
\end{equation}
and consider the connected liquidation region $[s_1,s_2]$,
the upper dividend barrier $b$, and the continuation set
\[
(-\infty,s_1)\cup(s_2,b).
\]

The two liquidation boundaries have distinct economic roles. Below $s_1$, the firm continues operating and liquidates upon recovery to $s_1$. A firm deteriorating from the upper continuation region liquidates upon entering the band. A downward jump may overshoot the band and land below $s_1$, returning the firm to the recovery region.

The construction of the candidate boundaries in this section reflects this asymmetry.  We first derive
the analytical expressions of the value associated with a fixed triplet $(s_1,s_2,b)$.  We then identify the candidate lower recovery-liquidation boundary $s_1^*$ from the associated up-crossing optimal liquidation problem, select the dividend barrier for each tentative upper boundary $s_2$, and determine $s_2$ using appropriate smooth-fit or continuous-fit conditions.  Finally, we give sufficient conditions under which the resulting single-band and barrier dividend candidate is globally optimal.  We also treat singleton stopping and the absence of voluntary liquidation separately. These results provide a complete verification of the candidate policies covered by the present construction.

To derive the explicit semi-analytical expression for this candidate value function utilizing $\omega$-scale functions, we recall several exit times associated with the underlying surplus dynamics. Let
\[
T_a^+:=\inf\{t\ge0:X_t\ge a\},
\qquad
\tau_a^-(b):=\inf\{t\ge0:Y_t^b<a\}.
\]
Since $X$ has no positive jumps, a path starting below $s_1$ enters the stopping band from below at $s_1$ continuously. Thus the lower continuation value is independent of the upper barrier and upper stopping boundary.

For fixed triplet $(s_1,s_2,b)$ we define an auxiliary function
\begin{equation}\label{eq:tilde-g}
\widetilde g_{s_1,s_2}(z):=
\begin{cases}
\displaystyle g(s_1)\frac{\mathcal{I}^{(r,q)}(z-y)}{\mathcal{I}^{(r,q)}(s_1-y)},&z<s_1,\\[2mm]
g(z),&s_1\le z\le s_2.
\end{cases}
\end{equation}
where $\mathcal{I}^{(r,q)}(\cdot)$ is the $\omega$-scale function given by \eqref{eq:Irqx-v2}.

The candidate value function under the triplet strategy $(s_1,s_2,b)$ is
\begin{equation}\label{eq:vc12bxy}
v_{s_1,s_2,b}(x;y)=
\begin{cases}
\displaystyle g(s_1)\frac{\mathcal I^{(r,q)}(x-y)}{\mathcal I^{(r,q)}(s_1-y)},&x<s_1,\\[2mm]
g(x),&s_1\le x\le s_2,\\[1mm]
\displaystyle
\mathbb E_x\left[e^{-A_{\tau_{s_2}^-(b)}^{b,y}}\widetilde g_{s_1,s_2}(Y_{\tau_{s_2}^-(b)}^b)
+\int_{[0,\tau_{s_2}^-(b)]}e^{-A_t^{b,y}}\,\dif D_t^b\right],&s_2<x<b,\\[2mm]
v_{s_1,s_2,b}(b;y)+x-b,&x\ge b.
\end{cases}
\end{equation}
The first line in \eqref{eq:vc12bxy} follows from the upward-exit identity \eqref{eq:onesideexitH}, with $x<s_1$. It covers targets on either side of $y$: if $s_1\le y$, the ratio is $e^{-\Phi(r+q)(s_1-x)}$, because the path stays below $y$ before the upward exit. All discounted exit terms are understood to be zero when the corresponding exit time is infinite.

\subsection{The reflected resolvent and the explicit candidate value function}

The reflected identities in \cite{li2018fluctuations} are stated under a continuity assumption on the killing function. Our killing rate $\omega_y=r+q\mathds{1}_{(-\infty,y)}$ has a single jump, hence, we shall first prove the following lemma.

\begin{lemma}\label{lem:step-omega-reflected}
The upper-reflected exit and resolvent identities of \cite[Theorems 2.3--2.4]{li2018fluctuations} extend to the bounded step function $\omega_y(x)=r+q\mathds{1}_{\{x<y\}}$, where the corresponding $\omega$-scale functions are the two-variable functions $\mathcal W^{(r,q)}$ and $\mathcal Z^{(r,q)}$ defined in Section \ref{sec:preliminaries}.
\end{lemma}

\begin{proof}
See Appendix \ref{App:C}.
\end{proof}

For $x\in[s_2,b)$, define the reflected $\omega$-resolvent
\begin{equation}\label{eq:resolvent}
R_{s_2,b}^{(r,q)}(x,\dif z)
:=\int_0^\infty
\mathbb E_x\left[e^{-A_t^{b,y}};
 t<\tau_{s_2}^-(b),\ Y_t^b\in\dif z\right]\dif t.
\end{equation}
By Lemma \ref{lem:step-omega-reflected},
\begin{equation}\label{eq:resolvent-density}
R_{s_2,b}^{(r,q)}(x,\dif z)
=\ell_{s_2}^{(r,q)}(x,b)\delta_b(\dif z)
+\ell_{s_2}^{(r,q)}(x,z)\dif z,
\end{equation}
where, for $z\in[s_2,b)$,
\begin{align}
\ell_{s_2}^{(r,q)}(x,z)
&=\frac{\mathcal W^{(r,q)}(x,s_2)}
        {\mathcal W_1^{(r,q)}(b,s_2)}
  \mathcal W_1^{(r,q)}(b,z)
  -\mathcal W^{(r,q)}(x,z),
\label{eq:ell-density}\\
\ell_{s_2}^{(r,q)}(x,b)
&=\frac{\mathcal W^{(r,q)}(x,s_2)}
        {\mathcal W_1^{(r,q)}(b,s_2)}W(0).
\label{eq:ell-atom}
\end{align}
The derivative conventions are given in Remark~\ref{rem:omega-derivatives}. At the reflecting point $b>\max\{y,z\}$, $\mathcal W_1^{(r,q)}(b,z)$ is an ordinary first-variable derivative. The resolvent has the atom \eqref{eq:ell-atom} in bounded variation; this atom is zero when $\sigma>0$.

For the dividend term in the third line of \eqref{eq:vc12bxy}, write
\[
h_b(x):=
\frac{\mathcal W^{(r,q)}(x,s_2)}
     {\mathcal W_1^{(r,q)}(b,s_2)},
\qquad x\le b.
\]
Since $s_2<b$ and $b>y$, the denominator is strictly positive. Indeed, if $s_2\ge y$, $\mathcal W_1^{(r,q)}(b,s_2)= W^{(r)\prime}(b-s_2)>0$, whereas if $s_2<y<b$,
\[
\mathcal W_1^{(r,q)}(b,s_2)= W^{(r)\prime}(b-s_2) + q\int_{s_2}^{y}
W^{(r)\prime}(b-z)
W^{(r+q)}(z-s_2)\,\dif z
>0.
\]
For $s_2\le x\le b$, the expected discounted dividends are
\begin{equation}\label{eq:dividend-reflected}
    \mathbb{E}_x\left[\int_{[0,\tau_{s_2}^-(b)]}e^{-A_t^{b,y}}\,\dif D_t^b\right] = \frac{\mathcal W^{(r,q)}(x,s_2)}
     {\mathcal W_1^{(r,q)}(b,s_2)}.
\end{equation}
In bounded variation, $\dif D_t^b=\delta_X\mathbf 1_{\{Y_t^b=b\}}\dif t$ for $t>0$. Since $W(0)=1/\delta_X$, integrating the atom \eqref{eq:ell-atom} proves \eqref{eq:dividend-reflected} for $s_2<x<b$.

When $\sigma>0$, $W(0)=0$, and $h_b$ extends continuously by zero below $s_2$. Extend it affinely above $b$ with slope $h_b'(b-)=1$. Its scale representation is $\omega_y$-harmonic on $(s_2,b)\setminus\{y\}$, with matching first derivatives at $y$ by Remark~\ref{rem:omega-derivatives}. Lemma~\ref{lem:Ymartingale} applies on $[s_2,b]$ up to the first hitting of $s_2$, which equals $\tau_{s_2}^-(b)$ almost surely. The local time at $s_2$ and the terminal value are both zero at this first hitting. The residual on survival to time $t$ is bounded by $e^{-rt}\sup_{s_2\le u\le b}h_b(u)$ and tends to zero. Lemma~\ref{lem:common-envelope} therefore gives \eqref{eq:dividend-reflected} in the infinite-horizon limit. For bounded-variation endpoint states, let $F(x)$ denote the expected dividends on the left of \eqref{eq:dividend-reflected}. Dividends begin at the first hitting of $b$, so the strong Markov property and \eqref{eq:exitup} give $F(x)=\mathcal W^{(r,q)}(x,s_2)F(b)/\mathcal W^{(r,q)}(b,s_2)$ for $s_2\le x<b$. The interior formula yields $F(b)=h_b(b)$ and $F(s_2)=h_b(s_2)$. In the Gaussian case the reflected martingale covers $x=b$, and both sides are zero at $x=s_2$.

We further define the transformed overshoot function, for $z\ge s_2$, by
\begin{equation}\label{eq:tilde-g-nu}
\widetilde g_\nu(z)
:=\int_{(z-s_2,\infty)}
\left(\widetilde g_{s_1,s_2}(z-h)-g(s_2)\right)\nu(\dif h).
\end{equation}
Assumption (\textbf{A.1}) and
$\int_{(0,\infty)}z\nu(\dif z)<\infty$ in (\textbf{A.2}) ensure that this integral is well defined.  Set
\begin{equation}\label{eq:Hgxc1}
H_{\widetilde g}(x,s_2)
:=g(s_2)\mathcal Z^{(r,q)}(x,s_2)
-\int_{s_2}^x
 \mathcal W^{(r,q)}(x,z)\widetilde g_\nu(z)\,\dif z.
\end{equation}
For $x>s_2$ and $x\ne y$, differentiation in the first variable gives the following identity. At $y$, the formula uses the corresponding one-sided derivatives.
\begin{align}\label{eq:Hg-derivative}
H_{\widetilde g,1}(x,s_2)
=g(s_2)\mathcal Z_1^{(r,q)}(x,s_2)
-W(0)\widetilde g_\nu(x) -\int_{s_2}^x
 \mathcal W_1^{(r,q)}(x,z)\widetilde g_\nu(z)\,\dif z.
\end{align}

\begin{proposition}\label{prop:fixed-triplet-value}
For $s_2<x<b$, the expectation in the third line of
\eqref{eq:vc12bxy} equals
\begin{equation}\label{eq:Hgxc2}
H_{\widetilde g}(x,s_2)
+\frac{\mathcal W^{(r,q)}(x,s_2)}
       {\mathcal W_1^{(r,q)}(b,s_2)}
 \left(1-H_{\widetilde g,1}(b,s_2)\right).
\end{equation}
Define
\begin{equation}\label{def:G-fun}
G(b,s_2)
:=\frac{1-H_{\widetilde g,1}(b,s_2)}
        {\mathcal W_1^{(r,q)}(b,s_2)},
\qquad b>\max\{s_2,y\}.
\end{equation}
Then
\begin{equation}\label{eq:vxy}
v_{s_1,s_2,b}(x;y)=
\begin{cases}
\displaystyle
g(s_1)\frac{\mathcal I^{(r,q)}(x-y)}
              {\mathcal I^{(r,q)}(s_1-y)},
&x<s_1,\\[2mm]
g(x),&s_1\le x\le s_2,\\[1mm]
H_{\widetilde g}(x,s_2)
+G(b,s_2)\mathcal W^{(r,q)}(x,s_2),
&s_2<x<b,\\[1mm]
x-b+v_{s_1,s_2,b}(b;y),&x\ge b.
\end{cases}
\end{equation}
In particular,
\begin{equation}\label{eq:barrier-smooth-auto}
\frac{\partial}{\partial x}
v_{s_1,s_2,b}(b-;y)=1.
\end{equation}
\end{proposition}

\begin{proof}
See Appendix~\ref{app:proof-fixed-triplet-value}.
\end{proof}
Under \textbf{(A.2)}, $W^{(r)}\in C^\infty(0,\infty)$; see \cite[Theorem 3]{loeffen2008optimality}. When $\sigma>0$, $W^{(r)\prime}(0+)=2/\sigma^2$ is finite and $W^{(r)\prime\prime}$ is locally integrable up to zero: $W^{(r)\prime\prime\prime}>0$ makes $W^{(r)\prime\prime}$ monotone, and the finite limit of $W^{(r)\prime}$ bounds its integral near zero. The overshoot integrand is bounded by $Ch$ for small jumps and by $C(1+h)$ for large jumps, so dominated convergence gives continuity of $\widetilde g_\nu$ on $[s_2,\infty)$. The Volterra and convolution representations yield piecewise $C^1$ regularity in bounded variation and piecewise $C^2$ regularity with locally integrable second derivatives when $\sigma>0$. The breakpoints comprise $y$, the policy boundaries, and the finitely many second-derivative breakpoints of $g$ allowed in \textbf{(A.1)}. Concavity gives locally bounded variation of $g'$, whose one-sided limits agree because $g\in C^1$. With the fit conditions below, the candidate is continuous at its policy boundaries. Local Lipschitz continuity, affine growth, and the finite first moment give finite, locally continuous nonlocal terms. These properties provide the regularity for verification.

For a fixed triplet policy $(s_1,s_2,b)$, no dividend is paid and no voluntary liquidation occurs while the controlled surplus lies in $(-\infty,s_1)\cup(s_2,b).$ The preceding representation therefore implies that the corresponding
candidate value function is $\omega_y$-harmonic on the interior of each of these regions; that is, $(\mathcal L-\omega_y)v_{s_1,s_2,b}(x)=0,$ away from the boundary points and, where necessary, the regulatory threshold $y$.

\begin{proposition}
\label{prop:fixed-triplet-harmonic}
Fix $s_1\le s_2<b$ with $b>y$, and let $v_{s_1,s_2,b}$ be given by \eqref{eq:vxy}. Define the continuation region of this candidate policy, in which neither dividends nor voluntary liquidation occur, by
\[
\mathscr C_{s_1,s_2,b} := (-\infty,s_1)\cup(s_2,b).
\]
Let $I=[a,c]$ be a compact interval such that $I\subset \mathscr C_{s_1,s_2,b},$ and define $T_I:=\inf\{t\ge0:X_t\notin(a,c)\}.$ Then, for every $x\in(a,c)$,
\begin{equation}\label{eq:fixed-policy-harmonic}
v_{s_1,s_2,b}(x;y)
=
\mathbb E_x\left[
e^{-A_{T_I}^{0,y}}
v_{s_1,s_2,b}(X_{T_I};y)
\right].
\end{equation}
Thus $v_{s_1,s_2,b}$ is $\omega_y$-harmonic in the
probabilistic sense on $\mathscr C_{s_1,s_2,b}$.

Moreover, at every point
$x\in\mathscr C_{s_1,s_2,b}\setminus\{y\}$ at which
$v_{s_1,s_2,b}$ has the regularity required by the generator, we have
\begin{equation}\label{eq:fixed-policy-generator}
(\mathcal L-\omega_y(x))
v_{s_1,s_2,b}(x;y) = 0.
\end{equation}
\end{proposition}

\begin{proof}
See Appendix~\ref{app:proof-fixed-triplet-harmonic}.
\end{proof}

\subsection{Candidate boundary selection}

The boundary construction proceeds in three steps in this subsection.  We first identify the lower recovery-liquidation boundary $s_1$, which is determined by the one-sided up-crossing problem below the stopping band.  Then, for each tentative upper liquidation boundary $s_2$, we select the dividend barrier $b$ by maximizing the coefficient $G(\cdot,s_2)$ given by \eqref{def:G-fun}.  Finally, the upper liquidation boundary $s_2$ is determined by the appropriate smooth-fit or continuous-fit condition.

\paragraph{Identification of the lower recovery-liquidation boundary}
We consider connected single-band candidates with a lower recovery-liquidation boundary $s_1<y$.
Let $\underline{k}:=\sup\{u\in\mathbb R:g(u)\le0\},$ and suppose throughout this subsection that
$\underline{k}<\infty$ and $\underline{k}<y$.  By continuity and monotonicity of $g$, we have $ g(\underline{k})=0$, and $g(u)>0$ for $u>\underline{k}$.  Then, for $u>\underline{k}$ define
\begin{equation}
\label{eq:F-lower}
F_y(u):=
\frac{g(u)}{\mathcal I^{(r,q)}(u-y)}.
\end{equation}

\begin{proposition}\label{prop:s1-identification}
Assume $q\ge0$ and
\begin{equation}\label{eq:threshold-regime-condition}
g'(\underline{k})>0,
\qquad
\frac{g'(y)}{g(y)}<\Phi(r+q).
\end{equation}
Then $F_y$ admits a unique maximizer $s_1^*$ on $(\underline{k},y]$. Moreover,
$
\underline{k}<s_1^*<y,
$
and $s_1^*$ is the unique solution on $(\underline{k},y)$ of
\begin{equation}\label{eq:s1*}
g'(u)-\Phi(r+q)g(u)=0.
\end{equation}
Hence, we have
\begin{equation}\label{def:s1*}
  \operatorname{argmax}_{\underline{k}<u\le y}F_y(u) = \{s_1^*\}.
\end{equation}
\end{proposition}
\begin{proof}
See Appendix \ref{App:D+1}.
\end{proof}

For later comparison, it is useful to distinguish the lower recovery-liquidation boundary $s_1^*$ we identified in Proposition \ref{prop:s1-identification} from the ``best'' liquidation level within the unrestricted family of up-crossing stopping strategies. Define
\begin{equation}
\label{eq:Fbar-global}
\overline F_y :=\sup_{u>\underline{k}}F_y(u) =
\sup_{u>\underline{k}}
\frac{g(u)}{\mathcal I^{(r,q)}(u-y)},
\end{equation}
and let $\mathcal M_F(y):=\left\{u>\underline{k}:F_y(u)=\overline F_y\right\}$. Whenever $\mathcal M_F(y)$ is nonempty and bounded, we define
\begin{equation}\label{eq:uF-global}
\widehat u_y:=\sup\mathcal M_F(y).
\end{equation}
\begin{remark}
 The quantities $s_1^*$ and $\widehat u_y$ serve different purposes.  The level $s_1^*$ is selected only among stopping levels below $y$ and represents the candidate lower recovery-liquidation boundary of a connected stopping component in such a triplet policy $(s_1,s_2,b)$. By contrast, $\widehat u_y$, when it exists, maximizes the
payoff within the larger family of strategies that commit to stopping at an up-crossing level.
Indeed, for any $u>\underline{k}$ and $x<u$, the
one-sided Omega-killed exit identity gives
\[
\mathbb E_x\!\left[
e^{-A_{T_u^+}^{0,y}}g(u)
\right]
=
\mathcal I^{(r,q)}(x-y)F_y(u).
\]
The coefficient $\overline F_y$ provides a benchmark for liquidation at an upward-exit target. If $\widehat u_y>y$, the corresponding policy
continues through the regulatory threshold before liquidating at $\widehat u_y$. The lower boundary $s_1^*$ remains selected from levels below $y$.

The global benchmark $\overline F_y$ will be used later when comparing voluntary liquidation with the no-liquidation candidate.
\end{remark}

\paragraph{Dividend barrier selection for a tentative upper liquidation boundary.
For $\max\{s_2,y\}<x<b$, equations~\eqref{eq:vxy} and \eqref{def:G-fun} give}
\begin{equation}\label{eq:gradient-G-identity}
\frac{\partial}{\partial x} v_{s_1,s_2,b}(x;y)-1 = \mathcal W_1^{(r,q)}(x,s_2) \bigl[G(b,s_2)-G(x,s_2)\bigr],\quad \max\{s_2,y\}<x<b.
\end{equation}
Whenever $\mathcal W_1^{(r,q)}(x,s_2)>0$, choosing $b$ as a global maximizer of $G(\cdot,s_2)$ yields
\[
\frac{\partial}{\partial x} v_{s_1,s_2,b}(x;y)\ge1, \qquad \max\{s_2,y\}<x<b.
\]
We first establish the regularity of $G(\cdot,s_2)$ needed for this
selection procedure.  We then give sufficient conditions ensuring
that a finite rightmost global maximizer of $G$ exists.

Fix $s_1$ and $s_2\ge s_1$, and put $m(s_2):=\max\{s_2,y\}$.
\begin{lemma}
\label{lem:G-continuity}
Under Assumptions \textbf{(A.1)}--\textbf{(A.2)}, for fixed
$s_1\le s_2$, the mappings
\[
z\longmapsto\widetilde g_\nu(z),
\qquad
x\longmapsto H_{\widetilde g,1}(x,s_2),
\qquad
x\longmapsto\mathcal W_1^{(r,q)}(x,s_2)
\]
are continuous on their respective domains
$[s_2,\infty)$, $(m(s_2),\infty)$, and
$(m(s_2),\infty)$.  Then, $b\longmapsto G(b,s_2)$ is continuous on $(m(s_2),\infty)$.
\end{lemma}

\begin{proof}
See Appendix \ref{App:D}.
\end{proof}
For any candidate $f$ that is affine with slope one above its barrier $b>y$, define the common dividend-region residual by
\begin{equation}\label{Jvb}
J_f(x):=(\mathcal L-r)f(x;y),\qquad x>b,
\end{equation}
This notation applies to single-band, singleton, and no-liquidation candidates. In particular, for $v_b:=v_{s_1,s_2,b}$, we have $v_b'(x;y)=1$ and $v_b''(x;y)=0$ for $x>b$, so
\begin{equation*}
%\label{eq:Jvb-expanded}
J_{v_b}(x)
=\gamma-rv_b(x;y)
+\int_{(0,\infty)}
\left[
 v_b(x-z;y)-v_b(x;y)
 +z\mathds{1}_{\{z<1\}}
\right]\nu(\dif z).
\end{equation*}

The next identity links the tail shape of $G$ to $J_{v_b}$.
\begin{lemma}
\label{lem:G-increment}
Fix $s_1\le s_2$, let $b>m(s_2)$, and let $c>0$.  Then
\begin{equation}\label{eq:G-increment}
\mathcal W_1^{(r,q)}(b+c,s_2)
\bigl(G(b+c,s_2)-G(b,s_2)\bigr)
=
\int_{[0,c)}
 J_{v_b}(b+c-u)\,W^{(r)}(\dif u).
\end{equation}
\end{lemma}

\begin{proof}
See Appendix~\ref{app:proof-G-increment}.
\end{proof}

\begin{corollary}
\label{cor:G-eventual-decrease}
If there exists $B>m(s_2)$ such that
\begin{equation}\label{eq:J-uniform-tail}
J_{v_b}(x)\le0,
\qquad x>b,\quad b\ge B,
\end{equation}
then $b\mapsto G(b,s_2)$ is non-increasing on $[B,\infty)$.
If, in addition, for every $b\ge B$ and $c>0$,
\[
\int_{[0,c)} \mathds{1}_{\{J_{v_b}(b+c-u)<0\}} \,W^{(r)}(\dif u)>0,
\]
then $G(\cdot,s_2)$ is strictly decreasing on $[B,\infty)$.
\end{corollary}

\begin{proof}
See Appendix~\ref{app:proof-G-eventual-decrease}.
\end{proof}

We can now define the dividend barrier for a fixed pair $(s_1,s_2)$.  Define
\begin{align*}
\overline G(s_1,s_2)
:=\sup_{b>m(s_2)}G(b,s_2), \qquad
\mathcal M_G(s_1,s_2)
:=\left\{
 b>m(s_2):
 G(b,s_2)=\overline G(s_1,s_2)
 \right\}.
\end{align*}
Note that $\overline G$ and $\mathcal M_G$  depends on $s_1$ through
$\widetilde g_{s_1,s_2}$ in $G$.  Whenever $\mathcal M_G(s_1,s_2)$ is nonempty and bounded, we define
\begin{equation}\label{eq:b*}
b^*(s_1,s_2):=\sup\mathcal M_G(s_1,s_2).
\end{equation}

\begin{lemma}
\label{lem:finite-band-barrier}
Let Assumptions \textbf{(A.1)}--\textbf{(A.2)} hold and fix
$s_1\le s_2$.  Suppose that there exists a finite level
$B_G(s_1,s_2)>m(s_2)$ such that
\begin{equation}\label{eq:G-eventual-decrease}
b\longmapsto G(b,s_2)
\quad\text{is strictly decreasing on }
[B_G(s_1,s_2),\infty),
\end{equation}
and that
\begin{equation}\label{eq:G-boundary-gap}
\sup_{b>m(s_2)}G(b,s_2)
>
\limsup_{b\downarrow m(s_2)}G(b,s_2).
\end{equation}
Then $G(\cdot,s_2)$ attains its global maximum at a finite interior point.  Therefore,
$\mathcal M_G(s_1,s_2)$ is nonempty and bounded and the rightmost maximizer $b^*(s_1,s_2)$ in \eqref{eq:b*} is finite.
\end{lemma}

\begin{proof}
By Lemma \ref{lem:G-continuity}, $G(\cdot,s_2)$ is continuous on $(m(s_2),\infty)$.  Condition \eqref{eq:G-boundary-gap} implies that there exists $\delta>0$ such that the global supremum cannot be approached on $(m(s_2),m(s_2)+\delta)$.  By
\eqref{eq:G-eventual-decrease},  $G(B_G(s_1,s_2),s_2) > G(b,s_2)$ for all $b>B_G(s_1,s_2)$.  Hence
\[
\sup_{b>m(s_2)}G(b,s_2)
=
\max_{b\in[m(s_2)+\delta,\,B_G(s_1,s_2)]}
G(b,s_2).
\]
Continuity on this compact interval gives an interior maximizer.
\end{proof}

For the selected candidate barrier, \eqref{eq:gradient-G-identity} immediately yields
\begin{equation}\label{eq:upper-gradient}
\frac{\partial}{\partial x}
v_{s_1,s_2,b^*(s_1,s_2)}(x;y)\ge1,
\qquad
m(s_2)<x<b^*(s_1,s_2).
\end{equation}

\paragraph{Selection of the upper liquidation boundary.}
With the candidate lower recovery-liquidation boundary $s_1=s_1^*$ fixed as in \eqref{def:s1*}, we next determine the right endpoint of a candidate connected liquidation region. For each tentative upper liquidation boundary $a>s_1^*$ for which Lemma \ref{lem:finite-band-barrier} applies, let $b^*(s_1^*,a)$ denote the associated dividend barrier defined by \eqref{eq:b*}.  The level $a$ separates the candidate
liquidation region $[s_1^*,a]$ from the upper continuation region $(a,b^*(s_1^*,a))$. Accordingly, the third line in \eqref{eq:vxy}, approached from the right of $a$, must be compatible with the
liquidation payoff $g(a)$ at the boundary.  When $\sigma>0$, regularity of the diffusion component leads to the smooth fit condition (see \eqref{eq:s2_unbounded} below) required at $a$, whereas in the bounded-variation case, we require the continuous fit condition (see \eqref{eq:s2_bounded} below).
Define
\begin{align}
\mathcal S_{\rm sf}
&:=\left\{
 a>s_1^*:
 H_{\widetilde g,1}(a+,a)
 +G(b^*(s_1^*,a),a)
  \mathcal W_1^{(r,q)}(a+,a)
 =g'(a)
 \right\},
\label{eq:s2_unbounded}\\
\mathcal S_{\rm cf}
&:=\left\{
 a>s_1^*:
 H_{\widetilde g}(a+,a)
 +G(b^*(s_1^*,a),a)
  \mathcal W^{(r,q)}(a+,a)
 =g(a)
 \right\}.
\label{eq:s2_bounded}
\end{align}
The root sets range over levels $a>s_1^*$ for which Lemma~\ref{lem:finite-band-barrier} defines $b^*(s_1^*,a)$. The set $\mathcal S_{\rm sf}$ applies when $\sigma>0$. Under \textbf{(A.2)}, $\sigma=0$ gives bounded variation and the continuous-fit set $\mathcal S_{\rm cf}$. In the bounded variation case,
\[
\mathcal W^{(r,q)}(a+,a)=W(0)>0, \qquad H_{\widetilde g}(a+,a)=g(a),
\]
and therefore the continuous-fit condition is equivalent to $G\bigl(b^*(s_1^*,a),a\bigr)=0$. Hence,
\[
\mathcal S_{\rm cf}
=
\left\{a>s_1^*:G\bigl(b^*(s_1^*,a),a\bigr)=0\right\}.
\]
Each root defines a candidate policy with liquidation on $[s_1^*,a]$ and dividends at $b^*(s_1^*,a)$. The next subsection gives sufficient conditions for its global optimality.

We select $s_2^*$ from the fit roots satisfying the global verification conditions below.

\subsection{The single-band candidate and verification}

To proceed, suppose an admissible root $s_2^*$ exists in \eqref{eq:s2_unbounded} (for $\sigma>0$) or \eqref{eq:s2_bounded} (for $\sigma=0$),  set
\begin{equation}\label{eq:bstar-final}
b^*:=b^*(s_1^*,s_2^*).
\end{equation}
The resulting candidate function is
\begin{equation}\label{eq:vxy*}
v^*(x;y)=
\begin{cases}
\displaystyle
g(s_1^*)\frac{\mathcal I^{(r,q)}(x-y)}
                 {\mathcal I^{(r,q)}(s_1^*-y)},
&x<s_1^*,\\[2mm]
g(x),&s_1^*\le x\le s_2^*,\\[1mm]
H_{\widetilde g}(x,s_2^*)
+G(b^*,s_2^*)\mathcal W^{(r,q)}(x,s_2^*),
&s_2^*<x<b^*,\\[1mm]
x-b^*+v^*(b^*;y),&x\ge b^*.
\end{cases}
\end{equation}

Define $x_g:=\inf\left\{ x\in\mathbb R: g'(u)\le1\ \text{for every }u\ge x \right\}$.
\begin{lemma}
\label{lem:obstacle-vstar}
Suppose Proposition \ref{prop:s1-identification} applies, $s_1^*<s_2^*<b^*$, $b^*>y$, and $x_g\le b^*$.
Assume further on the upper continuation region:
\begin{equation}\label{eq:upper-majorant-band}
v^*(x;y)\ge g(x),
\qquad x\in(s_2^*,b^*].
\end{equation}
Then
\begin{equation}\label{eq:derived-obstacle}
v^*(x;y)\ge g(x),
\qquad x\in\mathbb R.
\end{equation}
\end{lemma}

\begin{proof}
See Appendix~\ref{app:proof-obstacle-vstar}.
\end{proof}

\begin{lemma}\label{lem:dividend-vstar}
Suppose $b^*$ is given by \eqref{eq:bstar-final}.  If $s_2^*>y$, assume additionally that $g'(s_2^*)\ge1$. Then
\begin{equation}\label{eq:derived-dividend}
v^*(x;y)-v^*(x-\eta;y)\ge\eta,
\qquad
x\ge y,\quad 0\le\eta\le x-y.
\end{equation}
\end{lemma}

\begin{proof}
By \eqref{eq:upper-gradient}, $v_x^*\ge1$ on
$(\max\{y,s_2^*\},b^*)$, while $v_x^*=1$ above $b^*$.  If $s_2^*>y$, concavity of $g$ and $g'(s_2^*)\ge1$ give
$v_x^* = g'(x)\ge1$ on $[y,s_2^*]$.  If $s_2^*\le y$, the stopping band has no nontrivial intersection with the range $[y,\infty)$, Since $s_1^*<y$ by Proposition \ref{prop:s1-identification} and $v^*$ is continuous and piecewise absolutely continuous on $[y,\infty)$,  \eqref{eq:upper-gradient} directly gives \eqref{eq:derived-dividend}.
\end{proof}

Let $\mathcal B_g$ be the finite set of possible second-derivative breakpoints of $g$ under \textbf{(A.1)}; set $\mathcal B_g=\varnothing$ in bounded variation. For $x\in(s_1^*,s_2^*)\setminus(\{y\}\cup\mathcal B_g)$, define
\begin{equation*}
%\label{eq:stopping-residual}
R_{v^*}^{\rm stop}(x) :=(\mathcal L-\omega_y(x))g(x)+\int_{(x-s_1^*,\infty)}\left[g(s_1^*)
\frac{\mathcal I^{(r,q)}(x-z-y)}{\mathcal I^{(r,q)}(s_1^*-y)}-g(x-z)\right]\nu(\dif z).
\end{equation*}

\begin{lemma}\label{lem:generator-vstar}
Suppose
\begin{equation}\label{eq:S2a}
R_{v^*}^{\rm stop}(x)\le0,
\qquad x\in(s_1^*,s_2^*)\setminus(\{y\}\cup\mathcal B_g),
\end{equation}
and
\begin{equation}\label{eq:S3}
J_{v^*}(x)\le0,
\qquad x>b^*.
\end{equation}
Then
\begin{equation}\label{eq:derived-generator}
(\mathcal L-\omega_y(x))v^*(x;y)\le0
\end{equation}
at every interior point outside $\{s_1^*,s_2^*,y,b^*\}\cup\mathcal B_g$.
\end{lemma}

\begin{proof}
See Appendix~\ref{app:proof-generator-vstar}.
\end{proof}

\begin{theorem}[Optimality of a single-band liquidation and barrier dividend policy]
\label{thm:optimality_condition}
Suppose Assumptions \textbf{(A.1)}--\textbf{(A.2)} hold,
$q\ge0$, and \eqref{eq:threshold-regime-condition} holds.  Suppose Lemma \ref{lem:finite-band-barrier} applies and the boundary-selection procedure above yields finite levels
\[
s_1^*<s_2^*<b^*, \qquad b^*>y.
\]
Assume $x_g\le b^*$ and \eqref{eq:upper-majorant-band} hold, and when $s_2^*>y$, assume additionally $g'(s_2^*)\ge1$.
Assume also \eqref{eq:S2a} and \eqref{eq:S3} hold.

If $\sigma=0$ and $s_2^*\le y$, assume in addition that
\begin{equation}\label{eq:BV-threshold-boundary}
\delta_X-rv^*(y;y)
+\int_{(0,\infty)}
\bigl(v^*(y-z;y)-v^*(y;y)\bigr)\nu(dz)
\le 0.
\end{equation}

Then, with $\tau^*:=\inf\{t\ge0:Y_t^{b^*}\in[s_1^*,s_2^*]\}$, one has
\begin{equation}\label{eq:band-global-value}
V(x;y)=V^{\Mmark}(x;y)=v^*(x;y),
\qquad x\in\mathbb R,
\end{equation}
and $(D^{b^*},\tau^*)$ is globally optimal over
$\mathcal D(y)\times\mathcal T$.
\end{theorem}

\begin{proof}
See Appendix \ref{App:E+1}.
\end{proof}

\begin{remark}
\label{rem:s2-selection}
If several roots in $\mathcal S_{\rm sf}$ (when $\sigma>0$) or $\mathcal S_{\rm cf}$ (when $\sigma=0$) satisfy Theorem~\ref{thm:optimality_condition}, their policies are all globally optimal and their value functions coincide with $V$. We use the smallest verified root as the representative upper liquidation boundary in the numerical implementation.
\end{remark}

\subsection{Singleton stopping and no-liquidation}

The liquidation band can collapse to a singleton, $s_1^*=s_2^*$, or the policy can choose to continue indefinitely. We next construct and verify these two cases.

\paragraph{Singleton candidate.}
We first consider the case in which the liquidation band collapses to a single point. Let $a=s_1^*<y$ be the level identified under the hypotheses of Proposition~\ref{prop:s1-identification}, and let $b_{\rm s}>\max\{a,y\}$ be a global maximizer of $G(\cdot,a)$ as in \eqref{eq:b*}, with its finiteness ensured by Lemma \ref{lem:finite-band-barrier} applied with
$s_1=s_2=a$. Define
\begin{equation}\label{eq:singleton-value}
v^{\rm s}(x;y):=v_{a,a,b_{\rm s}}(x;y),
\end{equation}
and let $\tau_{\rm s}:=\inf\{t\ge0:Y_t^{b_{\rm s}}=a\}.$ Note that $v^s$ has the same analytical expression as given by \eqref{eq:vxy} except that the liquidation band degenerates to one point.
 When $\sigma>0$, according to the expression of $v_{a,a,b_{\rm s}}(x;y)$, we directly have
\[
v_x^{\rm s}(a-;y) =\Phi(r+q)g(a) =g'(a),
\]
where the second equality follows from the first-order condition defining the maximizer $a$ of $F_y$ (see \eqref{eq:F-lower}). Then, \eqref{eq:verification-local-time} in Condition (d) of Theorem \ref{prop:direct-verification} may reduce to
\begin{equation}\label{eq:singleton-derivative-condition}
v_x^{\rm s}(a+;y)\le g'(a).
\end{equation}
When $\sigma=0$, $v^{\rm s}(a+;y)=g(a)+G(b_{\rm s},a)W(0)$ with $W(0)>0$. Thus continuous fit, $v^{\rm s}(a+;y)=v^{\rm s}(a;y)=g(a)$, is equivalent to
\begin{equation}\label{eq:singleton-cont-fit}
G(b_{\rm s},a)=0.
\end{equation}
Equation~\eqref{eq:singleton-cont-fit} is an additional fit condition on the selected barrier.
\begin{lemma}\label{lem:singleton-dividend}
Let $a=s_1^*<y$ be the lower recovery-liquidation boundary identified in Proposition \ref{prop:s1-identification}, and let $b_{\rm s}>y$ be a global maximizer of $G(\cdot,a)$. Then
\begin{equation}\label{eq:singleton-dividend}
v^{\rm s}(x;y)-v^{\rm s}(x-\eta;y)\ge\eta,
\qquad
x\ge y,\quad 0\le\eta\le x-y.
\end{equation}
\end{lemma}
\begin{proof}
Since $a<y$ and $b_{\rm s}$ is a global maximizer of $G(\cdot,a)$, the argument in Lemma \ref{lem:dividend-vstar}, applied with $s_1=s_2=a$ and $b=b_{\rm s}$, gives the result.
\end{proof}

\begin{theorem}[Optimality of a singleton liquidation and barrier dividend policy]
\label{thm:singleton}
Suppose \textbf{(A.1)}--\textbf{(A.2)} and the hypotheses of Proposition~\ref{prop:s1-identification} hold. Let $a=s_1^*<y$ and let Lemma~\ref{lem:finite-band-barrier} supply the finite barrier $b_{\rm s}>y$ for $s_1=s_2=a$. Let $v^{\rm s}$ be defined by \eqref{eq:singleton-value}. Assume
\begin{equation}\label{eq:singleton-majorant}
v^{\rm s}(x;y)\ge g(x),
\qquad x\in(a,b_{\rm s}],
\end{equation}
and $x_g\le b_{\rm s}$. Assume further
\begin{equation}\label{eq:singleton-tail}
J_{v^{\rm s}}(x)\le0, \qquad x>b_{\rm s}.
\end{equation}
At the singleton liquidation level $a$, assume \eqref{eq:singleton-derivative-condition} when $\sigma>0$, and \eqref{eq:singleton-cont-fit} when $\sigma=0$.

In the bounded-variation case, assume additionally that
\begin{equation}\label{eq:singleton-BV-y}
\delta_X-rv^{\rm s}(y;y)
+\int_{(0,\infty)}
\bigl(v^{\rm s}(y-z;y)-v^{\rm s}(y;y)\bigr)\nu(dz)
\le0.
\end{equation}
Then $V(x;y)=V^{\Mmark}(x;y)=v^{\rm s}(x;y)$, for $x\in\mathbb R$, and $(D^{b_{\rm s}},\tau_{\rm s})$ is globally optimal.
\end{theorem}

\begin{proof}
See Appendix \ref{App:E+2}.
\end{proof}

\paragraph{No voluntary liquidation.}
For the policy $\tau=\infty$, the firm continues indefinitely and pays dividends at $b>y$. Its value is expressed through $\mathcal I^{(r,q)}(\cdot-y)$.

For $b>y$, define
\begin{equation}\label{eq:C0b}
C_0(b)
:=
\frac{1}{\mathcal I^{(r,q)\prime}(b-y)},
\end{equation}
and
\begin{equation}\label{eq:no-liq-value-b}
v_b^{0}(x;y)
:=
\begin{cases}
\displaystyle
C_0(b)\mathcal I^{(r,q)}(x-y),
&x<b,\\[2mm]
\displaystyle
C_0(b)\mathcal I^{(r,q)}(b-y)+x-b,
&x\ge b.
\end{cases}
\end{equation}
The normalization in \eqref{eq:C0b} gives  $(v_b^0)'(b-;y)=1,$ while the second line in \eqref{eq:no-liq-value-b} corresponds to the immediate payment of any surplus above $b$.

For $x\le b$, put $h(x)=C_0(b)\mathcal I^{(r,q)}(x-y)$ and extend $h$ affinely above $b$. The upward-exit identity \eqref{eq:onesideexitH} and the strong Markov property give its probabilistic $\omega_y$-harmonicity below $b$; its scale-function representation gives the corresponding generator equalities away from $y$. Its first derivatives agree at $y$ when $\sigma>0$, and $h'(b-)=1$. In bounded variation, the left limit of $(\mathcal L-r)h$ at $b>y$ is zero. Apply Lemma~\ref{lem:Ymartingale} on $[-n,b]$, for $-n<\min\{x,y\}$, stopped at $t\wedge\inf\{s\ge0:Y_s^b\le-n\}$. The terminal term is bounded by $e^{-rt}\sup_{u\le b}h(u)$ if the lower exit has not yet occurred, and by $C_0(b)e^{\Phi(r+q)(-n-y)}$ otherwise. Both bounds tend to zero as $t,n\to\infty$, and Lemma~\ref{lem:common-envelope} justifies removing localization. Since these lower-exit times increase to infinity, nonnegative dividends also allow monotone convergence. We obtain
\[
\mathbb E_x\left[\int_{[0,\infty)}e^{-A_t^{b,y}}\,\dif D_t^b\right]=\frac{\mathcal I^{(r,q)}(x-y)}{\mathcal I^{(r,q)\prime}(b-y)}=v_b^0(x;y).
\]
If $x>b$, the strategy first pays $x-b$ and thereafter starts from $b$, giving the second line of \eqref{eq:no-liq-value-b}. Hence
\begin{equation}\label{eq:no-liq-fixed-policy}
J(x;y;D^b,\infty)=v_b^0(x;y),
\qquad x\in\mathbb R.
\end{equation}

For $q>0$ and $z>0$, we have $\mathcal I^{(r,q)\prime}(z)=\Phi(r+q)\mathcal I^{(r,q)}(z)-qW^{(r)}(z)$. When $q=0$,
$
\mathcal I^{(r,0)}(z)=e^{\Phi(r)z}.
$
The derivative $\mathcal I^{(r,q)\prime}$ is continuous and positive on $(0,\infty)$, so $C_0$ is positive and continuous on $(y,\infty)$. At zero, $\mathcal I^{(r,q)\prime}(0-)=\Phi(r+q)$ and $\mathcal I^{(r,q)\prime}(0+)=\Phi(r+q)-qW^{(r)}(0)$.

We select the candidate dividend barrier by maximizing $C_0$. In particular, if $b$ is a global maximizer of $C_0$, then for
$y<x<b$,  $C_0(b)\ge C_0(x)$, and therefore
\begin{equation}\label{eq:no-liq-gradient-general}
(v_b^0)'(x;y)
=
\frac{\mathcal I^{(r,q)\prime}(x-y)}
     {\mathcal I^{(r,q)\prime}(b-y)}
\ge1.
\end{equation}
Define
\begin{equation}\label{eq:no-liq-barrier-set}
\mathcal M_0
:=
\operatorname{argmax}_{b>y}C_0(b).
\end{equation}

\begin{lemma}
\label{lem:finite-no-liq}
Suppose
\begin{equation}\label{eq:C0-boundary-gap}
\sup_{b>y}C_0(b)
>
\limsup_{b\downarrow y}C_0(b).
\end{equation}
Then $\mathcal M_0$ is nonempty and bounded.
\end{lemma}

\begin{proof}
See Appendix~\ref{app:proof-finite-no-liq}.
\end{proof}

By Lemma \ref{lem:finite-no-liq}, we define
\begin{equation}\label{eq:b0*}
b_0^*
:=
\max\mathcal M_0.
\end{equation}
For this choice,
\eqref{eq:no-liq-gradient-general} gives
\begin{equation}\label{eq:no-liq-gradient}
(v_{b_0^*}^0)'(x;y)
=
\frac{\mathcal I^{(r,q)\prime}(x-y)}
     {\mathcal I^{(r,q)\prime}(b_0^*-y)}
\ge1,
\qquad y<x<b_0^*,
\end{equation}
while $(v_{b_0^*}^0)'(x;y)=1$, for $x>b_0^*$.
Since $v_{b_0^*}^0$ is continuous and piecewise absolutely continuous
on $[y,\infty)$, it follows that
\begin{equation}\label{eq:no-liq-dividend}
v_{b_0^*}^0(x;y)
-
v_{b_0^*}^0(x-\eta;y)
\ge\eta,
\qquad
x\ge y,\quad 0\le\eta\le x-y.
\end{equation}
It remains to compare $v_{b_0^*}^0$ with the liquidation payoff.
Recall the global upward-exit benchmark $\overline F_y$ in \eqref{eq:Fbar-global}. More generally, its definition is $\sup_{u:g(u)>0}g(u)/\mathcal I^{(r,q)}(u-y)$, with an empty supremum set to zero; this agrees with \eqref{eq:Fbar-global} under the finite-$\underline{k}$ assumption. Set
\begin{equation}\label{eq:Delta0}
\Delta_0
:=
C_0(b_0^*)-\overline F_y.
\end{equation}

\begin{lemma}
\label{lem:no-liq-obstacle}
Suppose $\Delta_0\ge0$ and $x_g\le b_0^*$. Then
\[
v_{b_0^*}^0(x;y)\ge g(x),
\qquad x\in\mathbb R.
\]
\end{lemma}

\begin{proof}
See Appendix~\ref{app:proof-no-liq-obstacle}.
\end{proof}

\begin{theorem}[Optimality of no liquidation and barrier dividend policy]
\label{thm:no-liquidation}
Suppose \textbf{(A.1)}--\textbf{(A.2)} hold and Lemma \ref{lem:finite-no-liq} applies,
$\Delta_0\ge0$, $x_g\le b_0^*$, and
\begin{equation}\label{eq:no-liq-tail}
J_{v_{b^*_0}^{0}}(x)\le0,
\qquad x>b_0^*.
\end{equation}
Then $ V(x;y)=V^{\Mmark}(x;y)=v_{b_0^*}^0(x;y)$, for $x\in\mathbb R$, and $(D^{b_0^*},\infty)$ is globally optimal.
\end{theorem}

\begin{proof}
See Appendix \ref{App:E+3}.
\end{proof}

\paragraph{Dividend barrier at $y$ when $q=0$.}
When $q=0$, $C_0(b)=e^{-\Phi(r)(b-y)}/\Phi(r)$ is strictly decreasing for $b>y$, with $C_0(y)=1/\Phi(r)$. Reflection at $b=y$ is admissible. Set $\Phi_r:=\Phi(r)$ and define
\begin{equation}\label{eq:q0-endpoint-value}
u_y(x):=\begin{cases}
e^{\Phi_r(x-y)}/\Phi_r,&x\le y,\\
x-y+1/\Phi_r,&x>y.
\end{cases}
\end{equation}
This is the performance of $(D^y,\infty)$: the exponential branch in \eqref{eq:q0-endpoint-value} is $r$-harmonic, its slope at $y$ is one, and a localized reflected martingale calculation, with the same vanishing-terminal bounds as above, yields the expected dividends.

\begin{proposition}\label{prop:q0-endpoint}
Suppose \textbf{(A.1)}--\textbf{(A.2)} hold, $q=0$, and $u_y(x)\ge g(x)$ for every $x\in\mathbb R$. Then $(D^y,\infty)$ is globally optimal and $V(x;y)=V^{\Mmark}(x;y)=u_y(x)$.
\end{proposition}
\begin{proof}
See Appendix~\ref{app:proof-q0-endpoint}.
\end{proof}

\subsection{Existence in the regular barrier class}\label{sec:regular-existence}

We first give sufficient conditions for the verified
candidates to belong to the regular barrier class of
Definition~\ref{def:regular-barrier-solution}.
We then establish existence in this class under explicit
Brownian parameter conditions.

\begin{proposition}\label{prop:candidate-regular}
Let $v$ satisfy the hypotheses of Theorem~\ref{thm:optimality_condition}, \ref{thm:singleton}, or \ref{thm:no-liquidation}, and let $b>y$ be its dividend barrier. In the single-band and singleton cases put $\ell=s_2^*$ and $\ell=a$, respectively, and set $m=\max\{\ell,y\}$; in the no-liquidation case put $m=y$. Assume
\begin{equation}\label{eq:regular-candidate-slack}
v(b;y)>g(b),\qquad v_x(x;y)>1\quad(m<x<b),
\end{equation}
and, whenever $y\notin\mathscr S_v$, assume $v_x(y+;y)>1$. Then $v$ is a regular barrier stochastic solution at $b$.
Hence, the regular barrier class is nonempty and
consists of $v$, and the policy $(D^b,\tau_v^b)$ is
globally optimal.
\end{proposition}
\begin{proof}
See Appendix~\ref{app:proof-candidate-regular}.
\end{proof}

We provide a Brownian motion case where the regular barrier solution is explicitly obtained.
\begin{proposition}[A Brownian case with a regular barrier solution]\label{prop:BM-regular-existence}
Let $X_t=x+ct+\sigma B_t$, where $c>0$ and $\sigma>0$, and let
\[
q>\frac{\sigma^2r^2}{2c^2},\qquad g(x)=\eta x-K,\qquad 0<\eta<1.
\]
Let $p>0>n$ be the roots of $cs+\sigma^2s^2/2=r$, put $k=\Phi(r+q)$, and define
\[
A=\frac{k-n}{p-n},\qquad B=\frac{p-k}{p-n},\qquad
z_0=\frac{1}{p-n}\log\frac{-Bn^2}{Ap^2},\qquad b_0=y+z_0.
\]
Then $z_0>0$. If $K\ge\eta b_0$, the function
\[
h(x)=
\begin{cases}
\displaystyle\frac{\mathcal I^{(r,q)}(x-y)}{\mathcal I^{(r,q)\prime}(z_0)},&x<b_0,\\[2mm]
\displaystyle\frac{c}{r}+x-b_0,&x\ge b_0,
\end{cases}
\]
is a regular barrier stochastic solution, $\mathscr S_h=\varnothing$, and $(D^{b_0},\infty)$ is globally optimal, with $h=V=V^{\Mmark}$ whenever the auxiliary family includes the zero and interior-barrier policies.
\end{proposition}
\begin{proof}
The zero L\'evy density satisfies \textbf{(A.2)}, and the affine payoff satisfies \textbf{(A.1)}. For $z\ge0$,
\[
\mathcal I^{(r,q)}(z)=Ae^{pz}+Be^{nz},\qquad A>0,\quad B<0,
\]
while $\mathcal I^{(r,q)}(z)=e^{kz}$ for $z<0$. The condition on $q$ is equivalent to $ck>r$. Thus
\[
\mathcal I^{(r,q)\prime\prime}(0+)=\frac{2(r-ck)}{\sigma^2}<0,
\qquad
\mathcal I^{(r,q)\prime\prime\prime}(z)=Ap^3e^{pz}+Bn^3e^{nz}>0.
\]
Therefore, $z_0$ is the unique positive zero of $\mathcal I^{(r,q)\prime\prime}$ and the unique minimum of its positive first derivative on $[0,\infty)$. This gives $h_x>1$ on $[y,b_0)$, $h_x(b_0-)=1$, and $h_{xx}(b_0-)=0$. Harmonicity at $b_0$ yields $h(b_0)=c/r$. The function is continuously differentiable at $y$ and twice continuously differentiable at $b_0$. It has the required piecewise regularity and affine growth.

For $x\le b_0$, $g(x)\le0<h(x)$. For $x>b_0$,
\[
h(x)-g(x)\ge\frac{c}{r}+(1-\eta)(x-b_0)>0.
\]
The function is $\omega_y$-harmonic below $b_0$, and above it
\[
(\mathcal L-r)h(x)=-r(x-b_0)\le0.
\]
The bound $h_x\ge1$ on $[y,\infty)$ gives the
finite-dividend inequality \eqref{eq:finite-dividend-ineq}.
Apply \eqref{eq:controlled-Meyer-Ito} with $\varphi=h$.
The generator inequality and the finite-dividend inequality
make the drift and dividend contributions nonpositive.
The identities
$h'(y-)=h'(y+)$ and
$h'(b_0-)=h'(b_0+)=1$
eliminate the local-time terms.
Lemma~\ref{lem:common-envelope} removes the localization
and gives the controlled stochastic supersolution property.
TThe finite-dividend quotient estimate in the proof of Proposition~\ref{prop:candidate-regular} supplies the strict condition $\Delta_yh>1$ on $[y,b_0)$. Hence $h$ is a regular barrier stochastic solution with empty contact set, and Theorem~\ref{thm:barrier-uniqueness} gives the conclusion.
\end{proof}

\section{Numerical Illustrations for Benchmark Models}
\label{sec:numerical}

In this section, we illustrate the analytical construction and verification results for two benchmark surplus models: a Brownian
motion with drift (BM) and a Cram\'er--Lundberg (CL) process with exponentially distributed claim sizes.  For each model, the classical
scale functions admit explicit exponential representations, which in turn permit efficient evaluation of the $\omega$-scale functions, the auxiliary function $H_{\widetilde g}$, and the barrier-selection function $G$.  We first derive these specific formulas for the scale functions and verify that $G$ is eventually decreasing  required in \eqref{eq:G-eventual-decrease}.  We then compute the candidate value functions, identify the corresponding policy regimes, and examine their sensitivity to the regulatory and model parameters.

\subsection{Scale functions and the tail shape of \texorpdfstring{$G$}{G} in the benchmark models}

We begin with the classical $r$-scale functions $W^{(r)}$ and
$Z^{(r)}$, from which the functions
$\mathcal I^{(r,q)}$, $\mathcal W^{(r,q)}$ and
$\mathcal Z^{(r,q)}$ introduced in Section
\ref{sec:preliminaries} are constructed.

\paragraph{Drifted Brownian motion.} Consider
\[
X_t=x+ct+\sigma B_t,
\]
where $c\in\mathbb R$ and $\sigma>0$.  Its Laplace exponent is $\psi(\theta)=c\theta+\frac{\sigma^2}{2}\theta^2.$
The scale-function formulas hold for $c\in\mathbb R$, including the positive-drift considered below in the numerical illustration. For $r>0$, let
\[
\Phi_r=\frac{-c+\sqrt{c^2+2r\sigma^2}}{\sigma^2}>0,\qquad R_r = \frac{-c-\sqrt{c^2+2r\sigma^2}}{\sigma^2}<0.
\]
Then
\begin{align}
W^{(r)}(x) &= \frac{e^{\Phi_r x}-e^{R_r x}} {\sqrt{c^2+2r\sigma^2}}, \qquad x\ge0, \label{eq:BM-W}\\
Z^{(r)}(x) &= 1+ \frac{r}{\sqrt{c^2+2r\sigma^2}} \left( \frac{e^{\Phi_r x}-1}{\Phi_r} - \frac{e^{R_r x}-1}{R_r} \right), \qquad x\ge0.
\label{eq:BM-Z}
\end{align}
For the step Omega-killing rate $\omega_y(x)=r+q\mathds{1}_{\{x<y\}}$, the Volterra equations \eqref{omega-scalew}--\eqref{omega-scalez}
give, for $\ell\le x$,
\begin{equation}
\mathcal W^{(r,q)}(x,\ell)
=
\begin{cases}
W^{(r)}(x-\ell),
& \ell\ge y,\\[1mm]
W^{(r+q)}(x-\ell),
& \ell\le x\le y,\\[1mm]
\displaystyle
W^{(r)}(x-\ell)
+
q\int_\ell^y
W^{(r)}(x-z)W^{(r+q)}(z-\ell)\,\dif z,
& \ell<y<x.
\end{cases}
\label{eq:BM-omega-W}
\end{equation}
Similarly,
\begin{equation}
\mathcal Z^{(r,q)}(x,\ell)
=
\begin{cases}
Z^{(r)}(x-\ell),
& \ell\ge y,\\[1mm]
Z^{(r+q)}(x-\ell),
& \ell\le x\le y,\\[1mm]
\displaystyle
Z^{(r)}(x-\ell)
+
q\int_\ell^y
W^{(r)}(x-z)Z^{(r+q)}(z-\ell)\,\dif z,
& \ell<y<x.
\end{cases}
\label{eq:BM-omega-Z}
\end{equation}
Since the classical scale functions are finite sums of
exponentials, all the integrals in
\eqref{eq:BM-omega-W}--\eqref{eq:BM-omega-Z} can either be
evaluated explicitly or computed  numerically.

For $q>0$ and $x>0$,
\begin{equation}
\mathcal I^{(r,q)}(x)
=
\frac{q}{\sqrt{c^2+2r\sigma^2}}
\left(
\frac{e^{\Phi_r x}}{\Phi_{r+q}-\Phi_r}
-
\frac{e^{R_r x}}{\Phi_{r+q}-R_r}
\right).
\label{eq:BM-I}
\end{equation}
Thus  $\Lambda(x) = \mathcal I^{(r,q)\prime}(x) /\mathcal I^{(r,q)}(x)$
is available explicitly, and the lower recovery-liquidation boundary $s_1^*$
can be computed directly from \eqref{eq:s1*}.

Because the Brownian model has continuous sample paths, there is no
jump overshoot across the liquidation band and $\widetilde g_\nu\equiv0$. Hence,
\begin{equation}
H_{\widetilde g}(x,s_2)= g(s_2)\mathcal Z^{(r,q)}(x,s_2).
\label{eq:BM-H}
\end{equation}

The preceding exponential representations also allow the tail condition on the barrier-selection function $G$ (see Lemma \ref{lem:finite-band-barrier} ) to be verified analytically.

\begin{proposition}
\label{prop:G-tail-BM}
Fix $s_1\le s_2$.  In the drifted Brownian model,
there exists a finite $B_G^{\rm BM}(s_1,s_2)>m(s_2)$ such that $ b\longmapsto G(b,s_2)$ is strictly decreasing on
$[B_G^{\rm BM}(s_1,s_2),\infty)$.
Hence condition \eqref{eq:G-eventual-decrease} holds.
\end{proposition}
\begin{proof}
See Appendix \ref{App:E}.
\end{proof}

\paragraph{Cram\'er--Lundberg model.}
Consider
\[
X_t = x+ct-\sum_{i=1}^{N_t}Y_i,
\]
where $c>0$, $N$ is a Poisson process with intensity $\lambda>0$, and the claim sizes are exponentially distributed
with density $\alpha e^{-\alpha h}$, for $h>0$.
Thus $\nu(\dif h)=\lambda\alpha e^{-\alpha h}\,\dif h$, and Assumption \textbf{(A.2)} is satisfied. Under the compensated convention of Section~\ref{sec:formulation}, the L\'evy triplet drift is
\[
\gamma=c-\frac{\lambda}{\alpha}\{1-(1+\alpha)e^{-\alpha}\},
\qquad \delta_X=c.
\]
The reported parameter values in below numerical illustration satisfy the standing condition $\gamma>0$.  The Laplace exponent is
\begin{equation}
\psi(\theta) =c\theta-\lambda\frac{\theta}{\alpha+\theta}.
\label{eq:CL-psi}
\end{equation}
Let the two roots of $\psi(\theta)=r$ be $\rho_1(r)=\Phi_r>0$ and $\rho_2(r)<0$. Then
\begin{equation}
W^{(r)}(x)=A_1e^{\rho_1(r)x}+A_2e^{\rho_2(r)x},\qquad x\ge0,
\label{eq:CL-W}
\end{equation}
where $A_i =1/\psi'(\rho_i(r))$, for $i=1,2$. The functions $\mathcal W^{(r,q)}$ and $\mathcal Z^{(r,q)}$ are obtained from \eqref{omega-scalew}--\eqref{omega-scalez}; because both $W^{(r)}$ and $W^{(r+q)}$ are finite sums of exponentials, the required convolutions are explicit and are evaluated numerically.

For $q>0$ and $x>0$,
\begin{equation}
\mathcal I^{(r,q)}(x)=q\left(\frac{A_1e^{\rho_1(r)x}}{\Phi_{r+q}-\rho_1(r)}+\frac{A_2e^{\rho_2(r)x}}{\Phi_{r+q}-\rho_2(r)}\right).
\label{eq:CL-I}
\end{equation}

For $q=0$, the normalized continuous extension of \eqref{eq:BM-I} and \eqref{eq:CL-I} is $\mathcal I^{(r,0)}(x)=e^{\Phi_r x}$ for all $x\in\mathbb R$, as defined in Section~\ref{sec:preliminaries}.

In addition,  from \eqref{eq:tilde-g} we have
\begin{equation}
\widetilde g_\nu(z)
=
\lambda\alpha
\int_{z-s_2}^{\infty}
\left[
\widetilde g_{s_1,s_2}(z-h)-g(s_2)
\right]
e^{-\alpha h}\,\dif h,
\qquad z\ge s_2.
\label{eq:CL-overshoot}
\end{equation}
Accordingly,
\begin{equation}
H_{\widetilde g}(x,s_2)
=
g(s_2)\mathcal Z^{(r,q)}(x,s_2)
-
\int_{s_2}^x
\mathcal W^{(r,q)}(x,z)
\widetilde g_\nu(z)\,\dif z.
\label{eq:CL-H}
\end{equation}

The exponential distribution for the claim size makes the tail structure of \eqref{eq:CL-overshoot} particularly simple.  This again permits direct verification of the decreasing condition of $G$.

\begin{proposition}
\label{prop:G-tail-CL}
Fix $s_1\le s_2$.  In the Cram\'er--Lundberg model with exponentially distributed claim sizes, there exists a finite $B_G^{\rm CL}(s_1,s_2)>m(s_2)$ such that $ b\longmapsto G(b,s_2)$ is strictly decreasing on $[B_G^{\rm CL}(s_1,s_2),\infty)$.
Hence condition \eqref{eq:G-eventual-decrease} holds.
\end{proposition}
\begin{proof}
See Appendix \ref{App:F}.
\end{proof}

\begin{corollary}[Regular barrier solutions in the benchmark models]
\label{cor:benchmark-regular}
In either benchmark model, let $v$ satisfy the hypotheses of
Theorem~\ref{thm:optimality_condition}, \ref{thm:singleton},
or \ref{thm:no-liquidation}.
If its selected dividend barrier $b>y$ satisfies $g'(b)<1$,
then $v$ is a regular barrier stochastic solution at $b$.
\end{corollary}

\begin{proof}
See Appendix~\ref{app:proof-benchmark-regular}.
\end{proof}

\subsection{Liquidation payoff}

For the numerical illustrations, we use an increasing and concave liquidation payoff defined as
\begin{equation}
g(x)
:=
\begin{cases}
-K+\beta\kappa x,
&x<0,\\[1mm]
\beta(1-e^{-\kappa x})-K,
&x\ge0,
\end{cases}
\label{eq:numerical-g}
\end{equation}
where $\beta>0$ determines the asymptotic liquidation value,
$\kappa>0$ controls its curvature, and $K>0$ is a fixed
liquidation cost.
Moreover, $g$  has at most affine growth, and $g'(x)=\beta\kappa e^{-\kappa x}\longrightarrow0$, as $x\to\infty$. Hence the structural requirements in Assumption \textbf{(A.1)} are
satisfied.

\subsection{Candidate values and policy boundaries}

We next illustrate the candidate value functions and the associated barrier/liquidation policies for the two benchmark models.  For the
Brownian model, we take
\[
c=1.2,\qquad \sigma=2.0,\qquad r=0.05,\qquad q=0.8,\qquad y=4.0.
\]
For the Cram\'er--Lundberg model, we take
\[
c=2.5,\qquad \lambda=1.0,\qquad \alpha=1.0,
\]
while keeping $r=0.05$, $q=0.8$, and $y=4.0$.  In both models the liquidation payoff parameters are
\[
\beta=15,\qquad \kappa=0.5,\qquad K=1.
\]
For each scenario, the candidate boundaries are obtained from the construction in Section \ref{sec:analytical}, and the corresponding sufficient verification conditions are checked using the analytic reductions and numerical diagnostics described below.

\begin{figure}[htbp]
    \centering
    \begin{subfigure}[b]{0.48\textwidth}
        \centering
        \includegraphics[width=\textwidth]{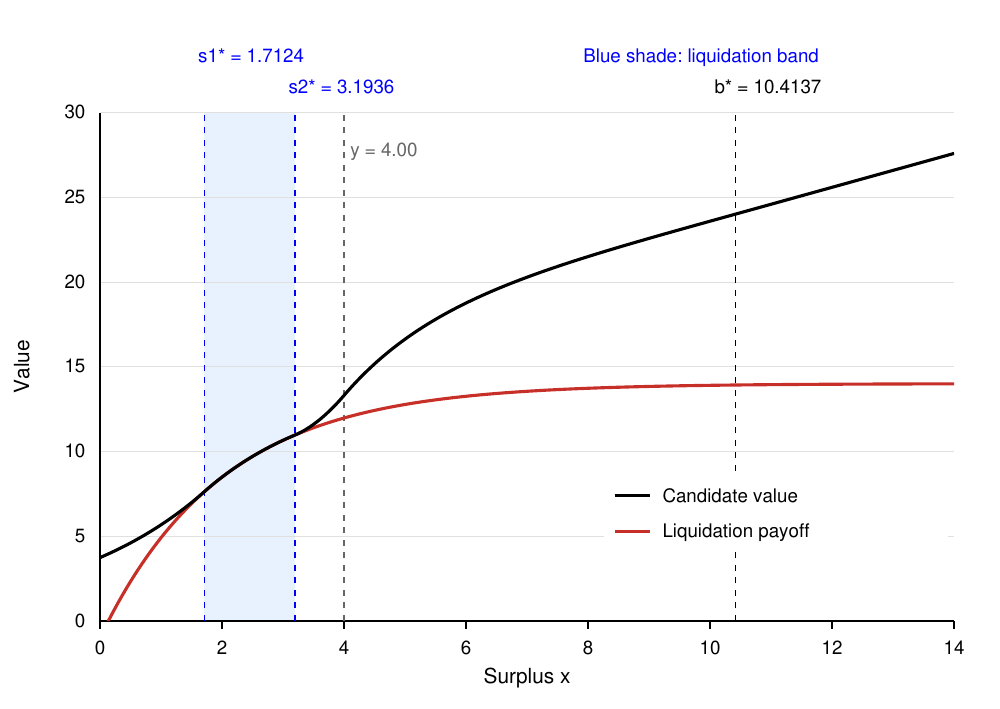}
        \caption{Brownian motion with drift}
        \label{fig:value_BM}
    \end{subfigure}
    \hfill
    \begin{subfigure}[b]{0.48\textwidth}
        \centering
        \includegraphics[width=\textwidth]{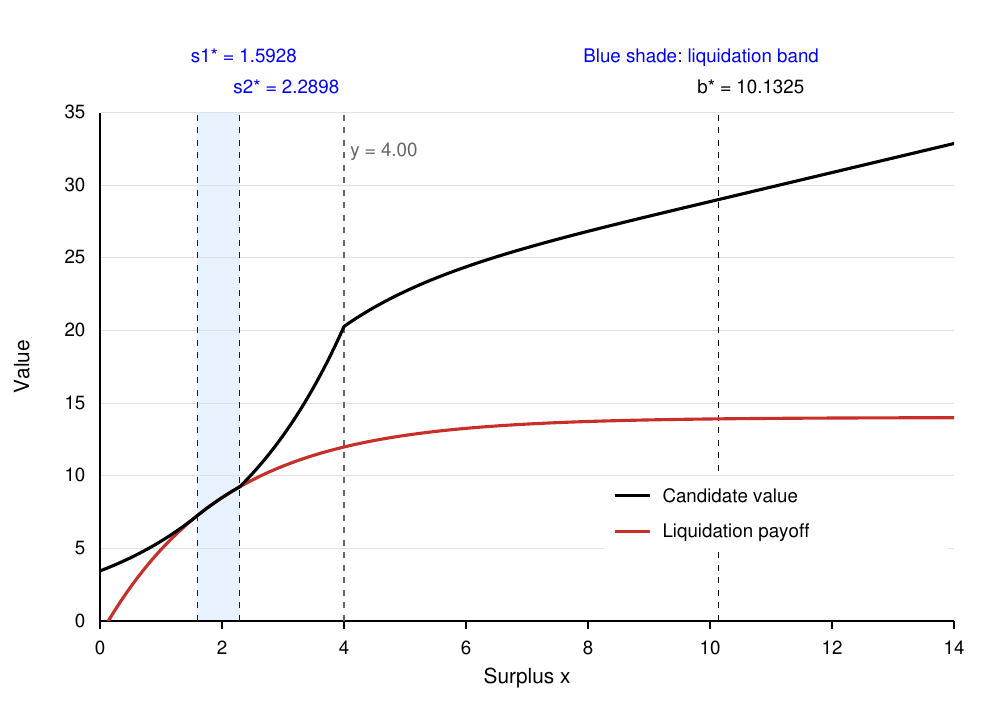}
        \caption{Cram\'er--Lundberg model}
        \label{fig:value_CL}
    \end{subfigure}
    \caption{Benchmark candidate values and liquidation payoffs. The shaded interval is the liquidation band $[s_1,s_2]$; dashed lines mark its endpoints, the regulatory threshold $y$, and the dividend barrier $b$.}
    \label{fig:value_functions_combined}
\end{figure}

Figure \ref{fig:value_functions_combined} displays the resulting value functions together with the continuation, liquidation, and dividend regions.  For the Brownian motion benchmark model, the candidate boundaries are approximately
\[
s_1^*=1.7124,\qquad s_2^*=3.1936,\qquad b^*=10.4137.
\]
Thus voluntary liquidation is optimal on the interval $[1.7124,~3.1936]$, while dividends are reflected at the upper barrier $10.4137$.

For the Cram\'er--Lundberg benchmark model, the corresponding boundaries are
\[
s_1^*=1.5928,\qquad s_2^*=2.2898,\qquad b^*=10.1325.
\]
Hence, voluntary liquidation is optimal on the interval $[1.5928,~2.2898]$ with dividend barrier at $10.1325$ for Cram\'er--Lundberg. Under these benchmark values, the liquidation region is narrower than in the Brownian motion case.  The mean increments per unit time are $1.2$ (BM) and $1.5$ (CL), and the variance rates are $4$ and $2$, respectively. The comparison therefore reflects both parameter scenario and path structure: the Brownian surplus is exposed to continuous fluctuations, whereas the Cram\'er--Lundberg surplus evolves deterministically upward between discrete downward claim jumps.  The sensitivity analysis below shows how these differences evolve when the regulatory and model parameters are varied.

By contrast, the selected dividend barriers are of similar magnitude in the two benchmark cases.  Hence, for the present setup, the difference between the two models is more pronounced in the voluntary liquidation region than in the dividend barrier.

\subsection{Sensitivity Analysis}
Finally, we investigate the sensitivity of the optimal policy to the regulatory solvency threshold $y$, the regulatory termination intensity $q$, and the main parameters governing the surplus dynamics.
Tables \ref{tab:sens_bm} and \ref{tab:sens_cl} report the candidate regimes and boundaries. A ``Band'' entry selects liquidation on $[s_1,s_2]$; ``No Liq.'' selects $\tau=\infty$ and reports only the dividend barrier. ``Inconclusive'' identifies a case for which none of the three candidate policy types passes the sufficient-condition checks.

The last column of Tables \ref{tab:sens_bm} and \ref{tab:sens_cl} reports $\widehat u_F$, the global maximizer of $F_y$ over unrestricted upcrossing stopping strategies. It is an auxiliary liquidation benchmark for this restricted stopping problem. The joint dividend and liquidation problem can select a different lower boundary; $\widehat u_F$ can also lie above $y$, as in Panel A of Table~\ref{tab:sens_cl} for $y=2$.

\paragraph{Brownian model.}
Table \ref{tab:sens_bm} shows that the dividend and liquidation
decisions respond differently to the regulatory and model
parameters.

Panel A considers the regulatory solvency threshold $y$.  At $y=2$,
the selected policy involves no voluntary liquidation, with a
dividend barrier $b^*=8.7239$. The upcrossing benchmark is $\widehat u_F=1.7124<y$, while the dividend and liquidation policy favors continuation. This contrast reflects the additional value of future dividends in the joint problem. A genuine liquidation band emerges when the solvency threshold is increased to $y=3$. Over
$y=3,\ldots,6$, the lower liquidation boundary remains essentially
unchanged at $s_1^*=1.7124$,
whereas the upper boundary increases from $2.0166$ to $5.3035$.
The constancy of the lower boundary is consistent with
\eqref{eq:s1*}: once the lower recovery-liquidation boundary satisfies $s_1^*<y$, its defining equation depends on $r+q$ and the liquidation payoff but not directly on $y$.  By
contrast, a higher solvency threshold enlarges the region exposed to
the Omega clock and moves the upper boundary of the stopping band upward.
The dividend barrier increases monotonically from $8.7239$ at $y=2$
to $12.1236$ at $y=6$.  Thus, a more stringent solvency threshold
affects both margins of the policy, which creates and subsequently
expands the voluntary-liquidation region, while also inducing the
firm to retain a larger capital buffer before paying dividends.

Panel B reports the effect of the regulatory termination rate $q$.
All parameter values considered in this panel generate a genuine
liquidation band.  As $q$ increases from $0.4$ to $1.2$, the lower
boundary decreases from $2.2770$ to $1.4389$, while the upper
boundary increases from $2.3390$ to $3.4623$.  The stopping band
therefore widens substantially in both directions.  This behavior is
consistent with the increased cost of remaining in the distress
region, i.e., a larger termination intensity lowers the continuation value
below the solvency threshold and makes voluntary liquidation
attractive over a wider range of surplus levels.  In comparison, the
dividend barrier responds only moderately, increasing from $10.1648$
to $10.5054$.  Hence, in this benchmark, the termination intensity
acts primarily through the liquidation decision rather than through
the dividend barriers.

Panel C displays the effect of the profitability drift $c$.  For
$c=0.8,\ldots,1.4$, the liquidation band shrinks  as the
drift increases; in particular, the lower boundary rises from $1.5615$ to $1.7883$,
whereas the upper boundary decreases from $3.9429$ to $2.4929$.
At $c=1.6$, voluntary liquidation disappears altogether and the
selected regime becomes no liquidation, with $b^*=10.2314$.
This transition is consistent with the increasing value of recovery
when the surplus has a stronger positive drift. In addition, the dividend barrier
is non-monotone over the full range: it increases from $7.6970$ at
$c=0.8$ to approximately $10.4458$ at $c=1.4$, before falling
slightly in the no-liquidation regime.  Thus profitability affects
the dividend decision through both the improved continuation value
and the endogenous change in the optimal liquidation regime.

Panel D shows a particularly clear effect of diffusion volatility.
As $\sigma$ increases from $1.5$ to $2.5$, the dividend barrier
rises from $8.5966$ to $11.8480$, indicating a substantially larger
precautionary reserve under greater surplus uncertainty.  Both
liquidation boundaries also move upward: $s_1^*$ increases from
$1.5519$ to $1.8676$ and $s_2^*$ from $3.0074$ to $3.3522$.
The width of the stopping band changes comparatively limited, so the
main liquidation effect of higher volatility in this calibration is
an upward shift of the region rather than a large expansion of its
width.

\begin{table}[ht]
  \centering
  \caption{Sensitivity of candidate policy regimes, boundaries, and the upcrossing benchmark for the Brownian model.}
  \begin{tabular}{l c cccc}
    \toprule
    Regime & Value & $\us^*$ & $\bs^*$ & $b^*$ & $\widehat u_F$\\
    \midrule
    \multicolumn{6}{l}{\textit{Panel A: Regulatory Solvency Threshold ($y$)}} \\
   No Liq. &2.00  & -- & -- & 8.7239 & 1.7124 \\
    Band & 3.00  & 1.7124 & 2.0166 & 9.6589  & 1.7124\\
    Band & 4.00  & 1.7124 & 3.1936 & 10.4137& 1.7124  \\
   Band  &5.00  & 1.7124 & 4.2669 & 11.2397 & 1.7124 \\
    Band  &6.00  & 1.7124 & 5.3035 & 12.1236 & 1.7124 \\
    \midrule
    \multicolumn{6}{l}{\textit{Panel B: Regulatory Termination Intensity ($q$)}} \\
    Band &  0.40  & 2.2770 & 2.3390 & 10.1648 & 2.2770 \\
  Band   & 0.60  & 1.9323 & 2.9204 & 10.3268 & 1.9323 \\
   Band  & 0.80  & 1.7124 & 3.1936 & 10.4137 & 1.7124 \\
   Band  & 1.00  & 1.5566 & 3.3551 & 10.4681 & 1.5566 \\
   Band  & 1.20  & 1.4389 & 3.4623 & 10.5054  & 1.4389\\
    \midrule
    \multicolumn{6}{l}{\textit{Panel C: Profitability Drift ($c$)}} \\
  Band  & 0.80 & 1.5615 & 3.9429 & 7.6970 & 1.5615 \\
   Band & 1.00  & 1.6366 & 3.6835 & 9.7492 & 1.6366  \\
   Band & 1.20  & 1.7124 & 3.1936 & 10.4137 & 1.7124 \\
  Band  & 1.40  & 1.7883 & 2.4929 & 10.4458 & 1.7883 \\
    No Liq. & 1.60  & -- & -- & 10.2314 & 1.8641 \\
    \midrule
    \multicolumn{6}{l}{\textit{Panel D: Volatility ($\sigma$)}} \\
   Band & 1.50  & 1.5519 & 3.0074 & 8.5966 & 1.5519  \\
  Band  & 1.75  & 1.6325 & 3.1017 & 9.5466  & 1.6325\\
  Band  & 2.00  & 1.7124 & 3.1936 & 10.4137 & 1.7124 \\
  Band &  2.25  & 1.7909 & 3.2778 & 11.1819 & 1.7909 \\
  Band  & 2.50  & 1.8676 & 3.3522 & 11.8480  & 1.8676\\
    \bottomrule
  \end{tabular}
  \label{tab:sens_bm}
\end{table}

\paragraph{Cram\'er--Lundberg model.}
Table \ref{tab:sens_cl} exhibits more pronounced regime changes.
The distinction is consistent with the path structure of the
Cram\'er--Lundberg surplus (i.e., between claim arrivals the process moves
deterministically upward, while deterioration occurs through
downward jumps).  Consequently, the value of waiting for recovery and
the value of voluntary liquidation can change sharply as the
regulatory or claim parameters vary.

Panel A shows the effect of the regulatory solvency threshold.  For
$y=2$ and $y=3$, the selected policy has no voluntary
liquidation, with dividend barriers $8.1580$ and $9.1580$,
respectively.  A single-band for liquidation appears at $y=4$ and
remains optimal for $y=5$ and $y=6$.  Over these three cases, $s_1^*=1.5928$
remains unchanged, while $s_2^*$ increases from $2.2898$ to
$4.7850$.  As in the Brownian model, the stability of the lower
boundary is consistent with its recovery-liquidation characterization,
whereas the upper stopping boundary responds directly to the extent
of the regulatory distress region.  The dividend barrier rises from
$8.1580$ at $y=2$ to $11.8692$ at $y=6$.  Thus increasing the
solvency threshold first induces a transition from pure continuation
to voluntary liquidation and subsequently expands the liquidation
band while raising the capital retained before dividend payments.

At $y=2$, the upcrossing benchmark $\widehat u_F=5.0720$ lies above $y$, while the selected policy has no voluntary liquidation. As in the Brownian case, future dividend opportunities change the comparison between liquidation and continuation.

Panel B gives a similar regime transition as the regulatory
termination intensity increases.  When $q=0.4$ and $q=0.6$, no
voluntary liquidation is optimal.  A single-band case appears
at $q=0.8$.  Increasing $q$ further to $1.2$ lowers $s_1^*$ from
$1.5928$ to $1.2673$ and raises $s_2^*$ from $2.2898$ to $3.0070$.
The liquidation band therefore widens as regulatory termination
becomes more severe.  By comparison, the dividend barrier changes
only modestly, from $9.5007$ at $q=0.4$ to $10.3037$ at $q=1.2$. Such observations are consistent with the Brownian model, where the Omega termination intensity has a stronger effect on voluntary liquidation than on the dividend barrier.

Panel C reveals the most pronounced regime dependence.  At $c=2.0$, the optimal stopping band is relatively wide,
$[1.3272,3.7552]$, whereas at $c=2.5$ it contracts to $[1.5928,2.2898]$.
When the premium rate is increased further to $c=2.7$ and $c=3.0$, the selected policy involves no voluntary liquidation.
Thus, over the range where the present verification applies, a
larger premium rate increases the value of waiting for recovery and eventually eliminates the incentive for voluntary liquidation.

At $c=1.8$, the positive safety-loading condition $c>\lambda/\alpha$ holds, but none of the single-band, singleton, or no-liquidation candidates passes the sufficient-condition checks. This case is therefore reported as ``Inconclusive''. Its optimal stopping geometry remains open within the present regular barrier candidate framework.

Panel D shows the effect of claim-arrival intensity.  For
$\lambda=0.5$ and $\lambda=0.75$, the optimal stopping set is empty.
A liquidation band emerges at $\lambda=1$ and expands as claim arrivals become more frequent.  In particular, increasing
$\lambda$ from $1$ to $1.5$ lowers the lower boundary from $1.5928$ to $1.4246$ and raises the upper boundary from $2.2898$ to $3.8248$.
More frequent claims therefore make voluntary liquidation attractive over a substantially larger set of surplus levels.  The dividend
barrier is non-monotone: it rises from $8.9710$ at $\lambda=0.5$ to $10.2154$ at $\lambda=1.25$, but then falls to $9.3287$ at
$\lambda=1.5$.  At moderate claim intensities, the firm responds by retaining a larger capital buffer, whereas at sufficiently high claim frequency, the deterioration in continuation value makes earlier
dividend distribution, together with a wider liquidation region, more attractive.

\begin{table}[ht]
  \centering
  \caption{Sensitivity of candidate policy regimes, boundaries, and the upcrossing benchmark for the Cram\'er-Lundberg model.}
  \begin{tabular}{lccccc}
    \toprule
    Regime & Value & $\us^*$ & $\bs^*$ & $b^*$  & $\widehat u_F$ \\
    \midrule
    \multicolumn{6}{l}{\textit{Panel A: Regulatory Solvency Threshold ($y$)}} \\
  No Liq.  & 2.00 &  -- & -- & 8.1580 & 5.0720  \\
   No Liq. & 3.00 &  -- & -- & 9.1580 & 1.5928\\
  Band & 4.00 &  1.5928 & 2.2898 & 10.1325& 1.5928  \\
   Band& 5.00 &  1.5928 & 3.6710 & 10.9891 & 1.5928 \\
  Band & 6.00 & 1.5928 & 4.7850 & 11.8692 & 1.5928 \\
    \midrule
    \multicolumn{6}{l}{\textit{Panel B: Regulatory Termination Intensity ($q$)}} \\
  No Liq.   & 0.40 &  -- & -- & 9.5007  & 2.2662 \\
   No Liq.  &0.60 &  -- & -- & 9.9037  & 1.8558\\
   Band &0.80 & 1.5928 & 2.2898 & 10.1325& 1.5928\\
    Band &1.00 &  1.4070 & 2.7636 & 10.2366 & 1.4070 \\
    Band &1.20 &  1.2673 & 3.0070 & 10.3037 & 1.2673 \\
    \midrule
    \multicolumn{6}{l}{\textit{Panel C: Premium Drift ($c$)}} \\
   Inconclusive &1.80  & -- &   -- & -- & 1.2145  \\
  Band  &2.00 & 1.3272 & 3.7552 & 8.7835 & 1.3272 \\
   Band &2.50 &  1.5928 & 2.2898 & 10.1325 & 1.5928 \\
  No Liq.  &2.70 &  -- & -- & 10.2097 & 1.6925 \\
  No Liq.  &3.00 &  -- & -- & 10.2455 & 1.8354\\
    \midrule
    \multicolumn{6}{l}{\textit{Panel D: Jump Intensity ($\lambda$)}} \\
  No Liq.  & 0.50 &  -- & -- & 8.9710 & 1.7687 \\
  No Liq. & 0.75 &  -- & -- & 9.6261  & 1.6802\\
   Band & 1.00 &  1.5928 & 2.2898 & 10.1325 & 1.5928 \\
   Band& 1.25 &  1.5074 & 3.2422 & 10.2154  & 1.5074\\
   Band & 1.50 &  1.4246 & 3.8248 & 9.3287 & 1.4246 \\
    \bottomrule
  \end{tabular}
  \label{tab:sens_cl}
\end{table}

\paragraph{Singleton liquidation regimes.}
No singleton liquidation regime appears at the parameter values in Tables \ref{tab:sens_bm}--\ref{tab:sens_cl}. In these cases, singleton liquidation occurs at the transition between no liquidation and a single band.

For example, keeping the Brownian model benchmark parameters fixed and varying only the solvency threshold, the transition occurs at
approximately
$y=2.7742$, with $a=1.7124, b=9.4981$. The transition is characterized by $\Delta_0=0$ given in \eqref{eq:Delta0}. There, the no-liquidation value touches $g$ at $a$ and agrees with the singleton candidate, giving alternative policies under Theorems~\ref{thm:singleton} and \ref{thm:no-liquidation}.  Similarly, in the Cram\'er--Lundberg benchmark model, varying only $q$ gives a singleton scenario at approximately $q=0.7257$, with $a=1.6790,  b=10.0748$, where the same obstacle-gap equation gives the CL singleton continuous-fit condition.

\section{Conclusion}\label{sec:conclude}
This paper studies joint dividend control and voluntary liquidation for a firm with spectrally negative L\'evy surplus and Omega-type regulatory termination. The regulatory threshold is exogenous, while shareholders choose voluntary liquidation and receive a state-dependent payoff. Recovery below the distress threshold, liquidation, and regulatory termination consequently play distinct roles.
The direct Meyer--It\^o argument establishes global optimality for attainable candidates satisfying the verification conditions over the full admissible class. Reflection and liquidation at the obstacle contact set attain every regular barrier stochastic solution, yielding uniqueness and identification with both the auxiliary value function $V^M$ and the value function $V$ of the original control--stopping problem.
Omega-scale functions provide semi-analytical values for single-band, singleton, and no-liquidation candidates. The regular-class existence criterion connects the verified candidates to the uniqueness theorem. In the Brownian and exponential-claim models, interior barrier stationarity implies the strict dividend inequality, while an explicit Brownian regime gives existence under primitive parameter conditions.
The benchmark computations reduce verification to exponential coefficient signs and scalar extrema, with exact residual formulas for the affine tail. Higher distress thresholds raise retained capital and can induce or enlarge a liquidation band. Greater profitability can remove voluntary liquidation, while higher Brownian volatility raises the dividend barrier in the reported range. More frequent CL claims can widen liquidation and produce a non-monotone dividend barrier. The comparisons reflect the stated calibrations and their different mean and variance rates.
The CL calibration $c=1.8$ remains inconclusive within this candidate framework. The analysis provides a tractable construction and global verification criteria for studying dividend distribution, voluntary liquidation, and regulatory distress jointly.

%
%\section*{Statements and Declarations}
%No competing interests.
%%%%%%%%%%%%%%%%%%%%%%%%%%%%%%%%%%%%%%%%%%%%%%%%%%%%%%%%%%%%%%%%%
%\section*{Acknowledgments}
%Ran Xu was supported by the National Natural Science Foundation of China [No.12201506 and 12371468].
%Zbigniew Palmowski was partially supported by the National Science Centre (Poland) under the grant 2021/41/B/HS4/00599.
%%%%%%%%%%%%%%%%%%%%%%%%%%%%%%%%%%%%%%%%%%%%%%%%%%%%%%%%%%%%%%%%%

\begin{appendices}

\makeatletter
\newcommand{\prooflocation}[1]{%
  \leavevmode

    \phantomsection
    \protected@edef\@currentlabel{\thesection}%
    \label{#1}%

  \ignorespaces
}
\makeatother

\section{Preliminaries and dynamic programming}
\label{app:preliminaries-dpp}

\begin{proof}[Proof of Proposition~\ref{prop:valuefunction}]
\prooflocation{App.1}
We divide the proof into six steps.
\paragraph{Step 1.} We first show that the negative voluntary liquidation payoffs can be excluded.
For any admissible pair $(D,\tau)$, define
$ B:=\{\tau<\infty,\ g(U_\tau)\ge0\}$, and $\widehat\tau:=\tau$ on $B$ and $\widehat\tau:=\infty$ on $B^c$.
Since $B\in\mathcal{F}_\tau$, $\widehat{\tau}$ is a $\mathbb{F}$-stopping time. On the event $\{\tau<\infty, g(U_\tau)<0\}$, the modified stopping strategy $\widehat{\tau}$ simply refrains from voluntary liquidation. Since dividend payments are nonnegative and regulatory termination yields zero payoff, continuing instead of exercising a negative voluntary-liquidation payoff cannot decrease the performance function. Hence
\begin{equation*}
%\label{prop3.3-1}
    J(x;y;D,\widehat{\tau})\ge J(x;y; D, \tau).
\end{equation*}
Consequently, without changing the value function, the supremum in the value function may be restricted to admissible strategies satisfying
\begin{equation}
\label{prop3.3-2}
     g(U_\tau)\ge0 \quad\text{on }\{\tau<\infty\}
\end{equation}
 Then, taking $D\equiv0$ and $\tau=\infty$ gives
\begin{equation}\label{prop3.3-3}
   V(x;y)\ge0
\end{equation}

 \paragraph{Step 2.}  Let $\overline{X}_t:=\sup_{0\le s\le t}X_s.$
 By the admissibility definition of $D$, dividends can increase only when the pre-dividend surplus is at least $y$, and a lump-sum dividend payment cannot reduce the surplus below $y$. Therefore
 \begin{equation}\label{prop3.3-4}
     D_t\le(\overline{X}_t-y)^+\le \overline{X}_t^+.
 \end{equation}
 Since $A^{D,y}_t\ge rt$,  integration by parts  and \eqref{prop3.3-4} yields
\begin{align*}
%\label{prop3.3-5}
   \int_{[0,t]}e^{-A_s^{D,y}}\,\mathrm{d}D_s &\le \int_{[0,t]}e^{-rs}\,\mathrm{d}D_s \nonumber \\
   &=e^{-rt}D_t+r\int_0^t e^{-rs}D_s\,\mathrm{d}s \nonumber\le r\int_0^\infty e^{-rs}\overline{X}_s^+\,\mathrm{d}s.
\end{align*}
Let $\mathbf{e}_r$ be an exponential random variable with parameter $r$, independent of $X$. Then
\begin{equation*}
%\label{prop3.3-6}
   \mathbb{E}_x\!\left[ r\int_0^\infty e^{-rs}\overline{X}_s^+\,\mathrm{d}s \right] =\mathbb{E}_x[\overline{X}_{\mathbf{e}_r}^+] \le x^++\frac{1}{\Phi(r)},
\end{equation*}
 where in the last inequality, we have applied the fact that for a spectrally negative L\'evy process, $\overline{X}_{\mathbf{e}_r}-x$ is exponentially distributed with parameter $\Phi(r)$. It follows that, uniformly over $D\in\mathcal{D}(y)$, we have
\begin{equation}\label{prop3.3-7}
\mathbb{E}_x\Bigg[ \int_{[0,t]}e^{-A_s^{D,y}}\,\mathrm{d}D_s\Bigg]\le x^+ + \frac{1}{\Phi(r)}.
\end{equation}
Moreover, by Assumption (\textbf{A.1}), there exists constants $c_0, c_1>0$ such that $g^+(z)\le c_0 + c_1 z^+$ for $z\in\mathbb{R}$. Since $D_t\ge 0$, $U_t = X_t-D_t \le X_t\le \overline{X}_t$, and $U^+_t\le \overline{X}^+_t$. Then, for any $\mathbb{F}$-stopping $\tau$, we have
\begin{align*}
    e^{-r\tau} \overline{X}^+_\tau \le \int_{\tau}^{\infty}re^{-rs}\overline{X}^+_s\dif s \le \int_{0}^{\infty}re^{-rs}\overline{X}^+_s\dif s,
\end{align*}
because $s\mapsto \overline{X}^+_s$ is non-decreasing. Thus
\begin{align}\label{prop3.3-8}
    \mathbb{E}_x\!\left[e^{-A_\tau^{D,y}}g^+(U_\tau)\right] &\le c_0+ c_1\mathbb{E}_x\!\left[e^{-r\tau}\overline{X}_\tau^+\right] \nonumber \\
    &\le c_0+ c_1\mathbb{E}_x\!\left[ r\int_0^\infty e^{-rs}\overline{X}_s^+\,\mathrm{d}s \right] \nonumber\\ &\le c_0+c_1\left(x^++\frac{1}{\Phi(r)}\right).
\end{align}
 Combining \eqref{prop3.3-7} and \eqref{prop3.3-8}, and by Step 1 (it is sufficient to consider stopping times for nonnegative liquidation payoff), we have
\begin{equation*}
%\label{prop3.3-9}
     0\le V(x;y) \le c_0+(1+c_1)\left(x^++\frac{1}{\Phi(r)}\right),
\end{equation*}
which proves finiteness, the assertion \textit{(iii)}, and \eqref{V-bound}.

\paragraph{Step 3.} Now, we prove the monotonicity with respect to the initial surplus. For every fixed $y>0$, fix $x_1>x_2$ and $\varepsilon>0$. By Step 1, choose an $\varepsilon$-optimal strategy $(D^\varepsilon,\tau^\varepsilon)\in (\mathcal{D}(y),\mathcal{T})$ for initial surplus $x_2$ such that $g(U^{(2)}_{\tau^\varepsilon})\ge 0$ on $\{\tau^\varepsilon<\infty\}$, and
\[
V(x_2;y)-\varepsilon \le J(x_2;y;D^\varepsilon,\tau^\varepsilon),
\]
where $U^{(i)}$ denotes the controlled surplus process with initial surplus $x_i$. Applying the same dividend process and stopping time to the process starting from $x_1$ gives
\[
U_t^{(1)}=U_t^{(2)}+(x_1-x_2).
\]
In addition, we have $\mathds{1}_{\{U^{(1)}<y\}}\le \mathds{1}_{\{U^{(2)}<y\}}$, hence, $A_t^{D^\varepsilon,y,(1)}\le A_t^{D^\varepsilon,y,(2)}$, where each functional is computed along $U^{(i)}$ and $g(U_{\tau^\varepsilon}^{(1)}) \ge g(U_{\tau^\varepsilon}^{(2)})\ge0$. Thus,
\[
V(x_1;y) \ge J(x_1;y;D^\varepsilon,\tau^\varepsilon) \ge J(x_2;y;D^\varepsilon,\tau^\varepsilon) \ge V(x_2;y)-\varepsilon.
\]
Letting $\varepsilon\downarrow0$ proves that $V(\cdot;y)$ is increasing on $\mathbb{R}$.  If $x_1>x_2\ge y$, and the amount $x_1-x_2$ is paid immediately as a lump-sum dividend (which is admissible), after which the firm starts from $x_2$. Therefore
\[
V(x_1;y)\ge x_1-x_2+V(x_2;y)-\varepsilon.
\]
If an $\varepsilon$-optimal continuation liquidates at time zero, first delay that liquidation to a deterministic $\eta>0$, make no additional dividend payment during the delay, and let $\eta\downarrow0$. Right continuity of the surplus and continuity of $g$, together with Lemma~\ref{lem:common-envelope}, make the loss tend to zero. The initial dividend is then included in the payoff, proving \eqref{eq:lowerbound-V}.

\paragraph{Step 4.} We next show the continuity with respect to the initial surplus. Fix $x<a$, and let
 \[
 T_a^+:=\inf\{t\ge0:X_t\ge a\}.
 \]
Since $X$ has no positive jumps,  we have $X_{T^+_a}=a$ on $\{T^+_a <\infty\}$. Consider the admissible strategy that pays no dividend and does not liquidate before $T^+_a$ and, conditional on reaching $a$, subsequently uses an $\varepsilon$-optimal strategy for initial surplus $a$.
Upward exit is continuous, so exit times to levels increasing to $a$ announce $T_a^+$ on $\{T_a^+<\infty\}$. The shifted control, including any initial payment at $T_a^+$, is predictable and satisfies the solvency constraint. The strong Markov property then gives

\begin{equation}
    \label{prop3.3-9}
     V(x;y)\ge \mathbb{E}_x\!\left[e^{-A_{T_a^+}^{0,y}}\mathds{1}_{\{T_a^+<\infty\}}\right] \bigl(V(a;y)-\varepsilon\bigr).
\end{equation}
By \eqref{eq:onesideexitH}
  \[
\mathbb{E}_x\!\left[e^{-A_{T_a^+}^{0,y}}\mathds{1}_{\{T_a^+<\infty\}}\right] = \frac{\mathcal I^{(r,q)}(x-y)}{\mathcal I^{(r,q)}(a-y)},
  \]
  which reduces to $e^{-\Phi(r+q)(a-x)}$ for $a\le y$. Letting $\varepsilon\downarrow 0$ in \eqref{prop3.3-9} yields
  \begin{equation}\label{prop3.3-10}
      0\le V(a;y)-V(x;y) \le \left( 1-\frac{\mathcal I^{(r,q)}(x-y)}{\mathcal I^{(r,q)}(a-y)} \right)V(a;y).
  \end{equation}
  The continuity of $\mathcal I^{(r,q)}$ and \textit{(iii)} imply the continuity of $V(x;y)$ for $x\in \mathbb{R}$.
\paragraph{Step 5.} We show the monotonicity with respect to the solvency threshold $y$. Fix $0<y_2<y_1$. Since the larger solvency threshold imposes a stricter dividend constraint, we have $\mathcal{D}(y_1)\subset \mathcal{D}(y_2)$. Let $\varepsilon>0$, and choose (by Step 1) an $\varepsilon$-optimal strategy $(D^\varepsilon, \tau^\varepsilon)\in \mathcal{D}(y_1)\times \mathcal{T}$ for $V(x;y_1)$ satisfying \eqref{prop3.3-2}. For the same control process associated with regulatory solvency level $y_2$, we have $\mathds{1}_{\{U_s<y_2\}}\le \mathds{1}_{\{U_s<y_1\}}$ and consequently $A^{y_2}_t\le A^{y_1}_t$ for all $t\ge 0$. Since both the dividend measure $\dif D^\varepsilon$ and the realized terminal liquidation payoff in \eqref{prop3.3-2} are nonnegative, we have
\begin{align*}
  V(x;y_1)-\varepsilon &\le J(x;y_1;D^\varepsilon,\tau^\varepsilon) \le J(x;y_2;D^\varepsilon,\tau^\varepsilon) \le V(x;y_2),
\end{align*}
where the last inequality holds since the strategy $(D^\varepsilon,\tau^\varepsilon)$ is admissible under $y_2$. Thus, by letting $\varepsilon\downarrow 0$ completes the proof of \textit{(ii)}.

\paragraph{Step 6.} Finally, we prove \textit{(iv)}. Set
\[
 Y_r:=r\int_0^\infty e^{-rs}\overline{X}_s^+\,\mathrm{d}s.
\]
 By \eqref{prop3.3-8}, $Y_r$ is integrable. From \eqref{prop3.3-7} we have $\int_{[0,t]}e^{-A_s^{D,y}}\,\mathrm{d}D_s\le Y_r$. Then, by \textit{(iii)} and $U^+_t\le \overline{X}^+_t$, we have
 \[
e^{-A_t^{D,y}}V(U_t;y) \le C e^{-rt}(1+\overline{X}_t^+) \le C(1+Y_r),
 \]
  for some $C>0$. Hence, the family in \eqref{VD-UI} is uniformly integrable, completing the proof of \textit{(iv)}.

\end{proof}

\begin{proof}[Proof of Corollary~\ref{cor:continuation-DPP}]

\prooflocation{App:A+1}
For the lower inequality, fix an admissible policy $D$ up to $h^D$. The assumed measurable approximate continuation selection and admissible concatenation yield, for every $\delta>0$, a payoff at least
\[
\mathbb E_x\left[\int_{[0,h^D]}e^{-A_s^{D,y}}\,\dif D_s
+e^{-A_{h^D}^{D,y}}\bigl(V^{\Mmark}(U_{h^D}^D;y)-\delta\bigr)\right].
\]
Taking the supremum and letting $\delta\downarrow0$ gives the lower inequality, with $D_0$ included in the integral. For the reverse inequality, use the continuation-set Bellman property in Assumption~\ref{ass:aux-dpp}. Before $h^D$, $V^{\Mmark}-g>\varepsilon$, so a $\delta$-optimal continuation with $\delta<\varepsilon$ dominates immediate liquidation; at $h^D$, $V^{\Mmark}\ge g$. The common stopping rule is compatible with this substitution, yielding \eqref{eq:continuation-DPP}.

\end{proof}

\begin{proof}[Proof of Proposition~\ref{prop:controlsuper}]
\prooflocation{app:proof-controlsuper}
Proposition~\ref{prop:Vm-supermartingale} gives
$V^{\Mmark}\in\mathscr L_I^{\varepsilon,y,\Mmark}$.
Fix $\phi\in\mathscr L_I^{\varepsilon,y,\Mmark}$ and
$x\in\operatorname{int}(I)\cap\mathscr N_\varepsilon^{\Mmark}$.
Put $h_t^D=t\wedge\rho^D$.
Optional sampling, with any initial dividend handled by
\eqref{eq:finite-dividend-ineq}, gives
\[
\phi(x)\ge
\mathbb E_x\left[
\int_{[0,h_t^D]}e^{-A_s^{D,y}}\,\dif D_s
+e^{-A_{h_t^D}^{D,y}}\phi(U_{h_t^D}^D)\right].
\]

On $\{\rho^D\le t\}$, the boundary condition gives
$\phi(U_{\rho^D}^D)\ge V^{\Mmark}(U_{\rho^D}^D;y)$.
On $\{\rho^D>t\}$, affine growth and
$|U_t^D|\le2X_t^*$ give
\[
\mathbb E_x\left[
e^{-A_t^{D,y}}
\left|\phi(U_t^D)-V^{\Mmark}(U_t^D;y)\right|
\mathds{1}_{\{\rho^D>t\}}\right]
\le C\mathbb E_x\left[e^{-rt}(1+X_t^*)\right],
\]
where $X_t^*=\sup_{s\le t}|X_s|$ and $C$ is independent
of $D$ and $t$.
Taking the supremum over $D$ and using the local DPP
\eqref{eq:continuation-DPP} therefore yields
\[
\phi(x)\ge V^{\Mmark}(x;y)
-C\mathbb E_x\left[e^{-rt}(1+X_t^*)\right].
\]
Letting $t\to\infty$ and applying
Lemma~\ref{lem:common-envelope} proves
$\phi(x)\ge V^{\Mmark}(x;y)$.
For the remaining initial states, this inequality follows
from the boundary condition.
Together with the membership of $V^{\Mmark}$, this proves
\eqref{eq:epsilon-min}.
\end{proof}

\section{Verification and regular barrier solutions}
\label{app:verification-barriers}

\begin{proof}[Proof of Theorem~\ref{prop:direct-verification}]

\prooflocation{App:A+2}
Fix an admissible pair $(D,\tau)$ and abbreviate
$U=U^D$ and $A=A^{D,y}$. Define
\[
H_t
:=
e^{-A_t}\varphi(U_t)
+\int_{[0,t]}e^{-A_s}\,\mathrm{d}D_s.
\]
Admissibility and condition (b) give
\[
H_0=D_0+\varphi(x-D_0)\le\varphi(x).
\]
Thus the initial dividend is included throughout the argument.

Put
\[
\delta_X
:=
\gamma+\int_{(0,1)}z\,\nu(\mathrm{d}z)>0,
\qquad
N_\varphi(u)
:=
\int_{(0,\infty)}
[\varphi(u-z)-\varphi(u)]\,\nu(\mathrm{d}z).
\]
Local Lipschitz regularity, affine growth, and the
finite-first-moment assumption imply that $N_\varphi$
is finite and continuous. On smooth pieces,
\[
\mathcal L\varphi
=
\delta_X\varphi'
+\frac{\sigma^2}{2}\varphi''
+N_\varphi.
\]

Since $D$ is predictable and $X$ is quasi-left-continuous,
they have no simultaneous nonzero jumps.
Apply the Meyer--It\^o formula when $\sigma>0$,
and the finite-variation chain rule when $\sigma=0$,
together with integration by parts. After localization, we have
\begin{align}
\label{eq:controlled-Meyer-Ito}
H_t-H_0
={}&M_t^{\mathrm{loc}}
+\int_0^t e^{-A_s}
(\mathcal L-\omega_y)\varphi(U_{s-})\,\mathrm{d}s
\nonumber\\
&+\int_{(0,t]}e^{-A_s}
[1-\varphi'_-(U_s)]\,\mathrm{d}D_s^c
\nonumber\\
&+\sum_{\substack{0<s\le t\\ \Delta D_s>0}}
e^{-A_s}
[\varphi(U_s)-\varphi(U_{s-})+\Delta D_s]
\nonumber\\
&+\frac12\sum_{a\in\mathcal B}
[\varphi'(a+)-\varphi'(a-)]
\int_0^t e^{-A_s}\,\mathrm{d}L_s^a(U).
\end{align}
Here $M^{\mathrm{loc}}$ is a local martingale and
$\varphi'_-$ denotes the spatial left derivative.
The last line is absent when $\sigma=0$.
In that case, at $a\in\mathcal B$, set
$\mathcal L\varphi(a)=\delta_X\varphi'_-(a)+N_\varphi(a)$;
when $\sigma>0$, its values on $\mathcal B$ may be
chosen arbitrarily.

At a dividend jump, admissibility gives
$U_s\ge y$ and $U_{s-}=U_s+\Delta D_s$.
Hence condition (b) makes every summand in
\eqref{eq:controlled-Meyer-Ito} nonpositive.
We next consider the continuous terms.

\medskip
\noindent
\emph{The case $\sigma>0$.}
Condition (b) gives $\varphi'_-(u)\ge1$ for $u>y$
and $\varphi'(y+)\ge1$.
Condition (d) then implies $\varphi'_-(y)\ge1$.
Since continuous dividends are supported on $\{U\ge y\}$,
their contribution is nonpositive.
The occupation-density formula implies that $U$ spends
zero Lebesgue time in the finite set $\mathcal B$.
Thus the generator term is nonpositive by (c),
and all local-time terms are nonpositive by (d).

\medskip
\noindent
\emph{The case $\sigma=0$.}
For $u>y$ outside $\mathcal B$, conditions (b) and (c)
give
\[
\delta_X+N_\varphi(u)-r\varphi(u)
\le
\delta_X\varphi'(u)+N_\varphi(u)-r\varphi(u)
\le0.
\]
By continuity, taking right limits yields
\begin{equation}
\label{eq:verification-BV-boundary}
\delta_X+N_\varphi(a)-r\varphi(a)\le0,
\qquad a\in\mathcal B\cap[y,\infty).
\end{equation}

For a finite-variation path, its continuous Stieltjes
measure vanishes on each fixed level in set $\mathcal B$.
Applying this fact to
$\mathrm{d}U_s^c=\delta_X\,\mathrm{d}s-\mathrm{d}D_s^c$
gives
\begin{equation}
\label{eq:BV-level-set}
\mathbf{1}_{\{U_s=a\}}\,\mathrm{d}D_s^c
=
\delta_X\mathbf{1}_{\{U_s=a\}}\,\mathrm{d}s.
\end{equation}
If $a<y$, admissibility makes the left-hand side zero,
so the occupation time at $a$ is zero.
If $a\in\mathcal B\cap[y,\infty)$, the generator and
continuous-dividend terms on this level combine into
\[
e^{-A_s}\mathbf{1}_{\{U_s=a\}}
[\delta_X+N_\varphi(a)-r\varphi(a)]\,\mathrm{d}s,
\]
which is nonpositive by
\eqref{eq:verification-BV-boundary}.
Away from $\mathcal B$, their signs follow directly
from (b), (c), and admissibility.

\medskip
Consequently, $H$ is a local supermartingale.
Let $(\rho_n)$ localize both the state and the local
martingale. Optional sampling gives
\[
\mathbb E_x[H_{t\wedge\tau\wedge\rho_n}]
\le H_0\le\varphi(x).
\]
The common integrable bound in
Lemma~\ref{lem:common-envelope} allows us to let
$n\to\infty$ at fixed $t$, obtaining
\[
\mathbb E_x[H_{t\wedge\tau}]\le\varphi(x).
\]

Finally, on $\{\tau<\infty\}$, condition (a) and
$\Delta D_\tau\ge0$ imply
\[
H_\tau
\ge
\int_{[0,\tau)}e^{-A_s}\,\mathrm{d}D_s
+e^{-A_\tau}g(U_\tau).
\]
By Lemma~\ref{lem:common-envelope}, the residual term
$e^{-A_t}\varphi(U_t)\mathbf{1}_{\{\tau>t\}}$
converges to zero in $L^1$.
Dominated convergence handles the finite terminal payoff,
and monotone convergence handles the discounted dividends.
Letting $t\to\infty$ therefore gives
\[
\varphi(x)\ge J(x;y;D,\tau),
\]
with zero terminal payoff when $\tau=\infty$.
The same argument includes $\tau=0$, since
$H_0\ge g(U_0)$.
Taking the supremum over $(D,\tau)$ proves
\eqref{eq:phi_ge_V}.
\end{proof}

\begin{proof}[Proof of Lemma~\ref{lem:Ymartingale}]
\prooflocation{App:B}
Apply \eqref{eq:controlled-Meyer-Ito} with $\varphi=h$.
For (i), take $D\equiv0$ and stop at
$\theta=t\wedge T_I$.
For (ii), take $D=D^b$ and stop at
$\theta=t\wedge\tau_I^b$, extending $h$ above $b$
by $h(b)+x-b$.
Since $x\le b$, we have $D_0^b=0$.
Moreover, $D^b$ is continuous on $(0,\infty)$,
so the dividend-jump sum in
\eqref{eq:controlled-Meyer-Ito} is absent.
Interior breakpoints have zero Lebesgue occupation time, by the
occupation-density formula when $\sigma>0$. When $\sigma=0$, \eqref{eq:BV-level-set} gives zero Lebesgue occupation time at every interior breakpoint, since
$\delta_X>0$ and no continuous dividends are paid there.
Thus, \eqref{eq:h-harmonic} eliminates the interior drift terms,
while, when $\sigma>0$, the continuity of $h'$ at the interior
breakpoints eliminates the corresponding local-time terms.. This proves that the process in (i)
is a local martingale.
For (ii), adding the discounted dividends gives the reflection
term
\[
\int_{(0,\theta]}e^{-A_s^{b,y}}
\bigl(1-h'(b-)\bigr)\,\dif D_s^b=0.
\]
At $b$, the affine extension satisfies
$h'(b+)=h'(b-)=1$, so there is no local-time
contribution at the reflecting boundary.
The occupation density eliminates the drift term when
$\sigma>0$, and when $\sigma=0$, the boundary drift vanishes by
\eqref{eq:barrier-BV-harmonic}. Hence the process in (ii)
is also a local martingale.
Finally, choose $C$ such that $|h(u)|\le C(1+|u|)$.
Lemma~\ref{lem:common-envelope} bounds the absolute values of
both processes, at all finite stopped times, by the integrable
random variable $(1+2C)K_r$. The martingale localizations
therefore converge in $L^1$, and both processes are uniformly
integrable martingales.
\end{proof}

\begin{proof}[Proof of Proposition~\ref{prop:barrier-attain}]

\prooflocation{App:C+1}
Fix $x\le b$ and write $\tau=\tau_h^b$.
If $x\in\mathscr S_h$, then $\tau=0$, $D_0^b=0$,
and the payoff is $g(x)=h(x)$.
Assume henceforth that $x\notin\mathscr S_h$.
The contact set is closed, so right continuity gives
$Y_\tau^b\in\mathscr S_h$ on $\{\tau<\infty\}$.
Before $\tau$, the reflected surplus lies in
$\mathscr K_{h,b}\cup\{b\}$.

Set $\rho_n:=\inf\{s\ge0:Y_s^b\le-n\}$. Apply \eqref{eq:controlled-Meyer-Ito} with
$\varphi=h$, $D=D^b$, and stopping time
$\theta=t\wedge\tau\wedge\rho_n$.
Below $b$, harmonicity in
Definition~\ref{def:regular-barrier-solution}(iii)
removes the generator term away from the breakpoints.
When $\sigma>0$, the matching first derivatives remove
the local-time terms, and the occupation-density formula
gives zero Lebesgue occupation time at the breakpoints.
When $\sigma=0$, the same zero-occupation conclusion
below $b$ follows from \eqref{eq:BV-level-set}, since
$D^b$ increases only at $b$ and $\delta_X>0$.

At $b$, the condition $h'(b-)=1$ cancels the
reflection and dividend terms, and the affine extension
matches the first derivative there.
The generator contribution at $b$ vanishes by zero
occupation time when $\sigma>0$, and by
\eqref{eq:barrier-BV-harmonic} when $\sigma=0$.
Thus the stopped discounted value plus accumulated
dividends is a local martingale.
Lemma~\ref{lem:common-envelope} provides an integrable
bound uniform over both the state and local-martingale
localizations.
Removing these localizations gives
\begin{equation}\label{eq:barrier-attain-local}
h(x)=\mathbb E_x\left[
\int_{[0,t\wedge\tau]}e^{-A_s^{b,y}}\,\dif D_s^b
+e^{-A_{t\wedge\tau}^{b,y}}h(Y_{t\wedge\tau}^b)
\right].
\end{equation}

Since $h(Y_\tau^b)=g(Y_\tau^b)$ when $\tau<\infty$,
the terminal term decomposes as
\[
e^{-A_{t\wedge\tau}^{b,y}}h(Y_{t\wedge\tau}^b)
=
e^{-A_\tau^{b,y}}g(Y_\tau^b)
\mathds{1}_{\{\tau\le t\}}
+
e^{-A_t^{b,y}}h(Y_t^b)
\mathds{1}_{\{\tau>t\}}.
\]
The second term tends to zero in $L^1$ by the final
assertion of Lemma~\ref{lem:common-envelope}.
The common integrable bound in that lemma also gives
dominated convergence for the first term.
The dividend integral converges monotonically.
Moreover, $D^b$ is continuous on $(0,\infty)$ because
$X$ has no positive jumps, and $D_0^b=0$ for $x\le b$.
Hence the dividend integrals over $[0,\tau]$ and
$[0,\tau)$ agree.
Letting $t\to\infty$ in \eqref{eq:barrier-attain-local}
and using \eqref{eq:perform} therefore gives
$h(x)=J(x;y;D^b,\tau)$.

For $x>b$, the initial dividend $x-b$ leaves surplus $b$.
The condition $\mathscr S_h\subset(-\infty,b)$ gives
$h(b)>g(b)$.
Continuity of $h-g$ and right continuity of the reflected
process ensure that liquidation occurs strictly after
time zero when starting from $b$.
Thus the initial dividend is included in the payoff, and
\[
J(x;y;D^b,\tau_h^b)
=x-b+J(b;y;D^b,\tau_h^b)
=x-b+h(b)=h(x),
\]
where the last equality uses the affine extension in
Definition~\ref{def:regular-barrier-solution}(i).
This proves \eqref{eq:h-performance}.
\end{proof}

\begin{proof}[Proof of Proposition~\ref{prop:candidate-regular}]
\prooflocation{app:proof-candidate-regular}
The candidate verification proofs give a finite breakpoint set
$\mathcal B_v$ for which the hypotheses of
Theorem~\ref{prop:direct-verification} hold.
They also establish $v\in C^1(\mathbb R)$ when $\sigma>0$.
A similar supermartingale calculation as in the proof of Theorem~\ref{prop:direct-verification},
together with Lemma~\ref{lem:common-envelope}, shows that $v$
is a controlled stochastic supersolution.

Since $x_g\le b$ and $v$ is affine with slope one above $b$,
$v(x;y)-g(x)\ge v(b;y)-g(b)>0$ for $x\ge b$.
Thus $\mathscr S_v\subset(-\infty,b)$.
The candidate normalization also gives $v_x(b-;y)=1$.

The finite-dividend inequality implies $\Delta_yv\ge1$.
Fix $x\in(y,b)\setminus\mathscr S_v$; then $x>m$.
By \eqref{eq:regular-candidate-slack} and continuity of $v_x$,
choose $0<\eta<x-m$ and $\varepsilon>0$ such that
$v_x\ge1+\varepsilon$ on $[x-\eta,x]$.
Combining this bound with the finite-dividend inequality gives
\[
\frac{v(x;y)-v(z;y)}{x-z}
\ge 1+\varepsilon\min\!\left\{1,\frac{\eta}{x-z}\right\}
\ge 1+\frac{\varepsilon\eta}{x-y}>1,
\qquad y\le z<x.
\]
Hence $\Delta_yv(x)>1$.
If $y\notin\mathscr S_v$, the assumption $v_x(y+;y)>1$
gives a uniform derivative bound above one on a right
neighborhood of $y$, and therefore $\Delta_yv(y)>1$.
The unit left derivative at $b$ and the affine expression
above $b$ give $\Delta_yv=1$ on $[b,\infty)$.
Consequently,
\[
\mathscr C_v^{\mathrm{str}}
=\mathscr K_{v,b}
=(-\infty,b)\setminus\mathscr S_v.
\]

The candidate is harmonic on
$\mathscr K_{v,b}\setminus\mathcal B_v$, since this set avoids
the prescribed liquidation region.
For a compact interval $I\subset\mathscr K_{v,b}$, choose
a compact interval $J$ such that
$I\subset\operatorname{int}(J)\subset J\subset\mathscr K_{v,b}$.
Lemma~\ref{lem:Ymartingale}(i), applied on $J$ and then stopped
at the exit time of $I$, gives the martingale condition in
Definition~\ref{def:stochastic_solution}.
Thus $v$ is a stochastic solution.

When $\sigma=0$, taking the left limit of the harmonic
equation at $b$, using continuity of the jump integral and
$v_x(b-;y)=1$, gives \eqref{eq:barrier-BV-harmonic}.
Thus $v$ satisfies
Definition~\ref{def:regular-barrier-solution}.
The uniqueness assertion and optimality of
$(D^b,\tau_v^b)$ follow from
Theorem~\ref{thm:barrier-uniqueness}.
\end{proof}

\section{Scale-function construction and boundary selection}
\label{app:scale-boundaries}

\begin{proof}[Proof of Lemma~\ref{lem:step-omega-reflected}]

\prooflocation{App:C}
Set
\[
\chi_n(u):=\begin{cases}1,&u\le-1/n,\\-nu,&-1/n<u<0,\\0,&u\ge0,\end{cases}
\qquad\omega_n(x):=r+q\chi_n(x-y).
\]
The rates are continuous, satisfy $r\le\omega_n\le r+q$, and converge pointwise to $\omega_y$, with $\omega_n(y)=\omega_y(y)=r$. Dominated convergence gives almost-sure convergence of the integrated rates along $Y^b$ at each finite time. Discounted exit expectations and bounded-function resolvents are dominated by $1$ and $e^{-rt}$, respectively. Lemma~\ref{lem:common-envelope} bounds the dividend Stieltjes integrals and terminal functions of affine growth.

Let $\mathcal W_n,\mathcal Z_n$ be the two-variable $\omega_n$-scale functions. With baseline $r$,
\[
\mathcal W_n(x,z)=W^{(r)}(x-z)+\int_z^xW^{(r)}(x-u)(\omega_n(u)-r)\mathcal W_n(u,z)\,\dif u,
\]
and analogously for $\mathcal Z_n$. On each fixed compact interval, Volterra iteration bounds these functions uniformly in $n$. Subtracting the limiting Volterra equation and applying Gronwall's inequality bounds the uniform difference by a constant times
$\int_z^B|\omega_n(u)-\omega_y(u)|\,\dif u\le q/n$.
Hence the two functions converge locally uniformly to $\mathcal W^{(r,q)}$ and $\mathcal Z^{(r,q)}$, respectively, identified by uniqueness of the Volterra equations.

For the reflected formulas we differentiate in the first variable at $b>\max\{y,z\}$. If $z\ge y$, the integral term is zero and the derivative is the ordinary $r$-scale-function derivative. If $z<y$, then
\[
\mathcal W_{n,1}(b,z)=W^{(r)\prime}(b-z)
+\int_z^yW^{(r)\prime}(b-u)(\omega_n(u)-r)\mathcal W_n(u,z)\,\dif u.
\]
Since $\omega_n(b)=r$, the derivative is given by the displayed integral. The bound $b-u\ge b-y>0$ makes its derivative kernel bounded on compact intervals under the standing scale-function regularity assumptions. Dominated convergence gives
$\mathcal W_{n,1}(b,z)\to\mathcal W_1^{(r,q)}(b,z)$, and the same argument gives convergence of $\mathcal Z_{n,1}(b,z)$.

Applying \cite[Theorems 2.3--2.4]{li2018fluctuations} to $\omega_n$ and taking the value and derivative limits gives Lemma~\ref{lem:step-omega-reflected}. In bounded variation, the limiting resolvent includes the $W^{(r)}(0)$ boundary atom.

\end{proof}

\begin{proof}[Proof of Proposition~\ref{prop:fixed-triplet-value}]
\prooflocation{app:proof-fixed-triplet-value}
We can rewrite the third line in \eqref{eq:vc12bxy} for $x \in (\bs, b)$ as
\begin{align*}
& \mathbb E_x\left[e^{-A_{\tau_{s_2}^-(b)}^{b,y}}\widetilde g_{s_1,s_2}(Y_{\tau_{s_2}^-(b)}^b)
+\int_{[0,\tau_{s_2}^-(b)]}e^{-A_t^{b,y}}\,\dif D_t^b\right] \\
    &\qquad  = \mathbb{E}_x\Bigg[e^{-A_{\tau_{s_2}^-(b)}^{b,y}}\Big(\widetilde g_{s_1,s_2}(Y_{\tau_{s_2}^-(b)}^b)- g(s_2)\Big)\mathds{1}_{\{Y^b_{\tau^-_{\bs}(b)}<\bs\}}\Bigg]  \\
    & \qquad + g(s_2)\mathbb{E}_x\Big[e^{-A_{\tau_{s_2}^-(b)}^{b,y}}\Big] + \mathbb{E}_x\Bigg[\int_{[0,\tau_{s_2}^-(b)]}e^{-A_t^{b,y}}\,\dif D_t^b\Bigg].
\end{align*}
The compensation formula for the downward jumps of $Y^b$ and the resolvent \eqref{eq:resolvent-density} give
\begin{align*}
\mathbb E_x\left[
 e^{-A_{\tau_{s_2}^-(b)}^{b,y}}
 \left(
 \widetilde g_{s_1,s_2}(Y_{\tau_{s_2}^-(b)}^b)-g(s_2)
 \right)
 \mathds{1}_{\{Y_{\tau_{s_2}^-(b)}^b<s_2\}}
 \right]=
\int_{[s_2,b]}
\widetilde g_\nu(z)R_{s_2,b}^{(r,q)}(x,\dif z).
\end{align*}
The reflected exit identity gives
\[
\mathbb E_x\left[
 e^{-A_{\tau_{s_2}^-(b)}^{b,y}}
 \right]
=\mathcal Z^{(r,q)}(x,s_2)
-\frac{\mathcal Z_1^{(r,q)}(b,s_2)}
       {\mathcal W_1^{(r,q)}(b,s_2)}
 \mathcal W^{(r,q)}(x,s_2).
\]
Substituting these identities together with
\eqref{eq:dividend-reflected}, and collecting the coefficient of
$\mathcal W^{(r,q)}(x,s_2)$, yields \eqref{eq:Hgxc2}.
Differentiating \eqref{eq:Hgxc2} at $b-$ and using
\eqref{def:G-fun} gives \eqref{eq:barrier-smooth-auto}.
\end{proof}

\begin{proof}[Proof of Proposition~\ref{prop:fixed-triplet-harmonic}]
\prooflocation{app:proof-fixed-triplet-harmonic}
Fix a compact interval $I=[a,c]\subset\mathscr C_{s_1,s_2,b}$ and $x\in(a,c)$. Under the prescribed $(s_1,s_2,b)$-policy, neither voluntary liquidation nor dividend payment occurs before $T_I$. In particular, the controlled surplus
coincides with the uncontrolled process $X$ on $[0,T_I]$. Note that the occupation-time functional is additive, i.e., for
every stopping time $T$ and $t\ge0$, $A_{T+t}^{0,y} = A_T^{0,y} + rt + q\int_0^t \mathds{1}_{\{X_{T+s}<y\}}\,\dif s$.
Hence, applying the strong Markov property of $X$ at $T_I$ gives
\[
v_{s_1,s_2,b}(x;y)
=
\mathbb E_x\left[
e^{-A_{T_I}^{0,y}}
v_{s_1,s_2,b}(X_{T_I};y)
\right],
\]
which proves \eqref{eq:fixed-policy-harmonic}.

Now fix a point
$x\in\mathscr C_{s_1,s_2,b}\setminus\{y\}$. Choose a compact interval $I$ containing $x$ in its interior. Applying \eqref{eq:fixed-policy-harmonic} to $I$ and Dynkin's formula to
$e^{-A_{t\wedge T_I}^{0,y}}
v_{s_1,s_2,b}(X_{t\wedge T_I};y)$ yields
\[
\mathbb E_x\left[ \int_0^{t\wedge T_I} e^{-A_s^{0,y}} (\mathcal L-\omega_y(X_s)) v_{s_1,s_2,b}(X_s;y)\,\dif s
\right]=0.
\]
Dividing by $t$ gives
\[
0=\mathbb E_x\left[\frac1t\int_0^t\mathbf 1_{\{s<T_I\}}e^{-A_s^{0,y}}(\mathcal L-\omega_y(X_s))v_{s_1,s_2,b}(X_s;y)\,\dif s
\right].
\]
Since $x$ is an interior point, $T_I>0$ a.s., and since $(\mathcal L-\omega_y)v_{s_1,s_2,b}$ is continuous at $x$ and
locally bounded on $I$, the right-continuity of $X$ and the dominated
convergence theorem yield
\[
\lim_{t\downarrow0}\frac1t\mathbb E_x\left[\int_0^{t\wedge T_I}e^{-A_s^{0,y}}(\mathcal L-\omega_y(X_s))
v_{s_1,s_2,b}(X_s;y)\,\dif s
\right]=(\mathcal L-\omega_y(x))v_{s_1,s_2,b}(x;y).
\]
Hence $(\mathcal L-\omega_y(x))v_{s_1,s_2,b}(x;y)=0$, which proves \eqref{eq:fixed-policy-generator}.
\end{proof}

\begin{proof}[Proof of Proposition~\ref{prop:s1-identification}]
\prooflocation{App:D+1}
For $u<y$,
$\mathcal I^{(r,q)}(u-y)=e^{\Phi(r+q)(u-y)}$.
Hence the sign of $F_y'(u)$ is the sign of
\[
h(u):=g'(u)-\Phi(r+q)g(u).
\]
By concavity and monotonicity of $g$, $h$ is non-increasing on
$(\underline{k},y)$.  Indeed, for
$\underline{k}<u<v<y$, we have
\[
h(v)-h(u)
=\bigl(g'(v)-g'(u)\bigr)
-\Phi(r+q)\bigl(g(v)-g(u)\bigr)\le0.
\]
The conditions in \eqref{eq:threshold-regime-condition} give $h(\underline{k}+)=g'(\underline{k})>0$,
and $h(y-)=g'(y)-\Phi(r+q)g(y)<0$.
Thus $h$ has at least one zero in $(\underline{k},y)$.  If
$u<v$ were two zeros, monotonicity would imply $h=0$ on
$[u,v]$.  Then
\[
0=h(v)-h(u)
=\bigl(g'(v)-g'(u)\bigr)
-\Phi(r+q)\bigl(g(v)-g(u)\bigr),
\]
where both terms on the right are non-positive.  Hence, $g'(v)=g'(u)$ and $g(v)=g(u)$.  Concavity and monotonicity would make $g$ constant on $[u,v]$, which contradicts the fact that $h(z)=-\Phi(r+q)g(z)<0$ for $z\in[u,v]$. Therefore $h$ has a unique zero $s_1^*\in(\underline{k},y)$. Since $h>0$ on $(\underline{k},s_1^*)$ and $h<0$ on $(s_1^*,y)$, $F_y$ is strictly increasing on $(\underline{k},s_1^*)$ and strictly decreasing on
$(s_1^*,y]$, hence, we complete the proof.

\end{proof}

\begin{proof}[Proof of Lemma~\ref{lem:G-continuity}]
\prooflocation{App:D}
Write $\nu(\dif h)=\pi(h)\dif h$.  By Assumption
\textbf{(A.2)}, $\pi$ is continuous on $(0,\infty)$.  The function
$\widetilde g_{s_1,s_2}$ is continuous and locally Lipschitz on $(-\infty,s_2]$, and has at most affine growth. On a compact set $z\in[s_2,M]$, for $0<h<1$ and $h>z-s_2$, local Lipschitz continuity gives
\[
|\widetilde g_{s_1,s_2}(z-h)-g(s_2)|\le C h.
\]
For $h\ge1$ the bound is $C(1+h)$. Both bounds are integrable against $\nu$ by \textbf{(A.2)}. The payoff difference vanishes at $h=z-s_2$, so the integrand, including $\mathbf 1_{\{h>z-s_2\}}$, is continuous in $z$. Dominated convergence proves continuity of $\widetilde g_\nu$ on $[s_2,\infty)$, including its finite right limit at $s_2$.

For $x>m(s_2)$, the first-variable derivatives of the $\omega$-scale functions may be obtained directly from their Volterra representations.  When $s_2<y<x$,
\begin{equation*}
%\label{eq:W1-above-y}
\mathcal W_1^{(r,q)}(x,s_2) =W^{(r)\prime}(x-s_2)+q\int_{s_2}^{y} W^{(r)\prime}(x-u) W^{(r+q)}(u-s_2)\,\dif u,
\end{equation*}
whereas for $s_2\ge y$, $\mathcal W_1^{(r,q)}(x,s_2) =W^{(r)\prime}(x-s_2)$.
Assumption \textbf{(A.2)} gives the required continuity of
$W^{(r)\prime}$ on $(0,\infty)$, so
$\mathcal W_1^{(r,q)}(\cdot,s_2)$ is continuous on
$(m(s_2),\infty)$; the same argument applies to
$\mathcal Z_1^{(r,q)}(\cdot,s_2)$.

Finally, use \eqref{eq:Hg-derivative}.  On compact subsets of
$(m(s_2),\infty)$, the part of the integral with $z<y$ is handled by the dominated convergence theorem.  For $z\ge y$ one has
$\mathcal W_1^{(r,q)}(x,z)=W^{(r)\prime}(x-z)$, so the remaining term is a convolution of the locally integrable kernel
$W^{(r)\prime}$ with the continuous function
$\widetilde g_\nu$.  Hence $H_{\widetilde g,1}(\cdot,s_2)$ is continuous.  The continuity of $G$ then follows from \eqref{def:G-fun} and the positivity of
$\mathcal W_1^{(r,q)}(b,s_2)$.

\end{proof}

\begin{proof}[Proof of Lemma~\ref{lem:G-increment}]
\prooflocation{app:proof-G-increment}
Set
$
w(x):=v_{s_1,s_2,b+c}(x;y)-v_{s_1,s_2,b}(x;y).
$
For $s_2<x\le b$, both candidates are in the same continuation region, and \eqref{eq:vxy} gives
\begin{equation}\label{eq:w-below-b}
w(x)
=\bigl(G(b+c,s_2)-G(b,s_2)\bigr)
 \mathcal W^{(r,q)}(x,s_2).
\end{equation}
For $b<x\le b+c$, the candidate with barrier $b+c$ remains in the continuation region, whereas $v_b$ is in the dividend region.  Since $b>y$, the Omega clock is inactive on this interval and $(\mathcal L-r)w(x)=-J_{v_b}(x)$.
Moreover, for $b<z\le x\le b+c$ both spatial arguments lie above $y$, and hence $\mathcal W^{(r,q)}(x,z)=W^{(r)}(x-z)$.

Let
$M(x):=\int_b^x W^{(r)}(x-z)J_{v_b}(z)\,\dif z$
for $x>b$, and set $M(x)=0$ for $x\le b$.
The Laplace transform identity for $W^{(r)}$ yields
$(\mathcal L-r)M(x)=J_{v_b}(x)$ for $x>b$.
Combining this identity with \eqref{eq:w-below-b} and,
when $\sigma>0$, $w'(b+)=w'(b-)$ gives
\begin{equation}\label{eq:w-voc}
w(x)
=\bigl(G(b+c,s_2)-G(b,s_2)\bigr)
 \mathcal W^{(r,q)}(x,s_2)
-\int_b^x W^{(r)}(x-z)J_{v_b}(z)\,\dif z.
\end{equation}
At $x=b+c$, both candidate value functions have unit left derivative, so $w'((b+c)-)=0$.  Differentiating \eqref{eq:w-voc} from the left and writing the derivative of the convolution as a Lebesgue--Stieltjes integral with respect to $W^{(r)}$ yields \eqref{eq:G-increment}.
\end{proof}

\begin{proof}[Proof of Corollary~\ref{cor:G-eventual-decrease}]
\prooflocation{app:proof-G-eventual-decrease}
Fix $b\ge B$ and $c>0$. By \eqref{eq:J-uniform-tail}, $J_{v_b}(b+c-u)\le0$, for $u\in[0,c)$, since $b+c-u>b$. As $W^{(r)}(\dif u)$ is a nonnegative Stieltjes measure,
\[
\int_{[0,c)}
J_{v_b}(b+c-u)\,W^{(r)}(\dif u)\le0.
\]
Together with
$\mathcal W_1^{(r,q)}(b+c,s_2)>0$, \eqref{eq:G-increment} therefore yields $G(b+c,s_2)\le G(b,s_2)$.
Since $b\ge B$ and $c>0$ are arbitrary,
$G(\cdot,s_2)$ is non-increasing on $[B,\infty)$.
Under the additional condition, the integrand is strictly negative
on a set of positive $W^{(r)}$-measure. Hence
\[
\int_{[0,c)}
J_{v_b}(b+c-u)\,W^{(r)}(\dif u)<0,
\]
and \eqref{eq:G-increment} gives $G(b+c,s_2)<G(b,s_2)$. Thus $G(\cdot,s_2)$ is strictly decreasing on $[B,\infty)$.
\end{proof}

\begin{proof}[Proof of Lemma~\ref{lem:finite-no-liq}]
\prooflocation{app:proof-finite-no-liq}
According to \cite[Lemma 2.6]{rodosthenous2018beating},
\[
\frac{\mathcal I^{(r,q)\prime}(z)}
     {\mathcal I^{(r,q)}(z)}
\longrightarrow\Phi(r)>0,
\qquad z\to\infty.
\]
Hence, for all sufficiently large $z$, $\frac{\mathcal I^{(r,q)\prime}(z)}{\mathcal I^{(r,q)}(z)}
\ge \frac{\Phi(r)}{2}.$
Integration shows that $\mathcal I^{(r,q)}(z)\to\infty$
exponentially, and consequently $\mathcal I^{(r,q)\prime}(z)\longrightarrow\infty$. Therefore, $C_0(b)\longrightarrow0$, as $b\to\infty$. Let $M:=\sup_{b>y}C_0(b)>0$. Condition \eqref{eq:C0-boundary-gap} implies that $C_0$ is bounded
strictly below $M$ on a sufficiently small right-neighborhood of
$y$, while $C_0(b)\to0$ implies the same for all sufficiently large
$b$.  Since $C_0$ is continuous, its supremum is attained in the compact. Hence $\mathcal M_0$ is nonempty and bounded.
\end{proof}

\section{Verification of candidate policies}
\label{app:policy-verification}

\begin{proof}[Proof of Lemma~\ref{lem:obstacle-vstar}]
\prooflocation{app:proof-obstacle-vstar}
For $x<s_1^*$, Proposition \ref{prop:s1-identification} gives
$v^*(x;y)\ge g(x)$ whenever $g(x)>0$, while
$v^*(x;y)>0\ge g(x)$ whenever $g(x)\le0$. Equality holds on
$[s_1^*,s_2^*]$ by \eqref{eq:vxy*}, and \eqref{eq:upper-majorant-band} covers
$(s_2^*,b^*]$.
It remains to consider $x>b^*$. By the concavity of $g$ and
$x_g\le b^*$, we have $g(x)-g(b^*) \le g'(b^*-)(x-b^*) \le x-b^*$. Since \eqref{eq:upper-majorant-band} gives
$v^*(b^*;y)\ge g(b^*)$ and $v^*$ is affine with slope one on $[b^*,\infty)$, then
\[
v^*(x;y)
=
v^*(b^*;y)+x-b^*
\ge g(b^*)+x-b^*
\ge g(x).
\]
This proves \eqref{eq:derived-obstacle}.
\end{proof}

\begin{proof}[Proof of Lemma~\ref{lem:generator-vstar}]
\prooflocation{app:proof-generator-vstar}
On the continuation regions $(-\infty,s_1^*)$ and
$(s_2^*,b^*)$, Proposition
\ref{prop:fixed-triplet-harmonic} gives
$
(\mathcal L-\omega_y(x))v^*(x;y)=0
$
on every connected component away from $y$. For $x\in(s_1^*,s_2^*)$, since
$v^*(x-z;y)=g(s_1^*)
\frac{\mathcal I^{(r,q)}(x-z-y)}
     {\mathcal I^{(r,q)}(s_1^*-y)}$ for $z>x-s_1^*$, we have
$(\mathcal L-\omega_y(x))v^*(x;y)=R_{v^*}^{\rm stop}(x)$. Hence \eqref{eq:S2a} gives the required inequality on the interior of the liquidation region.
Finally, for $x>b^*>y$,
\[
(\mathcal L-\omega_y(x))v^*(x;y)
=
(\mathcal L-r)v^*(x;y)
=
J_{v^*}(x)\le0
\]
by \eqref{eq:S3}.  This proves the assertion.
\end{proof}

\begin{proof}[Proof of Theorem~\ref{thm:optimality_condition}]
\prooflocation{App:E+1}
Use the finite breakpoint set $\mathcal B=\{s_1^*,s_2^*,y,b^*\}\cup\mathcal B_g$, where $\mathcal B_g$ contains the finitely many second-derivative breakpoints of $g$ when $\sigma>0$ and is empty when $\sigma=0$. The scale-function regularity established after Proposition~\ref{prop:fixed-triplet-value}, together with the imposed fit, gives the local regularity of Theorem~\ref{prop:direct-verification}. In particular, the second distributional derivative is a locally finite signed measure in the Gaussian case. Every nonlocal term is finite by local Lipschitz continuity, affine growth, and the finite first moment.
By Lemma \ref{lem:obstacle-vstar}, $x_g\le b^*$ together with
\eqref{eq:upper-majorant-band} gives $v^*(x;y)\ge g(x),~ x\in\mathbb R$. By Lemma \ref{lem:dividend-vstar}, and, when $s_2^*>y$, the additional condition $g'(s_2^*)\ge1$, we have $v^*(x;y)-v^*(x-\eta;y)\ge\eta$, for $x\ge y$, $0\le\eta\le x-y$. Finally, Lemma \ref{lem:generator-vstar}, together with
\eqref{eq:S2a} and \eqref{eq:S3}, gives $(\mathcal L-\omega_y(x))v^*(x;y)\le0$ at every interior point away from $s_1^*,s_2^*,b^*$ and $y$.

It remains to check the boundary requirements in Theorem \ref{prop:direct-verification}.

Suppose first that $\sigma>0$. At $s_1^*$, since
$s_1^*<y$, the lower continuation branch is proportional to
$e^{\Phi(r+q)x}$, and the defining equation $g'(s_1^*)=\Phi(r+q)g(s_1^*)$ gives the smooth fit. The smooth fit at $s_2^*$ is imposed by \eqref{eq:s2_unbounded}. At $b^*$,
the definition of $G$ gives $v_x^*(b^*-;y)=1=v_x^*(b^*+;y)$.
At the regulatory threshold $y$, if $y$ lies in the interior of the
liquidation band, then $v^*=g$ on a neighborhood of $y$; if $y$ lies in
a continuation component, the scale-function representation is
continuously differentiable across $y$; and if $y=s_2^*$, the smooth-fit condition at $s_2^*$ applies directly. Hence
\eqref{eq:verification-local-time} holds at every boundary point.
At the additional points of $\mathcal B_g$, the two first derivatives agree because $g\in C^1$, so their local-time coefficients are zero as well.

When $\sigma=0$, define the boundary residual by
\[
B_f(a):=\delta_X-rf(a;y)+\int_{(0,\infty)}[f(a-z;y)-f(a;y)]\nu(\dif z),\qquad a\ge y.
\]
For either finite one-sided derivative, this is $(\mathcal L_{\pm}-r)f(a;y)+\delta_X(1-f_x(a\pm;y))$, with $\mathcal L_{\pm}$ defined using that same derivative. If $s_2^*>y$, the stopping-side limit of \eqref{eq:S2a} gives $(\mathcal L_- -r)v^*(s_2^*;y)\le0$, and $v_x^*(s_2^*-;y)=g'(s_2^*)\ge1$. Hence $B_{v^*}(s_2^*)\le0$. Since $g'(y)\ge g'(s_2^*)\ge1$, the limit of \eqref{eq:S2a} from above $y$, where the killing rate is $r$, similarly gives $B_{v^*}(y)\le0$. At $b^*$ the dividend-side derivative is one. The right limit of \eqref{eq:S3}, justified by the same integrable nonlocal bound, gives $B_{v^*}(b^*)\le0$. Thus \eqref{eq:verification-BV-boundary} holds at these boundary points.

If instead $s_2^*\le y$, the condition at $b^*$ is verified as above,
while \eqref{eq:BV-threshold-boundary} is precisely the required
bounded-variation boundary condition at $y$; when $s_2^*=y$, it also
covers the upper liquidation boundary.

Therefore, $v^*$ satisfies all the conditions of Theorem \ref{prop:direct-verification}, and hence $V(x;y)\le v^*(x;y)$ for $x\in\mathbb R$.

Conversely, $v^*$ is the performance of the admissible Markov policy
$(D^{b^*},\tau^*)$. Hence, $V(x;y)\ge V^{\Mmark}(x;y)\ge J(x;y;D^{b^*},\tau^*)
=v^*(x;y)$.
Combining the two inequalities yields $V(x;y)=V^{\Mmark}(x;y)=v^*(x;y)$,
and $(D^{b^*},\tau^*)$ is globally optimal.

\end{proof}

\begin{proof}[Proof of Theorem~\ref{thm:singleton}]
\prooflocation{App:E+2}
Take $\mathcal B=\{a,y,b_{\rm s}\}$. The regularity argument after Proposition~\ref{prop:fixed-triplet-value} applies with $s_1=s_2=a$.
In particular, the finite limit $\widetilde g_\nu(a+)$ and
$W^{(r)\prime}(0+)=2/\sigma^2$ imply finite one-sided first
derivatives at $a$ when $\sigma>0$; continuous fit in bounded
variation follows from \eqref{eq:singleton-cont-fit}.
For $x<a$, the maximizing property of $F_y$ gives
$v^{\rm s}(x;y)\ge g(x)$. Together with
\eqref{eq:singleton-majorant} and $x_g\le b_{\rm s}$, this gives $v^{\rm s}(x;y)\ge g(x),~ x\in\mathbb R$.
Lemma \ref{lem:singleton-dividend} gives the finite-dividend inequality required by Theorem \ref{prop:direct-verification}. Away from $a$, $y$, and $b_{\rm s}$, Proposition \ref{prop:fixed-triplet-harmonic},
applied with $s_1=s_2=a$ and $b=b_{\rm s}$, gives
\[
(\mathcal L-\omega_y(x))v^{\rm s}(x;y)=0,
\qquad
x<b_{\rm s},\quad x\ne a,y,
\]
while \eqref{eq:singleton-tail} gives the required generator inequality condition for $x>b_{\rm s}$.

It remains to verify the boundary conditions. Suppose first that
$\sigma>0$. Since $v^{\rm s}\ge g$ with $v^{\rm s}(a;y)=g(a)$, one necessarily has $v_x^{\rm s}(a+;y)\ge g'(a)$, then by \eqref{eq:singleton-derivative-condition}, we have $v_x^{\rm s}(a+;y) =v_x^{\rm s}(a-;y) =g'(a)$. By the $\omega$-scale function representation, $v^s$ is  continuously differentiable at $y$, and at $b_{\rm s}$ Proposition \ref{prop:fixed-triplet-value} gives
$v_x^{\rm s}(b_{\rm s}-;y) =1 =v_x^{\rm s}(b_{\rm s}+;y)$. Hence the condition (d) of Theorem
\ref{prop:direct-verification} is satisfied at every boundary point.

When $\sigma=0$, continuity at $a<y$ follows from \eqref{eq:singleton-cont-fit}. The boundary checks in \eqref{eq:verification-BV-boundary} concern $y$ and $b_{\rm s}$. At $b_{\rm s}$, we have $v_x^{\rm s}(b_{\rm s}-;y)=1$.
The right-side derivative at $b_{\rm s}$ is also one. Local Lipschitz continuity and the finite first moment make the nonlocal term continuous there. Taking the limit of \eqref{eq:singleton-tail} from above $b_{\rm s}$ therefore gives directly
\[
\delta_X-rv^{\rm s}(b_{\rm s};y)+\int_{(0,\infty)}[v^{\rm s}(b_{\rm s}-z;y)-v^{\rm s}(b_{\rm s};y)]\nu(\dif z)\le0,
\]
which is \eqref{eq:verification-BV-boundary} at $b_{\rm s}$.
At the regulatory threshold $y$,
\eqref{eq:singleton-BV-y} is precisely the required
bounded-variation boundary condition.
Thus all the conditions of Theorem
\ref{prop:direct-verification} are satisfied, and hence $V(x;y)\le v^{\rm s}(x;y)$, for $x\in\mathbb R$.

On the other hand, Proposition \ref{prop:fixed-triplet-value},
applied with $s_1=s_2=a$ and $b=b_{\rm s}$, gives $J(x;y;D^{b_{\rm s}},\tau_{\rm s})
=
v^{\rm s}(x;y)$.
Since this is an admissible Markov policy, $V(x;y)\ge V^{\Mmark}(x;y)
\ge J(x;y;D^{b_{\rm s}},\tau_{\rm s})
=
v^{\rm s}(x;y)$.
Therefore $V(x;y)=V^{\Mmark}(x;y)=v^{\rm s}(x;y)$, and $(D^{b_{\rm s}},\tau_{\rm s})$ is globally optimal.

\end{proof}

\begin{proof}[Proof of Lemma~\ref{lem:no-liq-obstacle}]
\prooflocation{app:proof-no-liq-obstacle}
For $x<b_0^*$ with $g(x)>0$,
\eqref{eq:Fbar-global} and $\Delta_0\ge0$ give
$
g(x)
\le
\overline F_y\mathcal I^{(r,q)}(x-y)
\le
C_0(b_0^*)\mathcal I^{(r,q)}(x-y)
=
v_{b_0^*}^0(x;y).
$
If $g(x)\le0$, the inequality follows from
$v_{b_0^*}^0(x;y)>0$. By continuity, it also holds at $x=b_0^*$.

For $x>b_0^*$, concavity of $g$ and $x_g\le b_0^*$ give
$
g(x)-g(b_0^*)
\le
g'(b_0^*-)(x-b_0^*)
\le
x-b_0^*.
$
Since $v_{b_0^*}^0$ is affine with slope one above $b_0^*$,
$
v_{b_0^*}^0(x;y)
=
v_{b_0^*}^0(b_0^*;y)+x-b_0^*
\ge g(x).
$
\end{proof}

\begin{proof}[Proof of Theorem~\ref{thm:no-liquidation}]
\prooflocation{App:E+3}
By Lemma~\ref{lem:no-liq-obstacle}, $v_{b_0^*}^0\ge g$ on $\mathbb R$, and \eqref{eq:no-liq-dividend} gives the finite-dividend inequality. Thus conditions (a)--(b) of Theorem~\ref{prop:direct-verification} hold. The scale-function representation is piecewise smooth, with the local integrability of its second derivatives in the Gaussian case explained after Proposition~\ref{prop:fixed-triplet-value}. It has at most affine growth after its affine extension above $b_0^*$, and its nonlocal terms are finite. Use $\mathcal B=\{y,b_0^*\}$.

The infinite-horizon policy-value argument preceding \eqref{eq:no-liq-fixed-policy} gives $\omega_y$-harmonicity below the barrier. Therefore $(\mathcal L-r-q)v_{b_0^*}^0=0$ on $(-\infty,y)$ and $(\mathcal L-r)v_{b_0^*}^0=0$ on $(y,b_0^*)$. Condition \eqref{eq:no-liq-tail} gives the generator inequality above $b_0^*$, so (c) holds.

At the dividend barrier, the normalization gives
\begin{equation}\label{eq:no-liq-barrier-slope}
(v_{b_0^*}^0)'(b_0^*-;y)=(v_{b_0^*}^0)'(b_0^*+;y)=1.
\end{equation}
If $\sigma>0$, $W^{(r)}(0)=0$, and
\[
\mathcal I^{(r,q)\prime}(0-)=\Phi(r+q),\qquad
\mathcal I^{(r,q)\prime}(0+)=\Phi(r+q)-qW^{(r)}(0)
\]
show that the derivatives at $y$ also agree. Thus every local-time coefficient in (d) is zero.

If $\sigma=0$, use the boundary residual
\[
B_f(a)=\delta_X-rf(a;y)+\int_{(0,\infty)}[f(a-z;y)-f(a;y)]\nu(\dif z).
\]
At $b_0^*$, take the limit from the left and use $(v_{b_0^*}^0)'(b_0^*-;y)=1$ to obtain $B_{v_{b_0^*}^0}(b_0^*)=0$. At $y$, take the  limit from the right, where the killing rate is $r$, to obtain
\[
B_{v_{b_0^*}^0}(y)=\delta_X\bigl(1-(v_{b_0^*}^0)'(y+;y)\bigr)\le0.
\]
The inequality follows from \eqref{eq:no-liq-gradient}, by taking its right limit at $y$. The integral limits are justified by local Lipschitz continuity and the finite first moment. Hence condition (d) holds.

Theorem~\ref{prop:direct-verification} now gives $V\le v_{b_0^*}^0$. Conversely, \eqref{eq:no-liq-fixed-policy} and the admissibility of the interior barrier policy give
\[
V(x;y)\ge V^{\Mmark}(x;y)\ge J(x;y;D^{b_0^*},\infty)=v_{b_0^*}^0(x;y).
\]
Consequently, $V=V^{\Mmark}=v_{b_0^*}^0$ and $(D^{b_0^*},\infty)$ is globally optimal. The common bound in Lemma~\ref{lem:common-envelope} handles all terminal and localization limits used here.

\end{proof}

\begin{proof}[Proof of Proposition~\ref{prop:q0-endpoint}]
\prooflocation{app:proof-q0-endpoint}
The finite-dividend inequality holds with equality on $[y,\infty)$, and $u_y$ is $r$-harmonic on $(-\infty,y)$. At $y$, $u_y'(y-)=u_y'(y+)=1$, so the diffusion local-time term vanishes. For $x>y$, let $J_{u_y}(x)=(\mathcal L-r)u_y(x)$. The finite first moment and $0<u_y'(z)\le1$ give
\[
J_{u_y}'(x)=-r+\int_{(0,\infty)}[u_y'(x-z)-1]\nu(\dif z)\le-r.
\]
The small-jump part is zero for $z<x-y$, which justifies differentiation. The identity $\psi(\Phi_r)=r$ gives
\[
J_{u_y}(y+)=-\frac{\sigma^2\Phi_r}{2}\le0.
\]
Thus the dividend-region generator is nonpositive. In bounded variation the same identity gives
\[
\delta_X-r u_y(y)+\int_{(0,\infty)}[u_y(y-z)-u_y(y)]\nu(\dif z)=0,
\]
which is the required boundary condition at $y$. All conditions of Theorem~\ref{prop:direct-verification} hold, so $V\le u_y$. The performance of the admissible policy $(D^y,\infty)$ gives the reverse inequality. Moreover, the interior barrier strategies belong to the auxiliary class and $v_b^0(x;y)\to u_y(x)$ as $b\downarrow y$, by \eqref{eq:no-liq-value-b} with $q=0$. Hence $V^{\Mmark}\ge\lim_{b\downarrow y}J(x;y;D^b,\infty)=u_y$; together with $V^{\Mmark}\le V=u_y$, this proves $V^{\Mmark}=V=u_y$.
\end{proof}

\section{Benchmark models}
\label{app:benchmarks}

\begin{proof}[Proof of Proposition~\ref{prop:G-tail-BM}]
\prooflocation{App:E}
For $b>m(s_2)=\max\{s_2,y\}$, the representations
\eqref{eq:BM-W}--\eqref{eq:BM-omega-W} imply
\begin{equation}
\mathcal W_1^{(r,q)}(b,s_2) = A_W e^{\Phi_r b} + B_W e^{R_r b},
\label{eq:BM-W1-exp}
\end{equation}
for constants $A_W=A_W(s_2)>0$ and $B_W=B_W(s_2)$.
Indeed, when $s_2<y$, the integral in
\eqref{eq:BM-omega-W} is over the fixed interval $[s_2,y]$
and therefore only changes the coefficients of
$e^{\Phi_r b}$ and $e^{R_r b}$; when $s_2\ge y$,
\eqref{eq:BM-W1-exp} follows directly from the classical
scale function.

The same argument applied to $\mathcal Z^{(r,q)}$, together
with \eqref{eq:BM-H}, gives
\begin{equation}
H_{\widetilde g,1}(b,s_2)
=
A_H e^{\Phi_r b}
+
B_H e^{R_r b}
\label{eq:BM-H1-exp}
\end{equation}
for constants $A_H$ and $B_H$.  Hence
\[
G(b,s_2)
=
\frac{
1-A_H e^{\Phi_r b}-B_H e^{R_r b}}
{
A_W e^{\Phi_r b}+B_W e^{R_r b}
}.
\]
Dividing numerator and denominator by $e^{\Phi_r b}$ yields
\begin{equation}
G(b,s_2)
=
-\frac{A_H}{A_W}
+
\frac{1}{A_W}e^{-\Phi_r b}
+
O\!\left(e^{(R_r-\Phi_r)b}\right),
\qquad b\to\infty.
\label{eq:BM-G-tail}
\end{equation}
Since $R_r<0$, the last term decays strictly faster than
$e^{-\Phi_r b}$.  Differentiating the exact exponential
representation, therefore, gives
\[
G_1(b,s_2)
=
-\frac{\Phi_r}{A_W}e^{-\Phi_r b}
+
o(e^{-\Phi_r b})
<0
\]
for all sufficiently large $b$.  This proves the result.

\end{proof}

\begin{proof}[Proof of Proposition~\ref{prop:G-tail-CL}]
\prooflocation{App:F}
First, make the change of variable $u=z-h$ in \eqref{eq:CL-overshoot}.  Since $h>z-s_2$ is equivalent to
$u<s_2$,
\begin{align}
\widetilde g_\nu(z) &= \lambda\alpha e^{-\alpha z} \int_{-\infty}^{s_2}
\left[ \widetilde g_{s_1,s_2}(u)-g(s_2) \right]e^{\alpha u}\,\dif u \nonumber\\
&=:C_\nu(s_1,s_2)e^{-\alpha z}.
\label{eq:CL-overshoot-exp}
\end{align}
The constant $C_\nu(s_1,s_2)$ is finite because $\widetilde g_{s_1,s_2}$ has at most affine growth.

For $b>m(s_2)$, the finite-exponential representation \eqref{eq:CL-W} and the step-$\omega$ convolution imply
\begin{equation}
\mathcal W_1^{(r,q)}(b,s_2) = A_W e^{\rho_1(r)b} + B_W e^{\rho_2(r)b}, \qquad A_W>0.
\label{eq:CL-W1-exp}
\end{equation}
Similarly,
$\mathcal Z_1^{(r,q)}(b,s_2)$ is a linear combination of $e^{\rho_1(r)b}$ and $e^{\rho_2(r)b}$.

Using \eqref{eq:Hg-derivative} and \eqref{eq:CL-overshoot-exp}, split the integral at $y$.
For $z\ge y$, $\mathcal W_1^{(r,q)}(b,z)=W^{(r)\prime}(b-z),$ whereas the part with $z<y$ is integrated over a fixed
interval.  Since $W^{(r)\prime}$ is itself a linear
combination of $e^{\rho_1(r)(b-z)}$ and
$e^{\rho_2(r)(b-z)}$, convolution with $e^{-\alpha z}$
gives
\begin{equation}
H_{\widetilde g,1}(b,s_2)
=
A_H e^{\rho_1(r)b}
+
B_H e^{\rho_2(r)b}
+
C_H e^{-\alpha b}
\label{eq:CL-H1-exp}
\end{equation}
for constants $A_H,B_H,C_H$ depending on $(s_1,s_2)$ and
the model parameters.

Substituting \eqref{eq:CL-W1-exp} and
\eqref{eq:CL-H1-exp} into \eqref{def:G-fun}, and recalling
that $\rho_1(r)=\Phi_r>0>\rho_2(r)$, gives
\begin{align}
G(b,s_2)
&=
-\frac{A_H}{A_W}
+\frac{1}{A_W}e^{-\Phi_r b}
+O\!\left(e^{(\rho_2(r)-\Phi_r)b}\right)
+O\!\left(e^{-(\alpha+\Phi_r)b}\right).
\label{eq:CL-G-tail}
\end{align}
Both remainder terms decay strictly faster than
$e^{-\Phi_r b}$.  Hence differentiation of the exact
finite-exponential representation yields
\[
G_1(b,s_2)
=
-\frac{\Phi_r}{A_W}e^{-\Phi_r b}
+
o(e^{-\Phi_r b})
<0
\]
for all sufficiently large $b$.  Thus $G(\cdot,s_2)$ is
eventually strictly decreasing.

\end{proof}

\begin{proof}[Proof of Corollary~\ref{cor:benchmark-regular}]
\prooflocation{app:proof-benchmark-regular}
Use the notation of Proposition~\ref{prop:candidate-regular}.
For the band and singleton candidates, the exponential
representations above make $G(\cdot,\ell)$ differentiable
at its interior maximizer $b$.
Equation~\eqref{eq:gradient-G-identity} therefore gives
\[
v_{xx}(b-;y)
=-\mathcal W_1^{(r,q)}(b,\ell)G_1(b,\ell)=0.
\]
For the no-liquidation candidate, interior maximization
of $C_0$ gives
$\mathcal I^{(r,q)\prime\prime}(b-y)=0$,
and hence again $v_{xx}(b-;y)=0$.

On $(m,b)$, the harmonic equation reduces to
\[
\begin{aligned}
\frac{\sigma^2}{2}v_{xx}+cv_x-rv&=0
&&\text{(Brownian model)},\\
cv_{xx}+(c\alpha-\lambda-r)v_x-\alpha rv&=0
&&\text{(Cram\'er--Lundberg model)}.
\end{aligned}
\]
In either case, let $p>0>n$ denote the roots of the
characteristic polynomial.
The boundary conditions $v_x(b-;y)=1$ and
$v_{xx}(b-;y)=0$ give
\begin{equation}\label{eq:benchmark-strict-slope}
v_x(x;y)
=\frac{-n e^{p(x-b)}+p e^{n(x-b)}}{p-n},
\qquad m<x<b.
\end{equation}
The right-hand side is strictly decreasing on $(m,b)$
and tends to one as $x\uparrow b$.
Thus $v_x(x;y)>1$ throughout this interval.
When $y\notin\mathscr S_v$, we have $m=y$,
and the same formula gives $v_x(y+;y)>1$.

Finally, $v\ge g$ and $v_x(b-;y)=1>g'(b)$ imply
$v(b;y)>g(b)$.
Proposition~\ref{prop:candidate-regular} now applies.
\end{proof}

\end{appendices}

%%%%%%%%%%%%%%%%%%%%%%%%%%%%%%%%%%%%%%%%%%%%%%%%%%%%%%%%%%%%%%%%

\bibliographystyle{apalike}
\bibliography{Reference_final}

\end{document}